\documentclass[oneside]{amsart}
\usepackage{amsmath,amssymb,color,amsthm,graphicx,cite}
\usepackage{color,bm,multicol}
\graphicspath{{./img_pdf/}}
\usepackage{enumerate} 
\usepackage[backref]{hyperref}
\RequirePackage{luatex85}
\usepackage[all]{xy}
\newtheorem{theorem}{Theorem}[section]
\newtheorem{proposition}[theorem]{Proposition}
\newtheorem{lemma}[theorem]{Lemma}
\newtheorem{corollary}[theorem]{Corollary}
\newtheorem{claim}{Claim}

\newtheorem{example}{Example}
\theoremstyle{definition}
\newtheorem{definition}{Definition}
\newtheorem{main}{Theorem}

\def\F{\mathcal{F} }
\def\G{\mathcal{G} }
\def\Z{\mathbb{Z} }
\def\R{\mathbb{R} }
\def\A{\mathbb{A} }
\def\T{\mathbb{T} }

\def\nbd{neighborhood }
\def\nbds{neighborhoods }
\def\SF{\mathop{\mathrm{Sing}}(\mathcal{F})}

\author{Tomoo Yokoyama}
\date{\today}
\address{Department of Mathematics, Faculty of Science, Saitama University, Shimo-Okubo 255, Sakura-ku, Saitama-shi, 338-8570 Japan\\}
\email{tyokoyama@rimath.saitama-u.ac.jp}
\thanks{The author was partially supported by JSPS Grant Number 24K06733}
\subjclass[2010]{}

\makeatletter
\@namedef{subjclassname@2020}{%
  \textup{2020} Mathematics Subject Classification}
\makeatother

\title[Structural stability and transitions of incompressible line fields]{Structural stability and generic transitions of ``incompressible'' line fields on surfaces}

\newcommand{\B}{\Box}
\usepackage{multirow}

\begin{document}
\maketitle

\begin{abstract}
Various line fields naturally arise on surfaces in both physical and biological contexts, and generic singularities frequently appear in the form of 1-prong (thorn-like) and 3-prong (tripod-like) configurations, which can be modeled by partial differential equations with specific parameter values. However, it remains open under which topologies such line fields are structurally stable and form an open dense subset. In this paper, we propose a new topological framework for describing line fields and their evaluations on surfaces that is suitable from both theoretical and applied perspectives. Specifically, we demonstrate that, under a topology defined by a ``cone'' structure, line fields with 1-prong and 3-prong singularities are generic when an ``incompressibility condition''  holds. We also introduce representations of complete invariants for generic line fields and their generic transitions. These representations enable the evolution of ``incompressible'' line fields---such as those observed in active nematics---to be encoded as walks on transition graphs, providing a combinatorial framework for their analysis.
\end{abstract}

\section{Introduction}\label{intro}

Singular foliations on surfaces appear as a model of objects in nature (e.g. fingerprints\cite{gu2004combination,kass1987analyzing,PENROSE1979topology,wang2007fingerprint}, nematic liquid crystals\cite{decamp2015orientational,de1993physics}, the pinwheel structure of the visual cortex\cite{petitot2008neurogeometrie}, and second-order tensor fields\cite{Delmarcelle1995fields,hesselink1995topology,Tricoche2006topology}).  
Generic singularities of singular foliations on surfaces are studied in various settings (see ~\cite{boscain2016generic} for details). 

\begin{figure}[t]
\begin{center}
\includegraphics[scale=0.375]{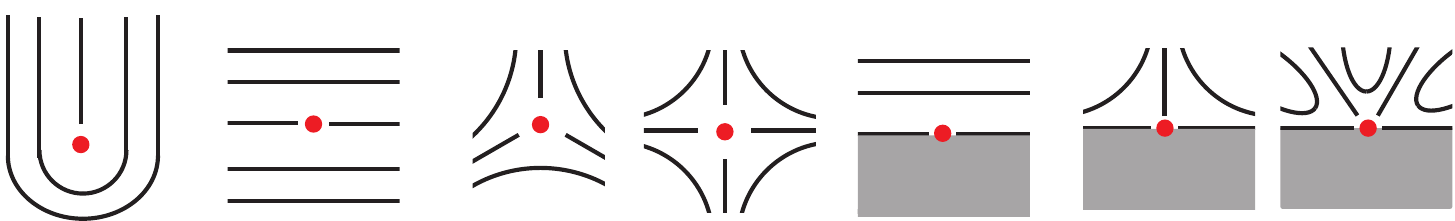}
\end{center} 
\caption{$1$-prong (thorn), $2$-prong (fake prong), $3$-prong (tripod), and $4$-prong outside of the boundary, and $2$-prong (fake prong), $3$-prong, and $4$-prong on the boundary.}
\label{prong}
\end{figure}

\begin{figure}[t]
\begin{center}
\includegraphics[scale=0.375]{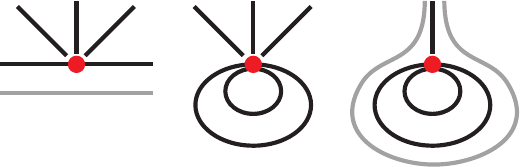}
\end{center} 
\caption{A sunset and two kinds of apples.}
\label{prong_02}
\end{figure}

As will be shown below, it can be seen that ``generic'' singularities and the structural stability of singular foliations on surfaces vary by topologies and class of singular foliations. 
In fact, it is pointed out that two types of singularities, thorns and tripods (also called stars) as in Figure~\ref{prong}, are usually observed in nature as shown in \cite[Figure~1]{boscain2016generic} (e.g. fingerprints, nematic liquid crystals, pinwheel structure of the orientation columns of the visual cortex), in the introduction of the paper. 
Sotomayor and Gutierrez showed that two types of thorns (called Lemons and Monstars \cite{berry1977umbilic,boscain2016generic}) and tripods are the structurally stable singularities of lines of principal curvature with respect to the Whitney $C^3$-topology of immersions of a surface in $\R^3$ in \cite[Theorem~2.4]{sotomayor1982structurally}. 
%
Bronshteyn and Nikolaev demonstrated that the structurally stable singularities of nonorientable singular foliations in the sense of Bronstein and Nikolaev are five types of singularities (i.e. tripods, sun-sets, thorns, and two kinds of apples as in Figure~\ref{prong_02})\cite[p.33 Theorem~1.4.3]{nikolaev2001foliations}(cf. \cite[Theorem~1]{Bronshteyn1998stable}) with respect to the Whitney $C^q$-topology~($r \geq 2$). 
In particular, the tripods and thorns, which are pronged singularities, are possible singularities in incompressible settings. 
Boscain, Sacchelli, and Sigalotti showed that singular foliations generated by generic proto-line-fields (i.e. the set leaves tangent to bisectors of a pair of vector fields), which are also called bisector line fields, are locally structurally stable and the singularities are lemons, monstars, and stars with respect to the Whitney $C^1$-topology\cite[Theorem~7]{boscain2016generic}.



Singular foliations with pronged singular points on compact surfaces are studied~\cite{fathi1979feuilletages02,levitt1982pantalons,levitt1987decomposition,lopez2007foliation,mendes1986interval}. 
For instance, Fathi classified measured foliations, which are some kinds of singular foliations, up to Whitehead equivalence~\cite{fathi1979feuilletages}. 
Moreover, Hounie classified the minimal sets of singular foliations on compact surfaces as well as vector fields  \cite[Theorem~2]{hounie1981minimal}. 
%
%
Rosenberg constructs counter-examples to the Poincar{\'e}-Bendixson theorem for singular foliations on the disk~\cite{rosenberg1983labyrinths}. 
On the other hand, the Poincar{\'e}-Bendixson theorem holds for singular foliations on compact surfaces with non-fake pronged singular points if the restriction to the complement of the singular point set is $C^2$ as a line field\cite{levitt1987differentiability}.

The singularities of a vector field and the Euler characteristic of the underlying manifold are related by the Poincar{\'e}-Hopf theorem. 
The result is generalized for one-dimensional singular foliations with finitely many singularities on surfaces \cite[p.113 Thoerem~2.2]{Hopf1989diff}, which is used in the physics literature  
 as mentioned above.


The structural stability and open denseness of Morse-Smale flows on orientable closed surfaces also hold for 
Morse-Smale foliations on orientable closed surfaces 
 with respect to the Whitney $C^q$-topology\cite[Theorem~2]{Bronshteyn1998stable}\cite[Theorem~2.2]{bronstein1997peixoto}.  

\subsection{Absence of a natural topological setting in the formulation of physical phenomena via line fields}

As mentioned above, it can be seen that ``generic'' singularities and structural stability of singular foliations on surfaces vary by topologies and class of singular foliations. 
Thus a suitable topology on the set of singular foliations, as partitions, is required to describe time evaluations of physical phenomena (e.g. nematic liquid crystal, Voronoi-like behaviors) and topological operations (e.g. Whitehead moves, bifurcations of prongs). 
On the other hand, as pointed out in \cite{boscain2016generic}, there is no natural topology on the set of one-dimensional singular foliations (i.e. line fields).
Moreover, from a data analysis point of view, one needs to analyze singular foliations as partitions of spaces directly, because such singular foliations are not equipped with any kind of structure (e.g. norms of line fields and smoothness) in general.  
However, there is no known class of ``incompressible'' line fields on surfaces equipped with a topology in which ``generic'' singularities are tripods and thorns and which can describe time evaluations of physical phenomena modeled by line fields. 
To describe such phenomena, we introduce ``cone'' topology, which is analogous to the topology of vector fields. 

\subsection{Statements of main results}

As the space of vector fields with $C^0$-topology is not structurally stable in general, we have the following absence of stability. 

\begin{main}\label{non-existence_structural_stability}
For any $q \in \Z_{>0}$ and any $p \in \R_{\geq 1} \sqcup \{ \infty \}$, any singular foliation tangent to the boundary on a compact surface $M$ is not structurally stable in $\mathcal{F}^q(M)$ with respect to the topology $\mathcal{O}^p_{w,\F^q(M)}$ (see Definition~\ref{def:topology} for the definition), where the symbol $\sqcup$ denotes a disjoint union. 
\end{main}

On the other hand, as the space of Hamiltonian vector fields (see Definition~\ref{def:Ham_vf}  in Appendix below) on a compact surface has an open dense subspace that consists of structurally stable Hamiltonian vector fields, the openness, density, and structural stability hold for ``Hamiltonian'' singular foliations, which are called levelable foliations (see Definition~\ref{def:lebelable} for the definition). 
Notice the fake multi-saddles in Hamiltonian vector fields are inessential from a bifurcation's point of view under the non-existence of annihilation and creation. 
Similarly, our density holds under the absence of ``inessential'' structures (i.e. fake prongs). 

To state it more precisely, denote by $\mathcal{L}^*(M)$ the set of levelable foliations with finitely many pronged singular points and without fake prongs (see \S~\ref{sec:fake} for the definition) on a surface $M$ and by $\mathcal{L}^{0}(M) \subset \mathcal{L}^*(M)$ the subspace of levelable foliations whose prong connection diagrams (see \S~\ref{sec:separatrix} for the definition) are semi-self-connected (see \S~\ref{sec:sscpc} for the definition).  
We have the following density under the non-existence of fake prongs. 

\begin{main}\label{th:stability03-}
For any surface $M$ and any $p \in \R_{\geq 1} \sqcup \{ \infty \}$, 
the subspace $\mathcal{L}^{0}(M)$ is dense in $\mathcal{L}^*(M)$ with respect to the relative topology of the topology $\mathcal{O}^p_{w,W,\mathcal{L}(M)}$ (see Definition~\ref{dec:levelable_top} for the definition). 
\end{main}

Though the space of Hamiltonian vector fields with $C^0$-topology is not structurally stable in general, the structural stability holds under the non-existence of annihilation and creations \cite[Theorem~5.2.]{yokoyama2025combinatorial}. 
Similarly, we have the open density and the structural stability for levelable foliations under the non-existence of annihilation and creation. 
To state it more precisely, denote by $\mathcal{L}^*_{i}(M) \subset \mathcal{L}^*(M)$ the subspace of levelable foliations whose sums of indices of pronged singular points with positive indices are $i$. 
We have the following open density and structural stability under the non-existence of annihilation and creation. 

\begin{main}\label{th:stability03}
For any compact surface $M$, any $p \in \R_{\geq 1} \sqcup \{ \infty \}$, and any $i \in \Z_{>1}$,  the subspace $\mathcal{L}^0(M) \cap \mathcal{L}^*_{i}(M)$ is open dense in $\mathcal{L}^*_{i}(M)$ and consists of structurally stable ones in $\mathcal{L}^*_{i}(M)$ with respect to the relative topology of the topology $\mathcal{O}^p_{w,W,\mathcal{L}(M)}$. 
\end{main}

%

Denote by $\mathcal{L}^1_{i}(M)$ the subspace of levelable foliation in $\mathcal{L}^*_{i}(M) \setminus \mathcal{L}^0(M)$ whose prong connection diagram is of coheight-one (see Definition~\ref{dec:levelable_top} for the definition of coheight). 
We demonstrate the existence of ``generic'' intermediate structures as follows. 

\begin{main}\label{th:stability04}
The following statements hold for for any compact surface $M$ and any positive integer $i \geq 2$: 
\\
{\rm(1)} The subspace $\mathcal{L}^1_{i}(M)$ is an open dense subset of $\mathcal{L}^*_{i} \setminus \mathcal{L}^0(M)$. 
\\
{\rm(2)} Every levelable foliation in $\mathcal{L}^1_{i}(M)$ is structurally stable in $\mathcal{L}^*_{i} \setminus \mathcal{L}^0(M)$. 
\end{main}

Note that the previous results can imply a conversion from ``generic'' time evaluations of a levelable foliation to walks on a graph. 
As the ``generic'' intermediate structures of Hamiltonian vector fields on a compact surface are characterized using the multi-saddle connection diagrams, ``generic'' intermediate structures of levelable foliations can be characterized using the prong connection diagrams. 
We call any connected component of the prong connection diagram a prong connection.
Moreover, we characterize ``generic'' structures and ``generic'' intermediate structures as follows. 

\begin{main}\label{main:char_0}
The following statements hold for any levelable foliation $\F$ on a compact surface $M$ with finitely many pronged singular points: 
\\
{\rm(1)} $\F \in \mathcal{L}^0(M)$ if and only if the prong connection diagram $D(\F)$ of $M$ consists of coheight-zero prong connections. 
\\
{\rm(2)} A prong connection in $D(\F)$ is coheight-zero if and only if it is contained in Figure~\ref{self_connected_prong} as a topological graph.
\end{main}

\begin{main}\label{main:char}
The following statements hold for any levelable foliation $\F$ on a compact surface $M$ with finitely many pronged singular points: 
\\
{\rm(1)} $\F \in \mathcal{L}^1(M)$ if and only if the prong connection diagram $D(\F)$ of $M$ consists of coheight-zero prong connections and a coheight-one prong connection. 
\\
{\rm(2)} A prong connection in $D(\F)$ is coheight-one if and only if it is contained in Figures~\ref{b_prong}, \ref{c_prong}, \ref{h_unstable_prong_connection}, and \ref{h_unstable_prong_transition} as a topological graph. 
\end{main}

In addition, we have the following complete invariants for generic (i.e. coheight-zero) levelable foliations and generic intermediate (i.e. coheight-one) foliations between them on a sphere and a disk. 

\begin{main}\label{th:cot_02}
For any $p \in \R_{\geq 1} \sqcup \{ \infty \}$, the COT representation (see \S~\ref{sec:cot_zero}--\ref{sec:cot_one}
 for the definition)
for line fields of $\mathcal{L}^0(\mathbb{S}^2) \sqcup \mathcal{L}^1(\mathbb{S}^2)$ with root components is a finite complete invariant. 
\end{main}

\subsubsection{Contents}

The present paper consists of seventeen sections.
We define a topology on the set of line fields corresponding to the $C^0$ topology of the set of flows and demonstrate the absence of structurally stable line fields in \S~2--7, and discuss ``Hamiltonian'' (i.e. levelable) foliations for line fields in \S~8--15. 
In fact, in the next section, we provide preliminaries by recalling fundamental concepts.
In \S 3, we introduce topologies on the space of one-dimensional foliations to describe ``generic'' singularities and structural stability of singular foliations on surfaces. 
In \S 4, we state the fundamental properties (e.g. persistence of non-existence/existence of singularities, invariance of closed transverse subsets). 
In \S 5, the properties of the line field (i.e. one-dimensional singular foliations) on surfaces are discussed. 
In \S 6, the invariance of line fields is described. 
In \S 7, we demonstrate Theorem~\ref{non-existence_structural_stability}, the absence of structural stability for general line fields on surfaces. 
In \S 8, an analogous concept of ``incompressible'' for line fields is introduced and studied.
In \S 9, ``Hamiltonian'' (i.e. levelable) foliations for line fields are introduced, and such properties are discussed. 
In \S 10, we define the topologies of the set of levelable foliations. 
In \S 11, we define the coheight of levelable foliations, and the structural stability of levelable foliations on surfaces is demonstrated. 
In particular, Theorem~\ref{main:char_0} is shown.
In \S 12, the properties with respect to the topology $\mathcal{O}^p_{w,W,\mathcal{L}(M)}$ is discussed. In particular, Theorem~\ref{th:stability03-} and Theorem~\ref{th:stability03} are demonstrated. 
In \S 13, we characterize coheight-one elements in the space of levelable foliations on surfaces and show Theorem~\ref{th:stability04} and Theorem~\ref{main:char}. 
In \S 14--15, we introduce COT representations, which are complete invariants, for generic (i.e. coheight-zero) levelable foliations and generic intermediate foliations between them on a sphere and a disk. 
In particular, we show Theorem~\ref{th:cot_02}. 
In \S 16, various examples of foliations are constructed. 
In particular, some examples illustrate the necessity and non-liftability to flows. 
In the final section, we discuss line fields, such as Voronoi partitions, and remark on line fields on non-compact surfaces.

\section{Preliminaries}\label{sec:prel}

\subsection{Notion of topology and combinatorics}

By a {\bf surface}, we mean a connected two-dimensional Riemannian manifold possibly with boundary.
The {\bf double} of a manifold $M$ with boundary is the resulting manifold $M\times \{0,1\}/\mathop{\sim} $, where $(x,0)\sim (x,1)$ if $x\in \partial M$, see Figure \ref{pic01}. 
%

For any subsets $A, B, C$, the union $A \cup B$ is denoted by $A \sqcup_C B$ when $A \cap B = C$.

\subsubsection{Curves, arcs, and loops}
A {\bf curve} is a continuous mapping $C: I \to X$ to a topological space $X$, where $I$ is a non-degenerate connected subset of a circle $\mathbb{S}^1$.
A curve is {\bf simple} if it is injective.
We also denote by $C$ the image of a curve $C$.
Denote by $\partial C := C(\partial I)$ the boundary of a curve $C$ if $C$ can be extended into a continuous map whose domain is $I \cup \partial I$, where $\partial I$ is the boundary of $I \subset \mathbb{S}^1$. 
Put $\mathop{\mathrm{int}} C := C \setminus \partial C$ if $\partial C$ is defined. 
A simple curve is a {\bf simple closed curve} if its domain is $\mathbb{S}^1$ (i.e. $I = \mathbb{S}^1$).
A simple closed curve is also called a {\bf loop}. 
An {\bf arc} is a simple curve with a non-degenerate interval domain. 

A $C^{q-1}$ curve $C \colon I \to M$ on a manifold $M$ is of {\bf piecewise $\bm{C^q}$} if there are finite subset $J \subset I$ such that the restriction $C \vert_{I - J }$ is of $C^q$ and that, for any point $t_0 \in J$, there the one-sieded limits 
\[
\lim_{s \to t_0+} \dfrac{C^{(q-1)}(s) - C^{(q-1)}(t_0)}{s - t_0}
\]
and 
\[
\lim_{s \to t_0-} \dfrac{C^{(q-1)}(s) - C^{(q-1)}(t_0)}{s - t_0}
\]
exist, where $C^{(q-1)}$ is the $(q-1)$-th derivative of $C$. 

%

\subsubsection{Hausdorff distance}

Let $(X,d)$ be a metric space. 
For any point $x$ on a metric space $X$ and any number $r \geq 0$, put 
\[
\bm{B_r(x)} = \bm{B(x,r)} := \{ y \in M \mid d(x,y) \leq r \}
\]
called the closed {\bf $\bm{r}$-ball} centered at $x$ in $X$. 
For any subset $A \subseteq X$ and a non-negative number $r \geq 0$, define the closed  {\bf $\bm{r}$-neighborhood} $B_r(A)$ as follows: 
\[
\bm{B_r(A)} := \bigcup_{a \in A} \{ x \in X \mid d(a,x) \leq r \} = \bigcup_{a \in A} B_r(a) 
\]
For any nonempty subset $A \subseteq X$ and any point $x \in X$, define $d(x, A)$ as follows: 
\[
d(x, A) := \inf_{a \in A} d(x,a)
\]
Write $d(x, \emptyset) := \infty$ for any point $x \in X$. 

The {\bf Hausdorff distance} $d_H(A,B)$ between nonempty subsets $A, B \subseteq X$ is defined as follows: 
\[
\bm{d_H(A,B)} := \inf \{ r \geq 0 \mid A \subseteq B_r(B), B \subseteq B_r(A) \}
\]

Notice that the Hausdorff distance is a pseudo-metric but need not be a metric.

\subsubsection{Distances and volumes on a Riemannian manifold}

The {\bf distance} $\bm{d_g}$ between two points on a Riemannian manifold $(M,g)$ is inferior to the length of the piecewise $C^1$-curves between them. 
Notice that $d_g(x,y) = \infty$ for any points $x,y \in M$ which are contained in different connected components. 
Write the Riemannian volume form $d \mathrm{vol}_g := {\star }(1)$, where $\star$ is the Hodge star. 
 Denote by $\mathrm{vol}_g(U) := \int_U d \mathrm{vol}_g$ the volume for any Borel sets $U \subseteq M$.

\subsubsection{Diameter}

Let $(X,d)$ be a metric space. 
The {\bf diameter} $\bm{{\mathrm{diam}(A)}}$ of a subset $A \subseteq X$ is defined by the superior $\sup_{x,y \in A}d(x,y)$ of distances of pair of two points in $A$. 
Notice that $\mathrm{diam}(\emptyset) = 0$.

\subsection{Singular foliations on manifolds}


A standard reference of $C^q$ (regular) foliations is \cite{HH1986A,HH1987B,CC2003I,CC2003II}. 
We refer to Stefan \cite{Stefan1974} and Sussmann \cite{Sussmann1973} for the definition of singular foliations on a manifold $M$ without boundary. 
Indeed, we have the following definitions of leaves and singular foliations. 

\begin{definition}\label{def:foliation02}
For any $k \in \Z_{\geq 0}$ and any $q \in \Z_{\geq 0} \sqcup \{ \infty \}$, a subset $L$ of a paracompact manifold $M$ is a \textbf{$\bm{C^q}$ $\bm{k}$-leaf} if $L$ satisfies the following two conditions: 
\\
{\rm(1)} $L$ is a connected immersed $C^q$ submanifold {\rm(i.e.} the image of a $C^q$ injective immersion of a connected manifold{\rm)}. 
\\
{\rm(2)} For any locally connected topological space $X$ and any continuous mapping $h \colon X \to M$ with $h(X) \subset L$, the restriction $h \colon X \to L$ is continuous. 
\end{definition}

For any $k$-leaf $L$, write $\bm{\dim L} := k$, and the integer $k$ called the {\bf dimension} of $L$.

\begin{definition}\label{def:foliation}
For any $q \in \Z_{\geq 0} \sqcup \{ \infty \}$, a partition $\F$ into $C^q$ leaves is a $C^q$ {\bf singular foliation} $\F$ on a $n$-manifold $M$ without boundary if for any $x \in M$ there are a non-negative integer $k$ and a local $C^q$ chart $\varphi \colon V \to U_\tau \times U_\pitchfork \subseteq \R^k \times \R^{n-k}$ on $M$ around $x$ with the following three conditions: 
\\
{\rm(1)} 
$U_\tau$ is an open neighborhood of $0$ in $\R^k$ and $U_\pitchfork$ is an open neighborhood of $0$ in $\R^{n-k}$. 
\\
{\rm(2)} $\varphi(x) = (0,0)$. 
\\
{\rm(3)} For any leaf $L \in \F$ with $L \cap V \neq \emptyset$, we have $\varphi(L \cap V) = U_\tau \times l_{V \cap L}$, where $l_{V \cap L} := \{ u_\pitchfork \in U_\pitchfork \mid \varphi^{-1}(0, u_\pitchfork) \in L \}$. 
\end{definition}
The above local chart $\varphi \colon V \to U_\tau \times U_\pitchfork$ is called a {\bf ($\bm{\F}$-)fibered chart} at $x$ with respect to $\mathcal F$. 
For any $x_\tau \in U_\tau$, the subset $\{ x_\tau \} \times U_\pitchfork$ is called a {\bf plaque}. 
\begin{definition}
A fibered chart $\varphi \colon V \to U_\tau \times U_\pitchfork$ is {\bf trivial} if $U_\tau = (-1,1)^k$ and $U_\pitchfork = (-1,1)^{n-k}$.  
\end{definition}

Note that fibered charts are not well defined if $M$ has the nonempty boundary. 
\begin{definition}
A partition $\mathcal F$ on a manifold $M$ with boundary $\partial M$ is a $C^q$ {\bf singular foliation} if the induced partition $\hat{\mathcal F}$ of $\mathcal F$ on the double of $M$ is a $C^q$ singular foliation. 
\end{definition}

Define the dimension of each leaf  $L\in \mathcal F$ through a boundary point as the dimension of the lifted leaf of  $L$ in $\hat{\mathcal F}$. 
%
The integer $\dim \F := \max \{ \dim L \mid L \in \mathcal F\}$ is the {\bf dimension} of $\mathcal{F}$, and $\dim M -  \dim \F$  the {\bf codimension} of $\mathcal{F}$. 
Denote by $\bm{\mathcal{F}(x)}$ the leaf of a foliation $\mathcal{F}$ through a point $x\in M$.

\subsubsection{Types of leaves}
A leaf $L$ is said to be {\bf regular} if $\dim L =\dim \F$ and {\bf singular} if $\dim L <\dim \F$. 
The union of singular leaves is the {\bf singular set} and is denoted by $\mathop{\mathrm{Sing}}(\F)$. 

\subsubsection{Regular foliations}
A foliation $\mathcal F$ is a {\bf regular foliation}  if $\mathop{\mathrm{Sing}} (\mathcal F ) =\emptyset $.
By the existence of fibered charts, the union of regular leaves is open, so the union of singular leaves is closed. 
%


%



\subsubsection{Leaf arcs for singular foliations}

A non-degenerate curve is a {\bf leaf arc} if it is contained in a leaf.  
A simple arc is a {\bf simple leaf arc} if it is contained in a leaf.  
A simple curve is a {\bf simple leaf curve} if it is either a simple leaf arc or a loop contained in a leaf.  

\subsubsection{Local differentiability at a point}

To local differentiability, we define induced foliations as follows. 

\begin{definition}
For an open subset $U \subseteq M$, the restriction $\F \vert_U$ is defined as the set of connected components of the intersections of $U$ and leaves. 
Then the restriction $\F \vert_U$ is called the {\bf induced foliation} on $U$. 
\end{definition}

In other words, a subset $L \subset U$ is a leaf of $\F \vert_U$ if and only if $L$ is a connected component of $L_0 \cap U$ for some leaf $L_0 \in \F$. 
By construction, notice the induced foliation on $U$ is a singular foliation.  
We define the local differentiability as follows. 

\begin{definition}
A point is {\bf locally $\bm{C^{q}}$} with respect to a singular foliation if there is its open \nbd on which the induced foliation is $C^{q}$. 
\end{definition}


\subsection{One-dimensional singular foliations on manifolds}

To handle foliations consisting of piecewise linear curves that appear in numerical analysis and related analyses, we define the following concepts of piecewise $C^1$ and smoother structures for foliations.

\subsubsection{Piecewise $C^r$ singular foliations}

The piecewise $C^q$ singular foliations are defined as follows. 

\begin{definition}
For any $q \in \Z_{\geq 1} \sqcup \{ \infty \}$, a one-dimensional $C^{q-1}$ singular foliation $\F$ on a manifold $M$ is of {\bf piecewise $\bm{C^q}$} if each leaf is piecewise $C^q$ and any singular point is locally $C^q$ with respect to $\F$. 
\end{definition}

For the above terminology, $C^0$ singular foliations are also called {\bf piecewise $\bm{C^0}$} foliations.
Note that the concept of ``piecewise $C^0$ singular foliation'' is used for the simplicity of statements, and it is entirely equivalent to the concept of $C^0$ singular foliation.

For any $q \in \Z_{\geq 0} \sqcup \{ \infty \}$, denote by $\bm{\F^q(M)}$ the set of one-dimensional piecewise $C^q$ singular foliations on $M$. 
The set $\F^1(M)$ is also denoted by $\bm{\F(M)}$ in this paper. 
%
Denote by $\bm{\F_{\mathrm{reg}}^q(M)}$ the set of one-dimensional piecewise $C^q$ regular foliations on $M$ and put $\bm{\F_{\mathrm{reg}}(M)} := \F_{\mathrm{reg}}^1(M)$.

\subsection{Transverse to the boundary}

A regular leaf $L$ in a one-dimensional piecewise $C^1$-foliation $\F$ on a surface $M$ 
is {\bf transverse at a point $x \in L \cap \partial M$ to the boundary $\partial M$} with respect to $\F$ if the lift $\hat{x}$ of the point $x$ in the double $\hat{M}$ is contained in an open arc $\gamma$ in the lift of $\partial M$ which is transverse to the induced partition $\hat{\mathcal F}$ of $\mathcal F$ on $\hat{M}$ (i.e. $\gamma \pitchfork \hat{\mathcal F}$). 
In other words, letting $p \colon \hat{M} = M \times \{0,1\}/\mathop{\sim} \to M$ be the canonical projection, there is an open interval $\gamma \subset p^{-1}(\partial M)$ with $\hat{x} = p^{-1}(x) \in \gamma$ which is transverse to the induced partition $\hat{\mathcal F}$. 
This paper applies the transversality to the boundary only in the case described above (see \cite{HH1986A,CC2003I}. for a more general definition). 
Denote by $\partial_{\pitchfork \F}$ the set of points in $\partial M$ such that the leaves containing the points are transverse at the points to the boundary $\partial M$ with respect to $\F$. 
In other words, a point $x \in \partial M$ belongs to $\partial_{\pitchfork \F}$ if and only if the leaf $\F(x)$ is transverse at $x$ to $\partial M$ with respect to $\F$.

\section{Topologies on the space of one-dimensional foliations}

Let $(M,g)$ be a Riemannian manifold with the induced distance $d_g$ and $\F$ a singular foliation on $M$. 
%
From now on, assume that any singular foliation on a Riemannian manifold is one-dimensional and satisfies that the leaves are of piecewise $C^1$ unless otherwise stated.  

\begin{figure}[t]
\begin{center}
\includegraphics[scale=0.75]{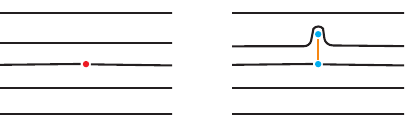}
\end{center} 
\caption{A bifurcation into a tripod and a thorn.}
\label{bifurcation_prong}
\end{figure}

As mentioned above, we would like to construct the foundation of time evaluations of divergence-free behaviors,
 which are described as divergence-free singular foliations (see \S~\ref{section:div-free} for definition). 
Though such time evaluations contain Whitehead moves as in Figure~\ref{Voronoi_prong} and a bifurcation, which is an analogous concept of a saddle-node bifurcation, as in Figure~\ref{bifurcation_prong}. 
We call such a bifurcation a {\bf tripod-thorn bifurcation}, because saddles (resp. centers) correspond to tripod (resp. thorn) in our setting. 
To describe such phenomena, we introduce a ``cone'' topology as below. 

\subsection{Topologies on the space of regular foliations}

We first define topologies on the space of regular foliations, and subsequently extend these definitions to the case of singular foliations.



\subsubsection{Polylines}

Fix a point $x \in M$. 
We identify tangent vectors $V \in T_x M$ as the closed unit intervals 
\[
\bm{[0_x, V]} := \{ s V \mid s \in [0,1]\} \subset T_x M
\]
 between the origin $0_x$ and $V$ in $T_x M$. 
Define $T_x M_{\pm}$ by 
\[
\bm{T_x M_{\pm}} := \{ \{ V, W \} \mid V,W \in T_x M \}
\]
and put $\bm{TM_\pm} := \bigsqcup_{y \in M} T_y M_\pm$. 
Then the set $T_x M_{\pm}$ is the set of two-point or one-point sets of tangent vectors of $T_x M$. 
Moreover, we also identify the two-point or one-point sets $\{ V, W \} \in T_x M_{\pm}$ with the unions $[0_x, V] \cup [0_x, W]$, which is called a {\bf polyline}.
Notice that these identifications of two-point sets preserve the Hausdorff distance between them. 
More precisely, we have that 
\[
d_H(\{V, W\}, \{V', W'\}) = d_H([0,V] \cup [0,W], [0,V'] \cup [0,W'])
\]
for any tangent vectors $V,W,V',W' \in T_x M \cong \R^n$, where $d_H$ is the Hausdorff distance on $T_x M$ with respect to the Riemannian distance $d_g$ and $[0,V] := \{ sV \mid s \in [0,1] \}$ for any $V \in \{V,W,V',W' \}$.

\subsubsection{Some notations related to foliated manifolds}

Denote by $T^1 M \subseteq TM$ the set of unit tangent vectors and by $T\F \subseteq TM$ the subbundle tangent to the leaves. 
Define $\bm{TM/\mathop{\pm}}$ by the quotient space by identifying any tangent vectors $V \in \bigsqcup_{x\in M}T_xM = TM$ and its inverses $-V \in T_xM$. 

\subsubsection{``Directed fields'' induced by regular foliations}

%
Let $\F$ be a $C^1$ regular foliation on $M$. 
Recall that the associated directed field 
\[
X_{\F} \colon M \to (T\F \cap T^1 M)/\pm \subset TM/\pm
\]
of $\F$ is the canonical section. 
Therefore the associated directed field $X_{\F}$ can be considered with the set $\{ \{-V_x, V_x \} \mid x \in M, V_x \in T_x\F \cap T_x^1 M \} \subset T_x M_\pm$ of pairs of the unit tangent vectors along leaves of $\F$. 

Similarly, let $\F$ be a regular foliation on $M$ whose leaves are of piecewise $C^1$. 
Define the associated directed field $X_{\F}$ as follows. 

\begin{figure}[t]
\begin{center}
\includegraphics[scale=0.95]{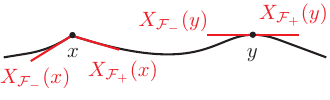}
\end{center} 
\caption{An example of pairs of vectors of the associated directed field.}
\label{fig:ass_dir_vf}
\end{figure}

\begin{definition}
A mapping 
\[
X_{\F} \colon M \to \{ \{X_{\F_-}(x), X_{\F_+}(x) \} \mid x \in M, X_{\F_-}(x), X_{\F_+}(x) \in  T_x^1 M \}  \subset TM_\pm
\] 
defined by $x \mapsto \{X_{\F_-}(x), X_{\F_+}(x) \} \in T_x M_\pm$ is called the {\bf associated directed field} of the regular foliation $\F$, where the pairs $\{X_{\F_-}(x), X_{\F_+}(x) \}$ of the distinct unit tangent vectors at any point $x \in M - \partial_{\pitchfork \F}$ of leaf arcs such that there are non-degenerate arc-length $C^1$ leaf arcs $C_{x,-}, C_{x,+} \colon [t_0,t_1] \to M$ starting from $x = C_{x,+}(t_0) = C_{x,+}(t_0)$ with 
\[
X_{\F_-}(x) = \lim_{s \to t_0+} \dfrac{C_{x,-}(s) - C_{x,-}(t_0)}{s - t_0} \in T_x M
\]
and 
\[
X_{\F_+}(x) = \lim_{s \to t_0+} \dfrac{C_{x,+}(s) - C_{x,+}(t_0)}{s - t_0} \in T_x M
\]
as in Figure~\ref{fig:ass_dir_vf}, and the singleton $\{X_{\F_-}(x), X_{\F_+}(x) \} = \{X_{\F_-}(x) \}$ of the unit tangent vector at any point $x \in \partial_{\pitchfork \F}$ of a leaf arc such that there is an non-degenerate arc-length $C^1$ leaf arc $C_{x,-} \colon [t_0,t_1] \to M$ starting from $x = C_{x,-}(t_0)$ satisfying the following equality:  
\[
X_{\F_-}(x) = \lim_{s \to t_0+} \dfrac{C_{x,-}(s) - C_{x,-}(t_0)}{s - t_0} \in T_x M
\]
\end{definition}

Recall that the pair $\{X_{\F_-}(x), X_{\F_+}(x) \} \subset T_x M_\pm$ is identified with a polyline $[0_x, X_{\F_-}(x)] \cup [0_x, X_{\F_+}(x)]$. 
By definitions, for any regular foliations $\F, \G$ on $M$ whose leaves are of piecewise $C^1$, we have $2 \geq \sup_{x \in M} d_H(X_{w,\F}(x), X_{w,\G}(x))$.  
Moreover, notice that if $\F$ is a $C^1$ regular foliation then the polyline $[0_x, X_{\F_-}(x)] \cup [0_x, X_{\F_+}(x)]$ for any point $x \in M - \partial_{\pitchfork \F}$ is a straight line $\{ s X_{\F_+}(x) \mid s \in [-1,1]\}$ with $-X_{\F_-}(x) = X_{\F_+}(x)$ and the polyline $[0_x, X_{\F_-}(x)] \cup [0_x, X_{\F_+}(x)]$ for any point $x \in \partial_{\pitchfork \F}$ is a straight line $\{ s X_{\F_-}(x) \mid s \in [0,1]\}$.
%

\subsubsection{Distance between regular foliations and their topology}

We define the following concepts. 

\begin{definition}
For any $p \in \R_{\geq 1}$ and any pair of two regular foliations $\F$ and $\mathcal{G}$, define the superior distance $\bm{D(\F, \mathcal{G})}$ and the $L^p$-distance $\bm{D^p(\F, \mathcal{G})}$ between them with respect to $(r, \alpha)$ as follows:  
\[
D(\F, \mathcal{G}) = D^\infty(\F, \mathcal{G}) := \sup_{x \in M} d_H(X_{\F}(x), X_{\G}(x))
\]
\[
D^p(\F, \mathcal{G}) := \left(\dfrac{1}{\mathrm{vol}_g(M)}  \int_M d_H(X_{\F}(x), X_{\G}(x))^p d \mathrm{vol}_g \right)^{1/p}
\]
\end{definition}

\begin{definition}
For any $p \in \R_{\geq 1}$, any $\delta \in \R_{> 0}$, and a regular foliation $\F$ on a  manifold $M$, define the {\bf open $\delta$-ball} $\bm{B^p_{D}(\F,\delta)}$ centered at $\F$ as follows:  
\[
B^p_{D}(\F,\delta) := \{ \mathcal{G} \in \F_{\mathrm{reg}}(M) \mid D^p(\F, \mathcal{G}) < \delta \}
\]
\end{definition}

\begin{definition}
We define the set $\mathcal{B}^p_{\F_{\mathrm{reg}}(M)}$ of open balls as follows: 
\[
\mathcal{B}^p_{\F_{\mathrm{reg}}(M)} := \{ B^p_{D}(\F,\delta) \mid  \delta \in \R_{> 0}, \F \in \F_{\mathrm{reg}}(M) \}
\]
\end{definition}

\begin{definition}
For any $p \in \R_{\geq 1} \sqcup \{ \infty \}$, the smallest topology on $\F_{\mathrm{reg}}(M)$ which contains the set $\mathcal{B}^p_{\F_{\mathrm{reg}}(M)} $ of open balls is called the {\bf topology induced by the local Hausdorff $L^p$-distance} and denoted by $\bm{\mathcal{O}^p_{\F_{\mathrm{reg}}^q(M)}}$. 
\end{definition}

Then the set $\mathcal{B}^p_{\F_{\mathrm{reg}}(M)} $ is a subbase of the topology $\mathcal{O}^p_{\F_{\mathrm{reg}}^q(M)}$. 
Notice that we can define a distance $D^\infty$ (resp. $D^p$) between $k$-dimensional singular foliations by the supremum (resp. integral) of the Hausdorff distances between the normed squares tangent to leaves as well as one between one-dimensional singular foliations. 
These constructions will be addressed in future work.

\subsection{Topologies on the space of singular foliations}

To extend the topologies on the space of regular foliations to those for singular foliations, we define a weight function for singularities.

\subsubsection{Weight for singularities}

To define ``directed fields'' induced by singular foliations, we define the weight for singularities as follows. 

\begin{definition}
A non-decreasing continuous function $w \colon [0,\infty] \to [0,1]$ is a {\bf distance weight} if $w^{-1}(0) = \{ 0 \}$ and $w(\infty) = 1$. 
\end{definition}

Note that any distance weight is uniformly continuous because it is bounded, monotone increasing, and continuous. 

\begin{definition}
For a distance weight $w \colon [0,\infty] \to [0,1]$, 
the {\bf weight} $\bm{\rho_{w}} \colon M \times \F(M) \to [0,1]$ with respect to $w$ is defined as follows:
\[
\rho_{w}(x, \F) = w(d_g(x, \mathop{\mathrm{Sing}}(\F)))
\]
\end{definition}

Put $\rho_{w,\F} := \rho_{w}(\cdot, \F) \colon M \to [0,1]$. 
Note that $\rho_{w}(\cdot, \F) \equiv 1$ for any regular foliation $\F$ and any distance weight $w$. 

\subsubsection{Associated ``directed fields'' of singular foliations}

Let $\F$ be a singular foliation on a Riemannian manifold $(M,g)$ whose leaves are of piecewise $C^1$. 
Since $\F_{\mathrm{reg}} := \F\vert_{\F - \mathop{\mathrm{Sing}}(\F)}$ is a regular foliation on a Riemannian manifold $M - \mathop{\mathrm{Sing}}(\F)$, denote by 
\[
X_{\F_{\mathrm{reg}}} \colon M  - \mathop{\mathrm{Sing}}(\F) \to \{ [0_x, X_{\F_-}(x)] \cup [0_x, X_{\F_{+}}(x)] \mid x \in M  - \mathop{\mathrm{Sing}}(\F) \} \subset TM_\pm
\]
the associated directed field of $\F_{\mathrm{reg}}$, where 
\[
X_{\F_{-}}(x) := X_{\F_{\mathrm{reg}-}}(x) \text{ and } X_{\F_{+}}(x) := X_{\F_{\mathrm{reg}+}}(x)
\]
for any $x \in M - \mathop{\mathrm{Sing}}(\F)$. 
For any $x \in \mathop{\mathrm{Sing}}(\F)$, put $X_{\F_{-}}(x)  = X_{\F_{+}}(x)  := 0 \in T_x M$. 
Define the associated directed field $X_{\F} \colon M \to TM_\pm$ of the foliation $\F$ by $X_{\F}\vert_{\mathop{\mathrm{Sing}}(\F)} \equiv 0$ and 
\[
X_{\F}(x) := X_{\F_{\mathrm{reg}}}(x) = [0_x, X_{\F_-}(x)] \cup [0_x, X_{\F_+}(x)],
\]
which need not be continuous. 
Notice that the polyline $X_{\F}(x) = [0_x, X_{\F_-}(x)] \cup [0_x, X_{\F_+}(x)]$ is identified with $\{X_{\F_-}(x), X_{\F_+}(x)\}$ for any $x \in M$. 

For the weight $\rho_{w,\F} \colon M \to [0,1]$ of a distance weight $w \colon [0,\infty] \to [0,1]$, put 
\[
X_{w,\F_-}(x) :=  \rho_{w,\F}(x)X_{\F_-}(x)
\hspace{20pt}
X_{w,\F_+}(x) :=  \rho_{w,\F}(x)X_{\F_+}(x)
\]
for any $x \in M$. 
We define the continuous associated directed field of $\F$ with respect to the weight $\rho_{w,\F} \colon M \to [0,1]$ as follows. 

\begin{definition}
Define the {\bf continous associated directed field} $X_{w,\F} \colon M \to TM_\pm$ of $\F$ with respect to the weight $\rho_{w,\F} \colon M \to [0,1]$ by 
\[
\begin{split}
X_{w,\F}(x) := \rho_{w,\F}(x) X_{\F}(x) & = \{ \rho_{w,\F}(x) X_{\F_-}(x), \rho_{w,\F}(x) X_{\F_+}(x)\}
\\
& = \{ X_{w,\F_-}(x), X_{w,\F_+}(x) \} 
\end{split}
\]
for any $x \in M$.
\end{definition}
By definition, any continuous associated directed field is continuous. 
Moreover, we obtain that $X_{w,\F}(x) = \{ 0_x \}$ for any $x \in \mathop{\mathrm{Sing}}(\F)$, where $0_x \in T_x M$ is the zero tangent vector, and that 
\[
\begin{split}
X_{w,\F}(x) = \rho_{w,\F}(x) X_{\F}(x) &= [0_x, \rho_{w,\F}(x)X_{\F_-}(x)] \cup [0_x, \rho_{w,\F}(x)X_{\F_+}(x)]
\\
&= [0_x, X_{w,\F_-}(x)] \cup [0_x, X_{w,\F_+}(x)]
\end{split}
\]
for any $x \in M - \mathop{\mathrm{Sing}}(\F)$.

\subsubsection{Distance between singular foliations and their topology}

Fix a distance weight $w \colon [0,\infty] \to [0,1]$. 
We define a distance as follows. 

\begin{definition}
Define the wighted metric $\bm{D_{w}} \colon \F(M) \times \F(M) \to [0, \infty]$, called the {\bf locally wighted Hausdorff superior distance} on the set $\F(M)$ of singular foliations on a  manifold $M$ for an ordered field $[0, \infty]$: 
\[
\begin{split}
\bm{D_{w}(\F, \mathcal{G})} = D^\infty_{w}(\F, \mathcal{G}) &:= \sup_{x \in M} d_H(X_{w,\F}(x), X_{w,\G}(x))
\end{split}
\]
\end{definition}

Notice that $2 \geq \sup_{x \in M} d_H(X_{w,\F}(x), X_{w,\G}(x))$ for any  singular foliations $\F, \G$ on $M$ whose leaves are of piecewise $C^1$, and that $D_{w}(\F, \mathcal{G}) = D(\F, \mathcal{G})$ for any regular foliations $\F, \mathcal{G} \in \F_{\mathrm{reg}}(M)$. 
%
Similarly, we define distances for any $p \in \R_{\geq 1}$ as follows. 

\begin{definition}
For any $p \in \R_{\geq 1}$,  
and for a regular foliation $\F$ on a  manifold $M$, define the wighted metric $\bm{D^p_{w}} \colon \F(M) \times \F(M) \to [0, \infty]$, called the {\bf locally wighted Hausdorff $L^p$-distance} on the set $\F(M)$ of singular foliations on a  manifold $M$ for an ordered field $[0, \infty]$: 
\[
\bm{D^p_{w}(\F, \mathcal{G})}:= \left(\dfrac{1}{\mathrm{vol}_g(M)}  \int_M d_H(X_{w,\F}(x), X_{w,\G}(x))^p d \mathrm{vol}_g \right)^{1/p}
\]
\end{definition}

\subsubsection{Open $\delta$-balls for singular foliations and its topology}
Fix a distance weight $w \colon [0,\infty] \to [0,1]$ and any number $p \in \R_{\geq 1} \sqcup \{ \infty \}$. 
We define open balls with respect to $w$ and $p$ as follows. 

\begin{definition}
For any $\delta \in \R_{> 0}$ and any $p \in \R_{\geq 1} \sqcup \{ \infty \}$, and a singular foliation $\F$ on a  manifold $M$, define the {\bf open $\delta$-ball} $\bm{B_{D^p_{w}}(\F,\delta)}$ centered at $\F$ with respect to $w$ as follows:  
\[
B_{D^p_{w}}(\F,\delta) := \{ \mathcal{G} \in \F(M) \mid D^p_{w}(\F, \mathcal{G}) < \delta \}
\]
\end{definition}

Moreover, we define the set $\bm{\mathcal{B}^p_{w,\F(M)}}$ of open balls with respect to $w$ as follows: 
\[
\mathcal{B}^p_{w,\F(M)} := \{ B_{D^p_{w}}(\F,\delta) \mid \delta \in \R_{> 0}, \F \in \F(M) \}
\]

\subsubsection{Topologies for singular foliations}

We define a topology on the set $\F(M)$ of one-dimensional singular foliations on a Riemannian manifold $M$. 

\begin{definition}\label{def:topology}
For any distance weight $w \colon [0,\infty] \to [0,1]$ and any $p \in \R_{\geq 1} \sqcup \{ \infty \}$, 
the smallest topology on $\F^q(M)$ which contains the set $\mathcal{B}^p_{w,\F(M)} $ of open balls  with respect to $w$ is called the {\bf topology induced by the local Hausdorff distance} with respect to $w$ and denoted by $\bm{\mathcal{O}^p_{w,\F^q(M)}}$. 
\end{definition}

Then the set $\mathcal{B}^p_{w,\F(M)} $ is a subbase of the topology $\mathcal{O}^p_{w,\F^q(M)}$.


\section{Fundamental properties}

From now on, 
fix any distance weight $w \colon [0,\infty] \to [0,1]$ unless otherwise stated.
Moreover, we equip the sets of one-dimensional singular foliations on Riemannian manifolds with the topology induced by the local Hausdorff distance with respect to the distance weight $w$ unless otherwise stated. 

First, we state the fundamental properties of topologies.  
Let $\F$ be a one-dimensional singular foliation on a Riemannian manifold $M$ with or without boundary whose leaves are of piecewise $C^1$. 
%
%
To define structural stability, we recall the topological equivalence.

\begin{definition}
A singular foliation $\F_1$ on a manifold $M_1$ is {\bf topologically equivalent} to a singular foliation $\F_2$ on a manifold $M_2$ if there is a homeomorphism $h \colon M_1 \to M_2$ such that the image of any leaf of $\F_1$ is a leaf of $\F_2$. 
\end{definition}

We define structural stability for singular foliations as follows. 

\begin{definition}
For any $q \in \Z_{\geq 1} \sqcup \{ \infty \}$ and a subspace $\chi(M) \subseteq \F^q(M)$, a singular foliation $\F \in \chi(M)$ is {\bf structurally stable} on $\chi(M)$ if there is its open neighborhood in $\chi(M)$ such that $\F$ is topologically equivalent to any singular foliation in the neighborhood. 
\end{definition}

To define locally topological equivalence, we recall the restriction as follows. 

\begin{definition}
The {\bf restriction $\F \vert_{A}$} of a partition $\F$ on a topological space $M$ to a subset $A \subseteq M$ is the set of connected components of the intersections of $A$ and leaves of $\F$. 
\end{definition}

Notice that a subset $l \subseteq A$ is an element of the restriction $\F \vert_{A}$ if and only if there is an element $L \in \F$ such that $l$ is a connected component of $L \cap A$. 
We define the locally topological equivalence as follows.

\begin{definition}
The restriction $\F_1 \vert_{A_1}$ of a partition $\F_1$ on a topogical space $M_1$ to a subset $A_1 \subseteq M_1$ is {\bf {\rm(}locally{\rm)} topologically equivalent} to the restriction $\F_2 \vert_{A_2}$ of a partition $\F_2$ on a topogical space $M_2$ to a subset $A_2 \subseteq M_2$ if there is a homeomorphism $h \colon A_1 \to A_2$ such that the image of any element of $\F_1 \vert_{A_1}$ is elements of $\F_2 \vert_{A_2}$. 
\end{definition}

Notice that restriction $\F_1 \vert_{A_1}$ of a singular foliation $\F_1$ on a manifold $M_1$ to a subset $A_1 \subseteq M_1$ is topologically equivalent to the restriction $\F_2 \vert_{A_2}$ of a singular foliation $\F_2$ on a manifold $M_2$ to a subset $A_2 \subseteq M_2$ if and only if there is a homeomorphism $h \colon A_1 \to A_2$ such that the image by $h$ of any connected components of the intersection $L_1 \cap A_1$ of a leaf $L_1$ of $\F_1$ and $A_1$ is a connected components of the intersection $L_2 \cap A_2$ of some leaf $L_2$ of $\F_2$ and $A_2$.

\subsection{Persistence of local non-existence of singular leaves}

For any distance weight $w \colon [0,\infty] \to [0,1]$, define continuous functions $\delta, \varepsilon_w \colon M \to [0,1]$ as follows: 
\[
\delta(x) := d(x,\mathop{\mathrm{Sing}}(\F))/2 > 0
\]  
\[
\varepsilon_w(x) =  w(d_g(x, \mathop{\mathrm{Sing}}(\F))/2)/2 =  w(\delta(x))/2
\]

When the distance weight $w$ is clear, the notation $\varepsilon_w(x)$ is abbreviated as $\varepsilon(x)$.
The local non-existence of singular leaves persists as follows. 

\begin{lemma}\label{lem:open_reg_pt}
For any distance weight $w \colon [0,\infty] \to [0,1]$, the following statements hold:
\\
{\rm(1)} $\varepsilon_w^{-1}(0) = \delta^{-1}(0) = \mathop{\mathrm{Sing}}(\F)$. 
\\
{\rm(2)} $B(x_0,\delta(x_0)) \cap \bigcup_{\mathcal{G} \in B_{D^\infty_{w}}(\F,\varepsilon_w(x_0))}  \mathop{\mathrm{Sing}}(\mathcal{G}) = \emptyset$. 
\\
{\rm(3)} For any points $x,y \in M$ with $d(x, \mathop{\mathrm{Sing}}(\F)) \leq d(y, \mathop{\mathrm{Sing}}(\F))$, we have that $\varepsilon_w (x) \leq \varepsilon_w (y)$. 
\\
{\rm(4)} For any sequence $(x_n)_{n \in \Z_{>0}}$ with $\lim_{n \to \infty} d(x_n, \mathop{\mathrm{Sing}}(\F)) = \infty$, we obtain that $\lim_{n \to \infty} \varepsilon_w (x_n) = 1/2$.  
\end{lemma}

\begin{proof}
Fix any distance weight $w \colon [0,\infty] \to [0,1]$.
If $M$ has a boundary, then fix an extension $M'$ of $M$, which is an open manifold containing $M$ in the interior. 
By $w^{-1}(0) = \{ 0 \}$, by definitions of $\delta$ and $\varepsilon_w$, we have $\varepsilon_w^{-1}(0) = \delta^{-1}(0) = \mathop{\mathrm{Sing}}(\F)$, which implies assertion (1).
From definition of $\varepsilon_w$, assertions (2) and (3) hold.

Define a function $\varepsilon'_w \colon M \to [0,1]$ by $\varepsilon'_w \vert_{\mathop{\mathrm{Sing}}(\F)} = 0$ and 
\[
\varepsilon'_w(x_0) := \inf \{ \rho_{w,\F}(x) \mid x \in B(x_0,\delta(x_0)) \}
\]
for any $x_0 \in M - \mathop{\mathrm{Sing}}(\F)$. 

\begin{claim}\label{claim:20}
$2 \varepsilon_w = \varepsilon'_w$. 
\end{claim}

\begin{proof}[Proof of Claim~\ref{claim:20}]
For any $x_0 \in \mathop{\mathrm{Sing}}(\F)$, by $\rho_{w,\F}(x_0) = w(d_g(x_0, \mathop{\mathrm{Sing}}(\F)))$ and $\delta(x) = d(x,\mathop{\mathrm{Sing}}(\F))/2$, 
assertion (1) implies the following equality: 
\[
\varepsilon'_w(x_0) = \inf \{ \rho_{w,\F}(x) \mid x \in B(x_0,\delta(x_0)) \} = \rho_{w,\F}(x_0) = 0 = 2\varepsilon(x_0)
\]

Fix any $x_0 \in M - \mathop{\mathrm{Sing}}(\F)$. 
From $\rho_{w,\F}(x_0) = w(d_g(x_0, \mathop{\mathrm{Sing}}(\F)))$ and $\delta(x) = d(x,\mathop{\mathrm{Sing}}(\F))/2$, we have the following equality: 
\[
\begin{split}
\varepsilon'_w(x_0) &= \inf \{ \rho_{w,\F}(x) \mid x \in B(x_0,\delta(x_0)) \}\\
&= \inf \{ w(d_g(x, \mathop{\mathrm{Sing}}(\F))) \mid d(x,x_0) \leq \delta(x_0) \}  \\
&= \inf \{ w(d_g(x, \mathop{\mathrm{Sing}}(\F))) \mid d(x,x_0) \leq d(x_0,\mathop{\mathrm{Sing}}(\F))/2 \}
\end{split}
\]
Therefore, the definition $d_g$ implies 
\[
\varepsilon'_w(x_0) =  w(d_g(x_0, \mathop{\mathrm{Sing}}(\F))/2) =  w(\delta(x_0)) = 2 \varepsilon_w(x_0)
\]
\end{proof}
%
%
%

\begin{claim}\label{claim:001}
$B(x_0,\delta(x_0)) \cap \mathop{\mathrm{Sing}}(\mathcal{G}) = \emptyset$ for any $\mathcal{G} \in B_{D^\infty_{w}}(\F,\varepsilon(x_0))$. 
\end{claim}

\begin{proof}[Proof of Claim~\ref{claim:001}]
Fix any point $x_0 \in M - \mathop{\mathrm{Sing}}(\F)$. 
For any point $x \in B(x_0,\delta(x_0))$, we have 
\[
\varepsilon'_w(x_0) X_{\F}(x) \subseteq \rho_{w,\F}(x) X_{\F}(x) = X_{w,\F}(x)
\]
as polylines, and so $d_H(X_{w,\F}(x), \{ 0_x \}) \geq d_H(\varepsilon'_w(x_0) X_{\F}(x), \{ 0_x \}) =  \varepsilon'_w(x_0) > \varepsilon(x_0)$.  
For any singular foliation $\mathcal{G}$ with $\mathop{\mathrm{Sing}}(\mathcal{G}) \cap B(x_0,\delta(x_0)) \neq \emptyset$, 
we have the following inequality: 
\[
\begin{split}
D^\infty_{w}(\F, \mathcal{G}) & = \sup_{x \in M} d_H(X_{w,\F}(x), X_{w,\G,x}) \\
& \geq \sup_{x \in \mathop{\mathrm{Sing}}(\mathcal{G}) \cap B(x_0,\delta(x_0))} d_H(X_{w,\F}(x) , \{ 0_x \}) \\
& \geq \sup_{x \in \mathop{\mathrm{Sing}}(\mathcal{G}) \cap B(x_0,\delta(x_0))} d_H(\varepsilon'_w(x_0) X_{\F}(x), \{ 0_x \})  = \varepsilon'_w(x_0)  > \varepsilon(x_0)
\end{split}
\]
This implies the claim. 
\end{proof}
The previous claim implies the assertion. 
\end{proof}

We observe the non-existence of singular points outside of \nbds of singular points under perturbations as follows.

\begin{corollary}\label{cor:nonexistence_sing}
For any distance weight $w \colon [0,\infty] \to [0,1]$ and any positive number $\delta_0 > 0$, there is a positive number $\varepsilon_0>0$ satisfying the following relation: 
\[
\bigcup_{\mathcal{G} \in B_{D^\infty_{w}}(\F,\varepsilon_0)}  \mathop{\mathrm{Sing}}(\mathcal{G}) \subseteq B_{\delta_0}(\mathop{\mathrm{Sing}}(\F)))
\]
\end{corollary}

\begin{proof}
Fix any distance weight $w \colon [0,\infty] \to [0,1]$ and any $\delta_0 > 0$. 
Put $\varepsilon_0 := w(\delta_0/2)/3$. 
Fix any $x_0 \in M - B_{\delta_0}(\mathop{\mathrm{Sing}}(\F))$. 
Then the following inequalities hold: 
\[
\delta_w(x_0) = d(x_0,\mathop{\mathrm{Sing}}(\F))/2 \geq \delta_0/2
\]  
\[
\begin{split}
\varepsilon_w(x_0) & = \min \{ \rho_{w,\F}(x) \mid x \in B(x_0,\delta_w(x_0)) \}/2 
\\
& \geq w(d(x_0,\mathop{\mathrm{Sing}}(\F))/2)/2
\\
& = w(\delta_w(x_0))/2 \geq w(\delta_0/2)/2 > \varepsilon_0
\end{split}
\]

By Lemma~\ref{lem:open_reg_pt},  we have 
\[
B(x_0,\delta_w(x_0)) \cap \bigcup_{\mathcal{G} \in B_{D^\infty_{w}}(\F,\varepsilon_0)}  \mathop{\mathrm{Sing}}(\mathcal{G}) = \emptyset
\] 
and so 
\[
\bigcup_{x \in M - B_{\delta_0}(\mathop{\mathrm{Sing}}(\F))} B(x,\delta_w(x)) \cap \bigcup_{\mathcal{G} \in B_{D^\infty_{w}}(\F,\varepsilon_0)}  \mathop{\mathrm{Sing}}(\mathcal{G}) = \emptyset
\] 
This implies the following equality: 
\[
(M - B_{\delta_0}(\mathop{\mathrm{Sing}}(\F))) \cap \bigcup_{\mathcal{G} \in B_{D^\infty_{w}}(\F,\varepsilon_0)}  \mathop{\mathrm{Sing}}(\mathcal{G}) = \emptyset
\]
Therefore, the assertion holds. 
\end{proof}

The local non-existence of singular leaves persists with respect to $\mathcal{O}^p_{w,\F^q(M)}$ ($p \in \R_{\geq 1}$) as follows. 

\begin{lemma}\label{lem:open_reg_pt_03}
Let $p \in \R_{\geq 1}$ be a number. 
Suppose that the volume of $M$ is finite. 
For any distance weight $w \colon [0,\infty] \to [0,1]$, there are continuous functions $\varepsilon_1, \delta_1 \colon M \to \R_{\geq 0}$ with $\varepsilon_1^{-1}(0) = \delta_1^{-1}(0) = \mathop{\mathrm{Sing}}(\F)$ satisfying the following equation for any $x_0 \in M - \mathop{\mathrm{Sing}}(\F)$: 
\[
B(x_0,\delta_1(x_0)) \cap \bigcup_{\mathcal{G} \in B_{D^p_{w}}(\F,\varepsilon_1(x_0))}  \mathop{\mathrm{Sing}}(\mathcal{G}) = \emptyset
\] 
\end{lemma}

\begin{proof}
If $M$ has a boundary, then fix an extension $M'$ of $M$, which is an open manifold containing $M$ in the interior. 
Fix any distance weight $w \colon [0,\infty] \to [0,1]$.
For any $x \in \mathop{\mathrm{Sing}}(\F)$, define $\varepsilon_1(x) = \delta_1(x) := 0$. 
Lemma~\ref{lem:open_reg_pt} implies that the continuous functions $\varepsilon_w, \delta_w \colon M \to \R_{\geq 0}$ with $\varepsilon_w^{-1}(0) = \delta_w^{-1}(0) = \mathop{\mathrm{Sing}}(\F)$ satisfy the following equation for any $x_0 \in M - \mathop{\mathrm{Sing}}(\F)$: 
\[
B(x_0,\delta_w(x_0)) \cap \bigcup_{\mathcal{G} \in B_{D^\infty_{w}}(\F,\varepsilon_w(x_0))}  \mathop{\mathrm{Sing}}(\mathcal{G}) = \emptyset
\] 
Then the functions $\delta_w/n \colon M \to \R_{\geq 0}$ are continuous for any $n \in \Z_{>1}$. 
For any integer $n \in \Z_{>1}$ and any $x_0 \in M - \mathop{\mathrm{Sing}}(\F)$, define positive numbers 
\[
d_{n}(x_0) := \inf_{x \in B(x_0,\delta_w(x_0)/n)} d_g(x,\mathop{\mathrm{Sing}}(\mathcal{F})) \geq  \delta_w(x_0) (n-1)/n > 0
\]
and $w_n(x_0) := w(d_{n}(x_0)) > 0$. 
For any integer $n \in \Z_{>1}$, putting $w_n \vert_{\mathop{\mathrm{Sing}}(\mathcal{F})} := 0$, by constructions of $d_n$ and $w_n$, the extended function $w_n \colon M \to \R_{\geq 0}$ is continuous. 

\begin{claim}\label{claim:01a}
For any point $x_0 \in M - \mathop{\mathrm{Sing}}(\F)$, there is a positive number $\delta_{1,x_0}$ such that $\inf_{x \in B(x_0,\delta_{1,x_0})} d_H(X_{w,\F}(x), X_{w,\G}(x)) >\sup_{x \in B(x_0,\delta_{1,x_0})} \rho_{w,\F}(x)/4$ for any $\mathcal{G} \in \mathcal{F}(M)$ with $\mathop{\mathrm{Sing}}(\mathcal{G}) \cap B(x_0,\delta_{1,x_0}) \neq \emptyset$. 
\end{claim}

\begin{proof}[Proof of Claim~\ref{claim:01a}]
Fix a point $x_0 \in M - \mathop{\mathrm{Sing}}(\F)$. 
Then we have $\lim_{n \to \infty} w_n(x_0) = w(d_g(x_0, \mathop{\mathrm{Sing}}(\F))) = \rho_{w,\F}(x_0) > 0$. 
Moreover, we have the following inequalities for any $n \in \Z_{>1}$: 
\[
\begin{split}
\inf_{x \in B(x_0,\delta_w(x_0)/n)} \vert X_{w,\F}(x) \vert & = \inf_{x \in B(x_0,\delta_w(x_0)/n)} \vert \rho_{w,\F}(x) \vert 
\\
& = \inf_{x \in B(x_0,\delta_w(x_0)/n)} w(d_g(x, \mathop{\mathrm{Sing}}(\F))) = w_n(x_0)
\end{split}
\]
Since $\lim_{n \to \infty} \sup_{x \in B(x_0,\delta_w(x_0)/n)} \rho_{w,\F}(x)/2 = \rho_{w,\F}(x_0)/2 <  \lim_{n \to \infty} w_n(x_0)$, there is a positive integer $N_1(x_0)$ such taht 
\[
\inf_{x \in B(x_0,\delta_w(x_0)/n)} \vert X_{w,\F}(x) \vert > \sup_{x \in B(x_0,\delta_w(x_0)/n)} \dfrac{\rho_{w,\F}(x)}{2} > 0
\]
for any $n \in \Z_{>N_1(x_0)}$. 
On the other hand, define a positive constant $M(n,x_0)$ as follows: 
\[
M(n,x_0) := \sup_{\mathcal{G} \in \F(M) \text{ with } \mathop{\mathrm{Sing}}(\mathcal{G}) \cap B(x_0,\delta_w(x_0)/n) \neq \emptyset} \sup_{x \in B(x_0,\delta_w(x_0)/n)} \vert X_{w,\G}(x) \vert
\] 
Since the norm $\vert X_{w,\G}(x) \vert$ for any $\mathcal{G} \in \F(M)$ and for any $x \in M$ satisfies $\vert X_{\G}(x) \vert = \vert \rho_{w,\G}(x) \vert$, we have the following equality: 
\[
\begin{split}
\lim_{n \to \infty} M(n,x_0) = & \lim_{n \to \infty} \sup_{\mathcal{G} \in \F(M) \text{ with } \mathop{\mathrm{Sing}}(\mathcal{G}) \cap B(x_0,\delta_w(x_0)/n) \neq \emptyset} \sup_{x \in B(x_0,\delta_w(x_0)/n)} \vert X_{w,\G}(x) \vert
\\
= & \lim_{n \to \infty} \sup_{\mathcal{G} \in \F(M) \text{ with } \mathop{\mathrm{Sing}}(\mathcal{G}) \cap B(x_0,\delta_w(x_0)/n) \neq \emptyset} \sup_{x \in B(x_0,\delta_w(x_0)/n)} \vert \rho_{w,\G}(x) \vert 
\\
 = & \, \, 0
\end{split}
\]
By $\lim_{r \to \infty} w(1/r)=0$, since $\lim_{n \to \infty} \inf_{x \in B(x_0,\delta_w(x_0)/n)} \rho_{w,\F}(x) = \rho_{w,\F}(x_0) > 0$, there is a positive integer $N_2(x_0) \geq N_1(x_0)$ such taht 
\[
\inf_{x \in B(x_0,\delta_w(x_0)/n)} \dfrac{\rho_{w,\F}(x_0)}{4} > M(n,x_0)
\]
for any $n \in \Z_{>N_2(x_0)}$. 
Fix a positive integer $n(x_0) \in \Z_{>N_2(x_0)}$. 
Put $\delta_{1,x_0} :=\delta(x_0)/n(x_0)$. 
Then the following inequalities holds for any $\mathcal{G} \in \mathcal{F}(M)$ with $\mathop{\mathrm{Sing}}(\mathcal{G}) \cap B(x_0,\delta_{1,x_0}) \neq \emptyset$:  
\[
\begin{split}
 \inf_{x \in B(x_0,\delta_{1,x_0})} \vert X_{w,\F}(x) \vert > \sup_{x \in B(x_0,\delta_{1,x_0})} \dfrac{\rho_{w,\F}(x)}{2} & > \inf_{x \in B(x_0,\delta_{1,x_0})} \dfrac{\rho_{w,\F}(x)}{4} 
 \\
 & > M(n(x_0),x_0) 
 \\
 & \geq \sup_{x \in B(x_0,\delta_{1,x_0})} \vert X_{w,\G}(x) \vert
\end{split}
\]
Therefore, we have the following inequality for any $\mathcal{G} \in \mathcal{F}(M)$ with $\mathop{\mathrm{Sing}}(\mathcal{G}) \cap B(x_0,\delta_{1,x_0}) \neq \emptyset$: 
\[
\inf_{x \in B(x_0,\delta_{1,x_0})} d_H(X_{w,\F}(x), X_{w,\G}(x)) > \sup_{x \in B(x_0,\delta_{1,x_0})} \rho_{w,\F}(x)/4
\]
\end{proof}

By the paracompactness of $M$ and the closedness of $\mathop{\mathrm{Sing}}(\F)$, there is a locally finite open cover $\{\mathop{\mathrm{int}} B(x_{0,\lambda},\delta_{1,x_{0,\lambda}}) \}_{\lambda \in \Lambda}$ of $M - \mathop{\mathrm{Sing}}(\F)$. 
For any $x \in M - \mathop{\mathrm{Sing}}(\F)$, the subset $\Lambda_x := \{ \lambda \in \Lambda \mid x \in \mathop{\mathrm{int}} B(x_{0,\lambda}, \delta_{1,x_{0,\lambda}}) \}$ is finite. 
Therefore, define a function $\delta_1 \colon M \to \R_{\geq 0}$ by $\delta^{-1}(0) = \mathop{\mathrm{Sing}}(\F)$ and 
\[
\delta_1(x) := \max \{ \delta_{1,x_{0,\lambda}} - d(x_{0,\lambda}, x) \mid \lambda \in \Lambda_x \}
\]
 for any $x \in M - \mathop{\mathrm{Sing}}(\F)$. 

\begin{claim}\label{claim:03a}
$\inf_{x \in B(x_0,\delta_1(x_0))} d_H(X_{w,\F}(x), X_{w,\G}(x)) \geq \rho_{w,\F}(x_{0})/4$ for any $x_0 \in M - \mathop{\mathrm{Sing}}(\F)$.
\end{claim}

\begin{proof}[Proof of Claim~\ref{claim:03a}]
By construction, the function $\delta_1$ is continuous. 
For any point $x_0 \in M - \mathop{\mathrm{Sing}}(\F)$, there is an index $\lambda \in \Lambda_{x_0} \subseteq \Lambda$ such that $B(x_0,\delta_1(x_0)) \subseteq B(x_{0,\lambda}, \delta_{1,x_{0,\lambda}})$ and so 
\[
\begin{split}
& \inf_{x \in B(x_0,\delta_1(x_0))} d_H(X_{w,\F}(x), X_{w,\G}(x)) 
\\
\geq & \inf_{x \in B(x_{0,\lambda}, \delta_{1,x_{0,\lambda}})} d_H(X_{w,\F}(x), X_{w,\G}(x))
\\
\geq & \sup_{x \in B(x_{0,\lambda}, \delta_{1,x_{0,\lambda}})} \rho_{w,\F}(x)/4 \geq \rho_{w,\F}(x_{0})/4
\end{split}
\]
because of Claim~\ref{claim:01a}.
\end{proof}
Define a continuous function $\varepsilon_1 \colon M \to [0,1/5]$ as follows: 
\[
\varepsilon_1(x_0) := \dfrac{\rho_{w,\F}(x_0)}{5} \left(\dfrac{1}{\mathrm{vol}_g(M)} \mathrm{vol}_g(B(x_0,\delta_1(x_0))) \right)^{1/p} 
\]

\begin{claim}\label{claim:02a}
For any point $x_0 \in M - \mathop{\mathrm{Sing}}(\F)$ and any $\mathcal{G} \in \mathcal{F}(M)$ with $\mathop{\mathrm{Sing}}(\mathcal{G}) \cap B(x_0,\delta_1(x_0)) \neq \emptyset$, the following relation holds:  
\[
\mathcal{G} \notin B_{D^p_{w}}(\F,\varepsilon_1(x_0)))
\] 
\end{claim}

\begin{proof}[Proof of Claim~\ref{claim:02a}]
Fix any point $x_0 \in M - \mathop{\mathrm{Sing}}(\F)$ and any singular foliation $\mathcal{G} \in \mathcal{F}(M)$ with $\mathop{\mathrm{Sing}}(\mathcal{G}) \cap B(x_0,\delta_1(x_0)) \neq \emptyset$. 
Then Claim~\ref{claim:03a} implies  
\[
\begin{split}
& D^p_{w}(\F, \mathcal{G}) \\
= & \left(\dfrac{1}{\mathrm{vol}_g(M)}  \int_M d_H(X_{w,\F}(x), X_{w,\G}(x))^p d \mathrm{vol}_g \right)^{1/p}\\
\geq & \left(\dfrac{1}{\mathrm{vol}_g(M)}  \int_{B(x_0,\delta_1(x_0))} d_H(X_{w,\F}(x), X_{w,\G}(x))^p d \mathrm{vol}_g \right)^{1/p}
\\
\geq & \left(\dfrac{1}{\mathrm{vol}_g(M)}  \int_{B(x_0,\delta_1(x_0))} \left( \dfrac{\rho_{w,\F}(x_0)}{4} \right)^p d \mathrm{vol}_g \right)^{1/p} \\
\geq & \dfrac{\rho_{w,\F}(x_0)}{4} \left(\dfrac{1}{\mathrm{vol}_g(M)}  \int_{B(x_0,\delta_1(x_0))}  d \mathrm{vol}_g \right)^{1/p} \\
= & \dfrac{\rho_{w,\F}(x_0)}{4} \left(\dfrac{1}{\mathrm{vol}_g(M)} \mathrm{vol}_g(B(x_0,\delta_1(x_0))) \right)^{1/p} > \varepsilon_1(x_0)
\end{split}
\]
because of $\varepsilon_1(x_0) = \dfrac{\rho_{w,\F}(x_0)}{5} \left(\dfrac{1}{\mathrm{vol}_g(M)} \mathrm{vol}_g(B(x_0,\delta_1(x_0))) \right)^{1/p}$. 
This means $\mathcal{G} \notin B_{D^p_{w}}(\F,\varepsilon_1(x_0)))$. 
\end{proof}
This completes the proof. 
\end{proof}

\subsubsection{Persistence of regular foliations}

We have the following persistence of regular foliations. 

\begin{lemma}\label{lem:open_reg}
If $\F$ is regular, then there is an open neighborhood of $\F$ in $\F(M)$ which consists of regular foliations with respect to the topology $\mathcal{O}^p_{w,\F^q(M)}$ for any $p \in \R_{\geq 1} \sqcup \{ \infty \}$. 
\end{lemma}

\begin{proof}
Fix any distance weight $w \colon [0,\infty] \to [0,1]$ and any number $p \in \R_{\geq 1}$. 
By $\lim_{r \to \infty} w(1/r) = 0$, there is a positive number $\varepsilon> 0$ such that $w(\varepsilon) < 1/2$. 
Fix a positive number $\delta := (\mathrm{vol}_g(B(x_0,\varepsilon))/\mathrm{vol}_g(M))^{1/p}/3$.
The regularity of $\F$ implies that $\rho_{\omega,\F} = 1$ is a constant function. 
Then $d_H(X_{w,\F}(x), \{ 0_x \} )= 1$ for any $x \in M$. 

Fix a singular foliation $\mathcal{G}$ with a point $x_0 \in \mathop{\mathrm{Sing}}(\G)$. 
Then $\{ 0_{x_0} \} = X_{\G}(x) = X_{w,\G}(x)$. 
Therefore $D^\infty_{w}(\F, \mathcal{G}) \geq d_H(X_{\F,x_0}, X_{w,\G}(x_0)) = d_H(X_{\F,x_0}, \{ 0_{x_0} \}) = 1 > 1/2$. 
This means that $B_{D^\infty_{w}}(\F,1/2)$ consists of regular foliations. 
We have that 
\[
\begin{split}
& D^p_{w}(\F, \mathcal{G}) \\
= & \left(\dfrac{1}{\mathrm{vol}_g(M)}  \int_M d_H(X_{w,\F}(x), X_{w,\G}(x))^p d \mathrm{vol}_g \right)^{1/p}\\
\geq & \left(\dfrac{1}{\mathrm{vol}_g(M)}  \int_{B(x_0,\varepsilon)} d_H(X_{w,\F}(x), X_{w,\G}(x))^p d \mathrm{vol}_g \right)^{1/p}\\
\geq & \left(\dfrac{1}{\mathrm{vol}_g(M)}  \int_{B(x_0,\varepsilon)} \left(\dfrac{1}{2}\right)^p d \mathrm{vol}_g \right)^{1/p}\\
\geq & \dfrac{1}{2} \left(\dfrac{1}{\mathrm{vol}_g(M)}  \int_{B(x_0,\varepsilon)} d \mathrm{vol}_g \right)^{1/p} =\dfrac{1}{2} \left(\dfrac{\mathrm{vol}_g(B(x_0,\varepsilon))}{\mathrm{vol}_g(M)} \right)^{1/p} > \delta 
\end{split}
\]
and so that $B_{D^p_{w}}(\F,\delta)$ consists of regular foliations. 
\end{proof}

The previous lemma implies the following openness. 

\begin{corollary}\label{cor:open_reg}
For any $p \in \R_{\geq 1} \sqcup \{ \infty \}$, the subset of regular foliations on $M$ is open in $\F(M)$ with respect to the topology $\mathcal{O}^p_{w,\F^q(M)}$.
\end{corollary}



\subsection{Persistence of singular points}

We describe the persistence of singular points.

\begin{lemma}\label{lem:singular_pt_nbd}
For any number $p \in \R_{\geq 1} \sqcup \{ \infty \}$, any distance weight $w \colon [0,\infty] \to [0,1]$, any number $\varepsilon>0$, and any point $x_0 \in \mathop{\mathrm{Sing}}(\F)$, there is a positive number $\delta > 0$ such that 
\[ 
B(x,\varepsilon) \cap \mathop{\mathrm{Sing}(\G)} \neq \emptyset
\]
for any $\mathcal{G} \in B_{D^p_{w}}(\F,\delta)$. 
\end{lemma}

\begin{proof}
Fix any distance weight $w \colon [0,\infty] \to [0,1]$, any $\varepsilon>0$, and any $x \in \mathop{\mathrm{Sing}}(\F)$. 
Take a large natural number $n \in \Z$ satisfying the following inequality: 
\[
w \left(\dfrac{(n-1)\varepsilon}{n} \right)  - w \left(\dfrac{\varepsilon}{n} \right) > 0
\]
Fix a positive number $\delta>0$ satisfying the following inequality: 
\[
 \delta < 
\left( w \left(\dfrac{(n-1)\varepsilon}{n} \right) - w \left(\dfrac{\varepsilon}{n} \right) \right) \min \left\{ 1, \left(\dfrac{\mathrm{vol}_g(B(x,\varepsilon/n))}{\mathrm{vol}_g(M)} \right)^{1/p} \right\}
\]

%
Assume that there is $\mathcal{G} \in B_{D^p_{w}}(\F,\delta)$ with $B(x,\varepsilon) \cap \mathop{\mathrm{Sing}(\G)} = \emptyset$. 
Then we have the following inequalities: 
\[
\begin{split}
D^\infty_{w}(\F, \mathcal{G}) & = \sup_{x \in M} d_H(X_{w,\F}(x), X_{w,\G}(x))  \\
& \geq \sup_{x \in B(x,\varepsilon/n)} d_H(X_{w,\F}(x), X_{w,\G}(x)) \\
& \geq \left( w \left(\dfrac{(n-1)\varepsilon}{n} \right) - w \left(\dfrac{\varepsilon}{n} \right) \right)  > \delta 
\end{split}
\]
and 
\[
\begin{split}
D^p_{w}(\F, \mathcal{G}) = & \left(\dfrac{1}{\mathrm{vol}_g(M)}  \int_M d_H(X_{w,\F}(x), X_{w,\G}(x))^p d \mathrm{vol}_g \right)^{1/p}\\
\geq  & \left(\dfrac{1}{\mathrm{vol}_g(M)}  \int_{B(x,\varepsilon/n)} d_H(X_{w,\F}(x), X_{w,\G}(x))^p d \mathrm{vol}_g \right)^{1/p}\\
\geq  & \left(\dfrac{1}{\mathrm{vol}_g(M)}  \int_{B(x,\varepsilon/n)} \left( w \left(\dfrac{(n-1)\varepsilon}{n} \right) - w \left(\dfrac{\varepsilon}{n} \right) \right)^p  d \mathrm{vol}_g \right)^{1/p}\\
= & \left( w \left(\dfrac{(n-1)\varepsilon}{n} \right) - w \left(\dfrac{\varepsilon}{n} \right) \right)  \left(\dfrac{\mathrm{vol}_g(B(x,\varepsilon/n))}{\mathrm{vol}_g(M)} \right)^{1/p}  > \delta 
\end{split}
\]
if $p \in \R_{\geq 1}$. 
The previous inequalities contradict to $\mathcal{G} \in B_{D^p_{w}}(\F,\delta)$. 
Thus $B(x,\varepsilon) \cap \mathop{\mathrm{Sing}(\G)} \neq \emptyset$ for any $\mathcal{G} \in B_{D^p_{w}}(\F,\delta)$. 
\end{proof}

\subsection{Descriptions near singular points}

We describe the behaviors near singular points to characterize the instability of some kinds of singular points. 

\begin{lemma}\label{lem:singular_pt_nbd03}
For any positive number $\varepsilon >0$, there is a positive number $\delta > 0$ satisfying the following inequality: 
\[ 
\sup_{x_0 \in \mathop{\mathrm{Sing}}(\F)} \sup_{x \in B(x_0,\delta)}d_H(X_{w,\F}(x), \{ 0_{x} \}) \leq \varepsilon
\]
\end{lemma}

\begin{proof}
Fix a positive number $\varepsilon >0$. 
Then $X_{w,\F}(x) = \{ 0_x \}$ for any $x \in \SF$. 
By $\lim_{r \to \infty} \sup_{x \in B(x_0,1/r)} \vert X_{w,\F}(x) \vert = 0$, there is a positive number $\delta \in (0, \varepsilon/2)$ such that $X_{w,\F}(x) \subseteq B(0_x, \varepsilon)$ for any $x_0 \in \mathop{\mathrm{Sing}}(\F)$ and any point $x \in B(x_0, \delta) $. 
On the other hand, we have $\{ 0_{x} \} \subset B(X_{w,\F}(x), \varepsilon)$ for any $x_0 \in \mathop{\mathrm{Sing}}(\F)$ and any $x \in B(x_0, \delta)$. 
This means that the assertion holds. 
\end{proof}

\begin{lemma}\label{lem:singular_pt_nbd02}
 Let $\mathcal{G} := \{ \{x \} \mid x \in M\}$ be the $0$-dimensional foliation on $M$. 
For any positive number $\varepsilon >0$, there is a positive number $\delta > 0$ satisfying the following inequality: 
\[ 
\sup_{x \in B(\mathop{\mathrm{Sing}}(\F),\delta)}d_H(X_{w,\F}(x), X_{w,\G}(x)) < \varepsilon
\]
Moreover, if $M = B(\mathop{\mathrm{Sing}}(\F),\delta)$, then $D^p_{w}(\F, \mathcal{G}) < \varepsilon$ for any $p \in \R_{\geq 1} \sqcup \{ \infty \}$. 
\end{lemma}

\begin{proof}
Fix a positive number $\varepsilon >0$. 
For any point $x \in M$, we obtain that $X_{w,\G}(x) = \{ 0_x \}$. 
Choose a positive number $\delta >0$ as in Lemma~\ref{lem:singular_pt_nbd03}. 
Fix a point $x \in B(\mathop{\mathrm{Sing}}(\F), \delta)$. 
If $x$ is singular, then $X_{w,\F}(x) = \{ 0_x \} = X_{w,\G}(x)$ and so $d_H(X_{w,\F}(x), X_{w,\G}(x)) = 0$. 
Thus we may assume that $x$ is not singular. 
By the choice of $\delta$, we obtain that $d_H(X_{w,\F}(x), X_{w,\G}(x)) = d_H(X_{w,\F}(x), \{ 0_x \}) < \varepsilon$. 

Suppose that $M = B(\mathop{\mathrm{Sing}}(\F),\delta)$. 
By defintion of $D^\infty_{w}$, we obtain that $D^\infty_{w}(\F, \mathcal{G}) < \varepsilon$. 
For any $p \in \R_{\geq 1}$, we have the following inequality: 
\[
\begin{split}
& D^p_{w}(\F, \mathcal{G}) \\
= & \left(\dfrac{1}{\mathrm{vol}_g(M)}  \int_M d_H(X_{w,\F}(x), X_{w,\G}(x))^p d \mathrm{vol}_g \right)^{1/p}\\
\leq & \left(\dfrac{1}{\mathrm{vol}_g(M)} \int_M \varepsilon^p d \mathrm{vol}_g \right)^{1/p} = \varepsilon \left(\dfrac{1}{\mathrm{vol}_g(M)} \mathrm{vol}_g(M) \right)^{1/p} = \varepsilon
\end{split}
\]
\end{proof}

\subsubsection{Transversality for singular foliations}

Recall that a subset $S$ is {\bf transverse} to a singular $C^1$-foliation $\F$ on a manifold $M$ if $T_x S + T_x \F = T_x M$ for any $x \in S$, where $T_x \F$ is the tangent space of a leaf of $\F$ at $x$. 
An arc is a {\bf transverse arc} to a codimension one foliation if it is transverse to the foliation. 

Although the transversality can be similarly defined for singular piecewise $C^1$-foliations on surfaces, we provide the details in Appendix~\ref{def:trans} due to its technical nature.

We define a transverse boundary of a trivial fibered chart as follows. 

\begin{definition}
Let $\F$ be a one-dimensional foliation on a manifold $M$. 
A connected transverse subset $T$ is a {\bf transverse boundary} of a trivial fibered chart $\varphi \colon V \to U_\tau \times U_\pitchfork = (-1,1) \times (-1,1)^{n-1}$ with respect to $\F$ if the trivial fibered chart $\varphi$ can be extended to a trivial fibered chart $\hat{\varphi} \colon \hat{V} \to (-1-\varepsilon,1+\varepsilon) \times U_\pitchfork$ for some positive number $\varepsilon$ such that $T = \hat{\varphi}^{-1}(\{-1,1\} \times U_\pitchfork)$. 
Then put $\partial_{\pitchfork} V := T$, $\partial_{\pitchfork}^-V := \hat{\varphi}^{-1}(\{-1\} \times U_\pitchfork)$, and $\partial_{\pitchfork}^+V := \hat{\varphi}^{-1}(\{1\} \times U_\pitchfork)$. 
\end{definition}

We define the first return to the transverse boundary of the trivial fibered chart. 

\begin{definition}
Let $V$ be a trivial fibered chart. 
The map $\partial_{\pitchfork}^-V \to \partial_{\pitchfork}^+V$ which maps the first returns of points in $\partial_{\pitchfork}^-V$ to points in $\partial_{\pitchfork}^+V$ along the leaf arcs in $V \sqcup \partial_{\pitchfork} V$ is called the {\bf first return with respect to the trivial fibered chart $V$}. 
\end{definition}

\subsubsection{Assumptions on the boundary}

In this paper, following the standard framework of foliation theory, we assume that every boundary component of a singular foliation is either transverse to the foliation or a union of leaves. 
More general settings  will be addressed in future work.

\subsection{Invariance of closed transverse subsets}

By definition of topologies, we have the following observations. 

\begin{lemma}\label{lem:trans}
Any closed transverse subset for a singular foliation on a Riemannian manifold $M$ is invariant under small perturbations with respect to $\mathcal{O}^\infty_{w,\F^q(M)}$. 
\end{lemma}

\begin{proof}
By definition of $\mathcal{O}^\infty_{w,\F^q(M)}$, the transversality of closed transverse subsets is invariant under small perturbations with respect to $\mathcal{O}^\infty_{w,\F^q(M)}$. 
\end{proof}

Notice that the previous lemma does not holds for any $p \in \R_{\geq 1}$ even for the restriction $\mathcal{O}^p_{\F_{\mathrm{reg}}^q(M)}$. 
More precisely, we have the following observation. 

\begin{lemma}
For any singular foliation $\F$ on a Riemannian manifold $M$ with a closed transverse subset $T$ such that $\dim M > \dim T$, any $p \in \R_{\geq 1}$, and any positive number $\varepsilon > 0$, 
there is a foliation $\mathcal{G}$ in $B_{D^p_{w}}(\F,\varepsilon)$ with $\mathop{\mathrm{Sing}}(\F) = \mathop{\mathrm{Sing}}(\G)$ which is obtain by a replacement of a small neighborhood $U$ of $T$ such that $\G$ is tangent to $T$ at a point in $T$. 
\end{lemma}

\begin{proof}
Fix a positive number $\varepsilon > 0$. 
Choose a neighborhood $U$ of $T$ with $\mathrm{vol}_g(U) < (\varepsilon/2)^p \mathrm{vol}_g(M)$. 
By $2 \geq \sup_{x \in M} d_H(X_{w,\F}(x), X_{w,\G}(x))$, the arbitrarily deformation on $U$ implies the following inequality for the resulting foliation $\mathcal{G}$: 
\[
\begin{split}
& D^p_{w}(\F, \mathcal{G}) \\
= & \left(\dfrac{1}{\mathrm{vol}_g(M)}  \int_M d_H(X_{w,\F}(x), X_{w,\G}(x))^p d \mathrm{vol}_g \right)^{1/p}\\
= & \left(\dfrac{1}{\mathrm{vol}_g(M)}  \int_U d_H(X_{w,\F}(x), X_{w,\G}(x))^p d \mathrm{vol}_g \right)^{1/p}\\
\leq & \left(\dfrac{1}{\mathrm{vol}_g(M)}  \int_U 2^p d \mathrm{vol}_g \right)^{1/p} = 2 \left(\dfrac{\mathrm{vol}_g(U)}{\mathrm{vol}_g(M)} \right)^{1/p} < \varepsilon
\end{split}
\]
Since the deformation on $U$ is an arbitrary replacement, we may assume that $\mathcal{G}$ has a tangency on the closed transverse subset. 
\end{proof}

Therefore, we use the topology $\mathcal{O}^\infty_{w,\F^q(M)}$ to discuss the structural stability of singular foliations on surfaces below.

\section{Singular foliations on surfaces}

In this section, we discuss singular foliations on surfaces. 

\subsection{Line fields}

We define line fields as follows.
\begin{definition} 
A singular foliation $\F$ is a {\bf line field}  if $\dim \F=1$. 
\end{definition} 
A line field is also called {\bf fields of line elements} (cf. \cite{Hopf1989diff}), {\bf director field} in nematic liquid crystals, and {\bf orientation field} in fingerprint recognition systems. 

From now on, let $\F$ be a line field on a surface $M$ unless otherwise stated.

\subsubsection{Types of leaves on codimension one singular foliations}

A leaf of a singular foliation on $M$ is {\bf proper} if it is an embedded submanifold in $M$, {\bf locally dense} if its closure has the nonempty interior, and {\bf exceptional} if it is neither proper nor locally dense. 

\subsubsection{Types of leaves on dimension one singular foliations}

As in the flow case (see \S~\ref{sec:flow} in Appendix), we define the following concepts. 

\begin{definition} 
A singular leaf of a singular foliation is a {\bf singularity} or a {\bf singular point} if it is a singleton (i.e. the dimension of the leaf is zero). 
\end{definition} 

\begin{definition} 
A leaf is {\bf periodic} if it is homeomorphic to a circle. 
\end{definition} 

Denote by $\bm{\mathop{\mathrm{Per}}(\F)}$ the union of periodic leaves. 

\begin{definition} 
A leaf is {\bf closed} if it is a singular or periodic leaf. 
\end{definition} 

Note that a leaf that is a closed interval on a closed disk is not a closed leaf. 
Moreover, any closed proper regular leaves on dimension one singular foliations on closed surfaces are periodic because any connected closed $1$-dimensional manifold is a circle.  
We define isolated property as follows. 

\begin{definition}
A singular point $x$ is {\bf isolated} if there is an open \nbd $U$ of $x$ with $U \cap \mathop{\mathrm{Sing}}(\F) = \{ x \}$. 
Then $U$ is called an isolated \nbd of $x$. 
\end{definition}




\subsubsection{The complement of the singular leaf set}

Notice that the compliment $\F - \mathcal S (\F)$ is a $C^q$ regular foliation, denoted by $\bm
{\F_{\mathrm{reg}}}$, of a surface $M - \mathop{\mathrm{Sing}}(\F)$, denoted by $\bm{M_{\mathrm{reg}}}$.  
We will show that the lift of the singular foliation with finitely many singular points to the ``$2$-fold covering'' of a compact surface is a flow, by similar arguments \cite[p.2 After the definition of orientability]{nikolaev2001foliations} and \cite[p.16, \S 2.3.4 and \S 2.3.5]{HH1986A}.

\subsubsection{Types of singular points}

We define the prongs, which are generalizations of multi-saddles, on the interiors of surfaces as follows. 

\begin{definition}
For any integer $k \in \Z_{\geq 2}$, an isolated singular point outside of the boundary $\partial S$ of $S$ is a {\bf $\bm{k}$-prong} or a {\bf $\bm{k}$-pronged singularity} if there is some complex coordinate $z$ on a neighborhood $U$ of the singularity such that $\mathcal F$ on $U$ corresponds to the set of the horizontal trajectories of the meromorphic quadratic differential $z^{k-2}dz^2$. 
\end{definition}

Notice that the complex coordinate in the previous definition need not be a chart on the surface. 
Note that a $3$-prong is also called a {\bf tripod}. 
We define the prongs on the boundaries of surfaces as follows. 

\begin{definition}
For any integer $k \in \Z_{\geq 2}$, an isolated singular point on $\partial S$ is a {\bf $\bm{k}$-prong} or a {\bf $\bm{k}$-pronged singularity} (cf. \cite{mosher1988tiling}) if there is a neighborhood $U$ of the singular point such that the double of $U$ is a neighborhood of a $2k-2$ prong (see Figure \ref{prong}). 
\end{definition}

\begin{figure}[t]
\begin{center}
\includegraphics[scale=0.275]{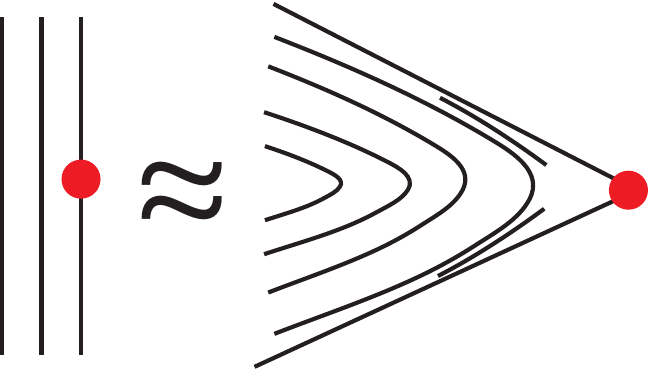}
\end{center} 
\caption{Hyperbolic sector.}
\label{hyperbolic_sector}
\end{figure}

\begin{figure}[t]
\begin{center}
\includegraphics[scale=0.25]{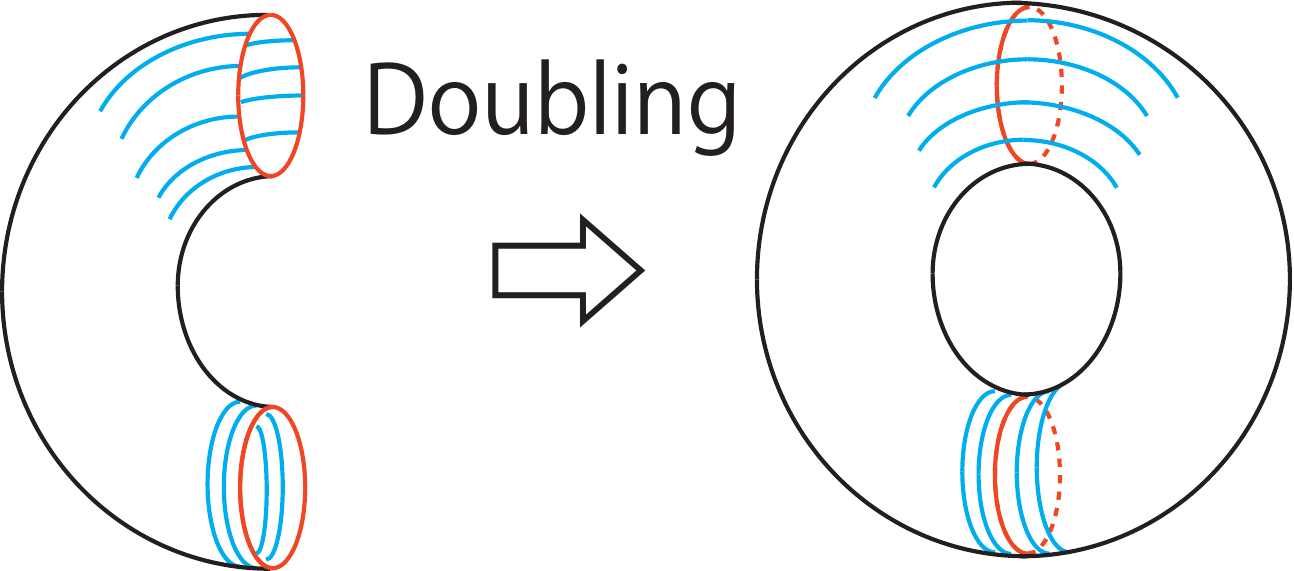}
\end{center} 
\caption{The lifted foliation on the double}
\label{pic01}
\end{figure}

Notice that an isolated singular point is a $k$-prong ($k \geq 2$) if and only if there is a neighborhood of it that consists of exactly $k$ hyperbolic sectors (see Figure \ref{hyperbolic_sector}).  

\begin{definition}
An isolated singular point outside of $\partial S$ is a $\bm{1}$-{\bf prong} (or {\bf thorn}) if the lift of some neighborhood $C$ of the singular point by $C \ni z\mapsto z^2 $ is a neighborhood of a $2$-prong. 
\end{definition}

\begin{definition}
An isolated singular point outside of $\partial S$ is a $\bm{0}$-{\bf prong} if it is a center. 
\end{definition}

%
%

\subsubsection{Singular foliations with pronged singular points}

Recall that \cite[Theorem~3]{cobo2010flows} implies that every singular point of a divergence-free flow (or more generally, a non-wandering flow (see Definition~\ref{def:nw})) with finitely many singular points on a compact surface is either a center or a multi-saddle. 
As an analogy to the case of incompressible (i.e. divergence-free) flow (see Definition~\ref{def:df}), we will introduce the concept of an “incompressible” singular foliation. 
To this end, we formulate a singular foliation with only prongs as follows, in analogy with the case where the foliation has only centers or multi-saddles.

\begin{definition}
A singular foliation is {\bf with pronged singular points} if any singular point is a prong and there is an open \nbd of the singular point set to which the restriction of the foliation is of $C^1$.  
\end{definition}

By definition of prong, any singular foliation with pronged singular points on compact surfaces has at most finitely many singular points. 
Notice that the set of orbits of a divergence-free flow with finitely many singular points on a compact surface is a singular foliation with pronged singular points, by applying \cite[Theorem 3]{cobo2010flows} to the lift of a singular foliation on the double of a compact surface.

\subsubsection{Fake prong}\label{sec:fake}
Recall that a singular point of a vector field on (resp. outside of) the boundary of a surface that consists of $2k$ (resp. $2k+1$) hyperbolic sectors is a multi-saddle, and that a singular point of a vector field on (resp. outside of) the boundary of a surface that consists of $2$ (resp. $1$) hyperbolic sectors is a fake multi-saddle. 
Notice that any fake multi-saddles are $2$-prongs. 
Similarly, we introduce fake prongs as follows. 

\begin{definition}
A $k$-prong of a singular foliation is {\bf fake prong} if $k =2$. 
\end{definition}

Notice that, by the assumption, each boundary component of a foliation with no degenerate
singular points on a surface is a circle that either is transverse to the foliation
or is a union of leaves. In particular, each center does not belong to the boundary.
We also note that any accumulation point of infinitely many non-degenerate singularities needs to be degenerate. 
In particular, if a singular foliation on a compact surface has no degenerate singularities, then there are at most finitely many non-degenerate singularities. 

\subsubsection{Tangencies of curves}

\begin{figure}
\begin{center}
\includegraphics[scale=0.175]{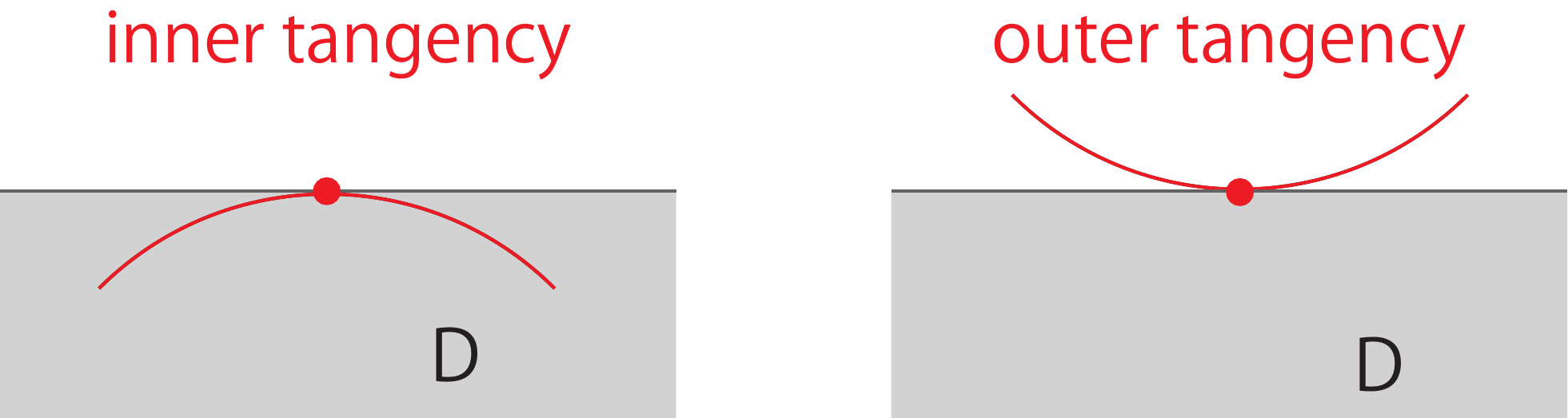}
\end{center}
\caption{Inner and outer tangencies.}
\label{fig:tangencies}
\end{figure}

A tangency $x$ of a curve $C$ to a loop bounding a closed disk $D$ is {\bf inner} (resp. {\bf outer}) if there is a small arc $I$ in $C$ containing $x$ such that the difference $I - \{ x \}$ is contained in the interior $\mathrm{int} D$ (resp. $(I - \{ x \}) \cap D = \emptyset$) as in left the (resp. right) on Figure~\ref{fig:tangencies}.

\subsubsection{Indices of prongs}

For a vector field on a surface $S$, recall that the {\bf index} of an isolated singular point outside of the boundary of $S$ is $(2 + n_{\mathrm{in}} - n_{\mathrm{out}})/2$, where $n_{\mathrm{in}}$ (resp. $n_{\mathrm{out}}$) is the number of inner (resp. outer) tangencies of a loop which is transverse at all but finitely many points and bounds an open disk containing the singular point. 
The {\bf index} of an isolated singular point $x$ on the boundary of $S$ is half of the index for the lift of the vector field at a lift of $x$ to the double of $S$. 
The index is known to be independent of the choice of loops that are transverse at all but finitely many points. 

%
Similarly, we define the index of an isolated singular point of a line field as follows. 

\begin{definition}
The {\bf index} $\mathop{\mathrm{ind}_{\F}} (x)$ of an isolated singular point $x$ of a line field on the boundary of $S$ is determined by $(2 + n_{\mathrm{in}} - n_{\mathrm{out}})/2$, where $n_{\mathrm{in}}$ (resp. $n_{\mathrm{out}}$) is the number of inner (resp. outer) tangencies of a loop which is transverse at all but finitely many points and bounds an open disk containing the singular point. 
\end{definition}

\begin{definition}
The {\bf index} $\mathop{\mathrm{ind}_{\F}} (x)$ of an isolated singular point $x$ on the boundary of $S$ is half of the index for the lift of the line field at a lift of $x$ to the double of $S$. 
\end{definition}

The index is also known to be independent of the choice of loops that are transverse at all but finitely many points. 
A singular point is orientable (resp. non-orientable) in the sense of Bronstein and Nikolaev \cite[p.6 Definition of (non-)orientability]{nikolaev2001foliations} if the index is (resp. is not) an integer. 
By definition of prongs, we have the following observation. 

\begin{lemma}
The following statements hold for any line field on a surface $S$ and for any integer $k \in \Z_{\geq 1}$: 
\\
{\rm(1)}  
The index of a $k$-prong 
 is $-(k -2)/2$. 
\end{lemma}

Denote by $\bm{\mathop{\mathrm{Sing}}_{-}(\F)}$ the set of singular points whose indices are negative.

\subsubsection{Indices of disks}\label{sec:index_disk}

We define the index of a disk with respect to a line field as follows. 

\begin{definition}
Let $D$ be a disk whose boundary is transverse to a line field $\F$ except finitely many tangencies. 
Then the {\bf index} $\mathop{\mathrm{ind}_{\F}} (D)$ of $D$ of the line field is determined by $(2 + n_{\mathrm{in}} - n_{\mathrm{out}})/2$, where $n_{\mathrm{in}}$ (resp. $n_{\mathrm{out}}$) is the number of inner (resp. outer) tangencies of the loop $\partial D$ which is transverse at all but finitely many points. 
\end{definition}

\begin{definition}
Let $D$ be a disk and $C \subset \partial D \cap \partial M$ a nonempty open arc such that $\partial D - C$ is transverse to a line field $\F$ except finitely many tangencies in $\mathop{\mathrm{int}}M$. 
The {\bf index} $\mathop{\mathrm{ind}_{\F}} (D, C)$ of the disk $D$ and the arc $C$ is half of the index of $\widetilde{D}$ of $\widetilde{\F}$, where $\widetilde{D}$ is the lift of $D$ to the double $\widetilde{M}$ of $M$ and $\widetilde{\F}$ is the lift of the line field $\F$ to $\widetilde{M}$. 
\end{definition}

Notice that the boundary $\partial \widetilde{D}$ in the previous definition is transverse to the line field $\F$ except finitely many tangencies.

\subsubsection{Separatrices of singular points}\label{sec:separatrix}

As a flow case, we define some kinds of separatrices as follows. 

\begin{definition}
A regular leaf is {\bf separatrix} if one of the connected components of the boundary of the leaf is a singular point. 
\end{definition}

\begin{definition}
A separatrix $L$ is a {\bf semi-prong separatrix} if one of the connected components of the difference $\overline{L} - L$ is a prong. 
\end{definition}

\begin{definition}
A separatrix $L$ is a {\bf prong separatrix} if the difference $\overline{L} - L$ consists of prongs. 
\end{definition}

The union of prongs and semi-prong separatrices is called the {\bf prong connection diagram}. 
A connected component of the prong connection diagram is called a {\bf prong connection}. 
Notice that the prong connection diagram and prong connections are analogous to the multi-saddle connection diagram and multi-saddle connections, respectively. 

\subsubsection{Orientability and transverse orientability for regular foliations on surfaces}

To define local orientability for prongs, recall the orientability of regular foliations as follows. 

\begin{definition}[cf. Definition~2.3.1(i) in \cite{HH1986A}]
A regular foliation on a surface is {\bf transversely orientable} if there is a foliated atlas whose restrictions of transition maps to the transverse directions are orientation-preserving. 
\end{definition}

\begin{definition}[cf. Definition~2.3.1(ii) in \cite{HH1986A}]
A regular foliation on a surface is {\bf (leafwise) orientable} if there is a foliated atlas whose restrictions of transition maps to the tangential directions are orientation-preserving. 
\end{definition}

Note that a $C^q$ regular foliation is leafwise (resp. transverse) orientable if and only if the tangent (resp. normal) bundle of the foliation is orientable. 
It is known that a regular foliation on a surface is orientable if and only if there is a flow on the surface whose set of orbits is the foliation (cf. \cite[Remark~2.3.2]{HH1986A}). 
Moreover, if the base surface of an orientable regular foliation is compact, then the flow is topologically equivalent to a $C^1$-foliation defined by Gutierrez's smoothing theorem~\cite{gutierrez1978structural}.

\subsubsection{Tangential orientable covering of a singular foliation on a surface}

Let $\F$ be a singular foliation on a surface $M$ with finitely many singular points and $\F_{\mathrm{reg}}$ the regular foliation on a surface $M_{\mathrm{reg}} = M - \mathop{\mathrm{Sing}}(\F)$. 
We define the tangent orientation branched covering $\pi \colon M^* \to M$ of $M$ with respect to $\F$ and the tangent orientation lift of $\F$ as follows.

As \cite[\S~2.3.5]{HH1986A}, we define a germ of the tangential orientation of the regular foliation $\F_{\mathrm{reg}}$ as follows. 

\begin{definition}
For any point $p \in M_{\mathrm{reg}}$, define an equivalence relation $\sim$ for fibered charts around $p$ as follows: 
Two fibered charts around $p$ are equivalent if the transition map locally preserves the orientations of the plaques of $\F$. 

Then each equivalent class is called a {\bf germ of the tangential orientation of $\bm{\F_{\mathrm{reg}}}$ at $\bm{p}$}. 
\end{definition}

Note that any fibered chart $(U,\varphi)$ of $\F_{\mathrm{reg}}$ determines exactly one germ $U^{\varphi,p,*}$ of the tangential orientation of $\F_{\mathrm{reg}}$ at a point $p \in U$. 
Notice that $\{ U^{\varphi,p,*} \mid p \in U \}$ is canonically homeomorphic to $U$. 
Moreover, for an nonempty open subset $U \in M_{\mathrm{reg}}$, the family $U^* := \{ U^{\varphi,p,*} \mid (U,\varphi) : \text{ fibered chart}, p \in U \}$ is the set of germs at points in $U$. 

For any foliated atlas $\{ (U_\lambda,\varphi_\lambda)\}_{\lambda \in \Lambda}$ of $\F_{\mathrm{reg}}$, the set $\{ U^*_\lambda \}_{\lambda \in \Lambda}$ form a basis for a topology on the set  $M_{\mathrm{reg}}^* := \bigcup_{\lambda \in \Lambda} U^*_\lambda$ of germs of the tangential orientation of $\F_{\mathrm{reg}}$. 
By construction, we have the following observation.

\begin{lemma}
The map $\pi \colon M_{\mathrm{reg}}^* \to M_{\mathrm{reg}}$, assigning to a germ of the tangential orientation of $\F_{\mathrm{reg}}$ at $p \in M - \mathop{\mathrm{Sing}}(\F)$ the point $p$, is a $2$-fold covering, called the {\bf tangent orientation branched covering} of $M_{\mathrm{reg}}$ with respect to $\F_{\mathrm{reg}}$. 
\end{lemma}

By \cite[Proposition~2.3.6]{HH1986A}, the lift $\F_{\mathrm{reg}}^*$ of $\F_{\mathrm{reg}}$ to the surface $M_{\mathrm{reg}}^*$ is a tangential orientable foliation. 
By definition of $M_{\mathrm{reg}}$ and by construction of the $2$-fold covering $\pi \colon M_{\mathrm{reg}}^* \to M_{\mathrm{reg}}$, adding finitely many points to $M_{\mathrm{reg}}^*$, we can extend $\pi$ to a map $\pi \colon M^* \to M$ from the resulting surface $M^*$. 
By construction, we have the following observations.

\begin{lemma}
The map $\pi \colon M^* \to M$ is the double covering of $M$ with ramification points contained in $\mathop{\mathrm{Sing}}(\F)$, called the {\bf tangent orientation branched covering} of $M$ with respect to $\F$. 
\end{lemma}

\begin{lemma}
The resulting decomposition on $M^*$ from $\F_{\mathrm{reg}}^*$ adding finitely many singular points is a singular foliation, denoted by $\bm{\F^*}$ and called the {\bf tangent orientation lift} of $\F$, such that the branched covering $\pi \colon M^* \to M$ maps leaves of $\F^*$ to leaves of $\F$.
\end{lemma}

\begin{figure}[t]
\begin{center}
\includegraphics[scale=0.75]{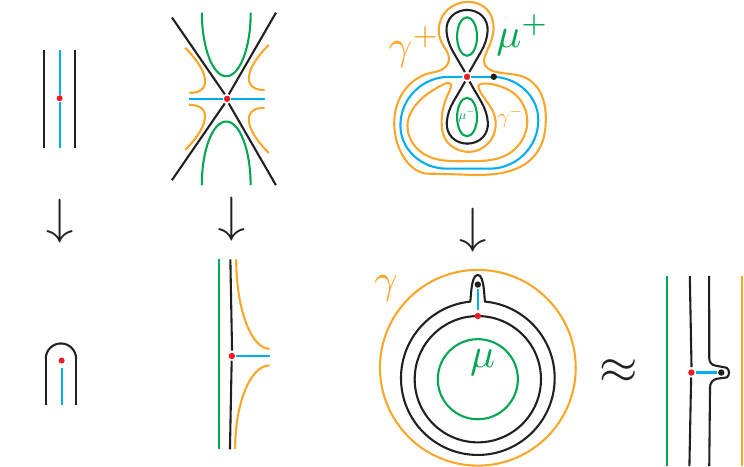}
\end{center} 
\caption{branched double covers with ramification points which are $1$-prongs and $3$-prongs.}
\label{bifurcation_prong_ex_01}
\end{figure}

We illustrate examples of local brunches as in Figure~\ref{bifurcation_prong_ex_01}. 
Notice that, for any $C^1$ singular foliation $\F$ on a surface $M$, the $2$-fold covering $\pi \colon M_{\mathrm{reg}}^* \to M_{\mathrm{reg}}$ can be identify with the canonical projection $\pi \colon T\F \cap T^1 M_{\mathrm{reg}} \to M_{\mathrm{reg}}$, where $T\F$ is the sub-bundle of the tangent bundle $TM$ consisting of tangent vectors to leaves of $\F$, and $T^1 M_{\mathrm{reg}}$ is the unit tangent bundle on $M_{\mathrm{reg}}$.

\subsubsection{Periodic annuli and periodic subsets}

As the flow case, we define the following concept. 

\begin{definition}
For a singular foliation $\F$, an annulus $A$ is a {\bf periodic annulus} if the restriction $\F|_A$ is topologically equivalent to a foliation $\{ \{t \} \times \mathbb{S}^1 \mid t \in I \}$ on an open annulus $I \times \mathbb{S}^1$ for some non-degenerate interval $I$. 
\end{definition}

More generally, we define the following concept. 

\begin{definition}
For a singular foliation $\F$, a subset $A$ is a {\bf periodic} if it consists of periodic orbits.
\end{definition}

\subsection{Comparison with similar concepts of singular foliations}

We compare our singular foliations with the following distinct definitions.

\subsubsection{Singular foliations with finitely many singular points on surfaces}

A triple $\F_{f,\theta} := (M, f, \theta)$ is $C^q$ foliation in the sense of Bronstein and Nikolaev~\cite[Definition~2.1]{bronstein1997peixoto} (cf. \cite{nikolaev2001foliations}) if $f \colon  \R \times M \to M$ is a $C^q$-flow on a compact surface $M$ and $\theta \colon M \to M$ is a smooth involution on $M$ satisfying the following conditions: 
\\
{\rm(i)} $f(-t,\theta(x)) = \theta (f(t,x))$ for any $(t,x) \in \R \times M$. 
\\
{\rm(ii)} The set of fixed points of $\theta$ consists of even number of singular points of $f$. 
\\
{\rm(iii)} If there are no fixed point of $\theta$, then $M$ is a disjoint union of two surfaces $M_1$ and $M_2$,  and the restriction $\theta|_{M_1} \colon M_1 \to M_2$ is a diffeomorphism. 

Notice that the first and second conditions in the original definition in \cite[Definition~2.1]{bronstein1997peixoto} are weaker than the above, but are required implicitly, because such conditions are required in \cite[Remark~2.1]{bronstein1997peixoto} and \cite[Proposition~0.1.1]{nikolaev2001foliations}. 
For instance, for a $z$-axial rotational flow $f \colon \R \times \mathbb{S}^2 \to \mathbb{S}^2$ with $2\pi$ period on the unit sphere $\mathbb{S}^2$ in $\R^3$ and the $\pi$-rotation $\theta:= f( \pi, \cdot) \colon \mathbb{S}^2 \to \mathbb{S}^2$, the triple $(\mathbb{S}^2, f, \theta)$ satisfies the original definition but not the first condition in the above definition. 
For an involution $\theta' \colon \mathbb{S}^2 \to \mathbb{S}^2$ defined by $(x,y,z) \to (-x,y,-z)$, which is a $\pi$-rotation with respect to the $y$-axis, the triple $(\mathbb{S}^2, f, \theta')$ satisfies the original definition but not the second condition in the above definition because the fixed point set of $\theta'$ consists of two fixed points $(0,1,0)$ and $(0,-1,0)$ which are not singular points of $f$. 

Notice that foliations in the sense of Bronstein and Nikolaev can not describe Whitehead moves and tripod-thorn bifurcations because Whitehead moves and tripod-thorn bifurcations change the base surfaces of foliations. 
In particular, the base surfaces of the foliations before and after such moves or bifurcations are not homeomorphic to each other in general. 
Thus, we employ a topology distinct from the uniform (Whitney) $C^q$ topology, as used in the Nikolaev book\cite{nikolaev2001foliations}, to describe the various deformations. 

Nikolaev listed all locally structurally stable non-orientable singularities for $C^q$ foliations in the sense of Bronstein and Nikolaev with respect to the uniform (Whitney) $C^q$ topology. 
In particular, if the flow is divergence-free, then such locally structurally stable non-orientable singularities are either thorns (i.e. $1$-prongs) or tripods (i.e. $3$-prongs). 

\section{Invariances of line fields}

In this section, we demonstrate the following invariance of sums of indices of singular leaves. 

\begin{proposition}\label{prop:inv_index}
Let $\F$ be a line field with finitely many singular points on a surface $M$. 
For any distance weight $w \colon [0,\infty] \to [0,1]$ and any $\varepsilon > 0$, there are a positive number $\delta > 0$ and pairwise disjoint disks $D_1, \ldots , D_k$ which is open and whose diameters are less than $\varepsilon$ and each of which contain exactly one singular point $p_i \in D_i$ satisfying the following relations: 
\[
\bigcup_{\mathcal{G} \in B_{D^\infty_{w}}(\F,\delta)} \mathop{\mathrm{Sing}}(\G) \subset \bigsqcup_{i=1}^k D_i
\]
\[
\mathop{\mathrm{ind}_{\F}} (p_i) = \sum_{p \in \mathop{\mathrm{Sing}}(\F) \cap U_i} \mathop{\mathrm{ind}_{\F}} (p) = \sum_{p \in \mathop{\mathrm{Sing}}(\G) \cap U_i} \mathop{\mathrm{ind}_{\G}} (p)
\]  
\end{proposition}

To show this, we have the following statements. 

\subsection{Persistence of leaf arcs between open transverse arcs}

We demonstrate the existence of leaf arcs between open transverse arcs under small perturbations as follows. 

\begin{lemma}\label{lem:persist_leaf_arcs}
Let $\F$ be a line field on a $(q+1)$-dimensional Riemannian manifold $M$ and fix any distance weight $w \colon [0,\infty] \to [0,1]$. 
Let $U \subset M$ be a closed fibered chart with a homeomorphism $h \colon U \times [0,1] \times [-1,1]^q$ such that, for any $y \in [-1,1]^q$, the inverse image $h^{-1}([0,1] \times \{ y \})$ is a $C^1$-leaf arc. 
For any distance weight $w \colon [0,\infty] \to [0,1]$ and any $\varepsilon > 0$, there is a positive number $\delta > 0$ such that the following statements hold: 
\\
{\rm(1)} The inverse images $T_0 := h^{-1}(\{ 0 \} \times [-1,1]^q) \subset \partial U$ and $T_1 := h^{-1}(\{ 1 \} \times [-1,1]^q) \subset \partial U$ are transverse to any line field $\mathcal{G} \in B_{D^\infty_{w}}(\F,\delta)$.
\\
{\rm(2)} For any $\mathcal{G} \in B_{D^\infty_{w}}(\F,\delta)$, there are a point $p^\G \in T_1 \cap B_\varepsilon (h^{-1}((0,1)))$ and a leaf arc $C^{\mathcal{G}} \subset U$ with respect to $\mathcal{G}$ between $h^{-1}(0)$ and $p^\G$ such that $d_H (C^{\mathcal{G}}, h^{-1}([0,1] \times \{ 0 \})) < \varepsilon$. 
\\
{\rm(3)} $U \cap \bigcup_{\G \in B_{D^\infty_{w}}(\F,\delta)} \mathop{\mathrm{Sing}}(\G) = \emptyset$.
\end{lemma}

\begin{proof}
Fix any $\varepsilon > 0$. 
Put $p^0 := h^{-1}(0)$, $C^\F := h^{-1}([0,1] \times \{ 0 \})$, $\mathop{\mathrm{int}} T_0 := h^{-1}(\{ 0 \} \times (-1,1)^q)$, and $\mathop{\mathrm{int}} T_ 1 := h^{-1}(\{ 1 \} \times (-1,1)^q)$. 
From Lemma~\ref{lem:open_reg_pt}, the compactness of $U$ implies that there is a positive number $\delta_0 > 0$ such that $U \cap \bigcup_{\G \in B_{D^\infty_{w}}(\F,\delta_0)} \mathop{\mathrm{Sing}}(\G) = \emptyset$.
By Lemma~\ref{lem:trans}, there is a positive number $\delta_1 \in (0, \delta_0)$ such that the arcs $T_0$ and $T_1$ are transverse to any line field $\G \in B_{D^\infty_{w}}(\F,\delta_1)$. 
Put $D_0 := \min_{p \in C_\F} \min_{p' \in \partial U - (\mathrm{int}T_0 \sqcup \mathrm{int}T_1)} d(p,p') > 0$. 

For any $\G \in B_{D^\infty_{w}}(\F,\delta_1)$, there uniquely exist a leaf arc $C^{\G} \subset U$ and a point $p^\G \in \partial U - \{ p_0 \}$ such that $C^\G \setminus \mathop{\mathrm{int}}U = \{ p^0, p^\G \}$. 
By the compactness of $C^\F$, the uniform continuity of $X_{w,\F}(x)\vert_{C^\F}$ implies that the sequenece $(\G_n)_{n \in \Z_{>0}}$ of any elements $\G_n \in B_{D^\infty_{w}}(\F,\delta_1/n)$ satisfies $\lim_{n \to \infty} d_H(C^{\G_n}, C^\F) = 0$ and $\lim_{n \to \infty} p^{\G_n} = h^{-1}((0,1))$. 
Therefore, there is a number $N > 0$ such that $\delta := \delta_1/N$ is desired. 
\end{proof}

The previous lemma implies the following observation. 

\begin{corollary}\label{cor:persist_leaf_arcs}
Let $\F$ be a singular foliation on a surface $M$ and $\gamma$ a loop which is a finite union of closed transverse arcs $T_0, \ldots , T_k$ and pairwise disjoint $C^1$ closed leaf arcs $C_0, \ldots , C_k$ with $\gamma = \left( \bigsqcup_{i = 0}^k T_i \right) \sqcup \left( \bigsqcup_{i = 0}^k \mathop{\mathrm{int}} C_i \right)$. 
For any distance weight $w \colon [0,\infty] \to [0,1]$ and any $\varepsilon > 0$, there are a positive number $\delta > 0$ and pairwise disjoint closed transverse arcs $\widetilde{T_0}, \ldots , \widetilde{T_k}$ with $T_i \subset \mathop{\mathrm{int}} \widetilde{T_i}$ for any $i \in \{0, \ldots , k \}$ such that the following statements hold for any line field $\mathcal{G} \in B_{D^\infty_{w}}(\F,\delta)$: 
\\
{\rm(1)} There is an open annular \nbd $A$ of $\gamma$ with $A \cap \bigcup_{\G \in B_{D^\infty_{w}}(\F,\delta)} \mathop{\mathrm{Sing}}(\G) = \emptyset$ and there is a loop $\gamma^\G \subset A$ which is a finite union of closed transverse arcs $T^\G_0 \subset \widetilde{T_0}, \ldots , T^\G_k \subset \widetilde{T_k}$ with respect to $\G$ and pairwise disjoint closed leaf arcs $C^\G_0, \ldots , C^\G_k$ with respect to $\G$ such that $\gamma^\G = \left( \bigsqcup_{i = 0}^k T^\G_i \right) \sqcup \left( \bigsqcup_{i = 0}^k \mathop{\mathrm{int}} C^\G_i \right)$, $d_H (T_i, T^{\G}_i) < \varepsilon$, and $d_H (C_i, C^{\G}_i) < \varepsilon$ for any $i \in \{0, 1, \ldots , k \}$.
\\
{\rm(2)} Suppose that there is a closed disk $D'^\F$ whose boundary is $\gamma$. 
Then we can make $\gamma^\G$ the boundary of a closed disk $D'^\G$ and there are disks $D^\F$ and $D^\G$ which are arbitraliry near to $D'^\F$ and $D'^\G$ respectively with respect to the Hausdorff distance such that $\mathop{\mathrm{ind}_\F}(D^\F) = (2-k)/2 = \mathop{\mathrm{ind}_\G}(D^\G)$ and $\partial D^\F \cup \partial D^\G \subset A$. 
Moreover, if $D'^\F \cap (\mathop{\mathrm{Sing}}(\F) \cup \mathop{\mathrm{Sing}}(\G)) < \infty$, then $\sum_{p \in \mathop{\mathrm{Sing}}(\F) \cap D'^\F } \mathop{\mathrm{ind}_{\F}} (p) = \sum_{p \in \mathop{\mathrm{Sing}}(\G) \cap D'^\F } \mathop{\mathrm{ind}_{\G}} (p)$.  
\end{corollary}

\begin{proof}
Applying Lemma~\ref{lem:persist_leaf_arcs} finitely many times, one can obtain assertion (1). 
Deforming any leaf arc $C_i$ between transverse arcs $C'_i$ into an arc transverse to a singular foliation except one tangency fixing the boundary of the leaf arc, the resulting loop $\left( \bigsqcup_{i = 0}^k T_i \right) \sqcup \left( \bigsqcup_{i = 0}^k \mathop{\mathrm{int}} C'_i \right)$ bounds a disk $D^\F$ near $D$ and is transverse to $\F$ except finitely many tangencies. 
By the same argument for $\gamma^\G$, one can obtain a loop a disk $D^\G$ near $D$ and is transverse to $\G$ except finitely many tangencies such that the number of inner (resp. outer) tangencies of $\F$ equals to one of inner (resp. outer) tangences of $\F$. 
This means that the assertion $(2)$ holds. 
\end{proof}

Similarly, we have the following statement for disks tangent to the boundary of a surface. 

\begin{corollary}\label{cor:persist_leaf_arcs_02}
Let $\F$ be a singular foliation tangent to the boundary $\partial M$ on a surface $M$ and $\gamma$ a closed arc between distinct two points on $\partial M$ which is a finite union of closed transverse arcs $T_0, \ldots , T_k$ and pairwise disjoint $C^1$ closed leaf arcs $C_0, \ldots , C_k$ with $\gamma = \left( \bigsqcup_{i = 0}^k T_i \right) \sqcup \left( \bigsqcup_{i = 0}^k \mathop{\mathrm{int}} C_i \right)$. 
For any distance weight $w \colon [0,\infty] \to [0,1]$ and any $\varepsilon > 0$, there is a positive number $\delta > 0$ such that the following statements hold for any line field $\mathcal{G} \in B_{D^\infty_{w}}(\F,\delta)$: 
\\
{\rm(1)} There is an open \nbd $A$ of $\gamma$ with $A \cap \bigcup_{\G \in B_{D^\infty_{w}}(\F,\delta)} \mathop{\mathrm{Sing}}(\G) = \emptyset$ such that the lift of $A$ to the double of $M$ is an open annulus, and there is a closed arc $\gamma^\G \subset A$ between distinct two points on $\partial M$ which is a finite union of closed transverse arcs $T^\G_0, \ldots , T^\G_k$ with respect to $\G$ and pairwise disjoint closed leaf arcs $C^\G_0, \ldots , C^\G_k$ with respect to $\G$ such that $\gamma^\G = \left( \bigsqcup_{i = 0}^k T^\G_i \right) \sqcup \left( \bigsqcup_{i = 0}^k \mathop{\mathrm{int}} C^\G_i \right)$, $d_H (T_i, T^{\G}_i) < \varepsilon$, and $d_H (C_i, C^{\G}_i) < \varepsilon$ for any $i \in \{0, 1, \ldots , k \}$.
\\
{\rm(2)} Suppose that there are an open arc $I \subset \partial M$ and a disk $D'^\F$ whose boundary is $\gamma \sqcup I$. 
Then there is a disk $D^\F$ arbitraliry near to $D'^\F$ whose lift to the double of $M$ is a disk, and we can make $\gamma^\G$ the boundary of a closed disk $D'^\G$ arbitraliry near to $D'^\F$ whose lift to the double of $M$ is a disk and whose boundary is $\gamma^\G \sqcup I^\G$ with $\mathop{\mathrm{ind}_\F}(D^\F) = (2-k)/2 = \mathop{\mathrm{ind}_\G}(D^\G)$ and $(\partial D^\F \cup \partial D^\G) \setminus I \subset A$, where $I^\G \subset \partial M$ is some open arc.
Moreover, if $D'^\F \cap (\mathop{\mathrm{Sing}}(\F) \cup \mathop{\mathrm{Sing}}(\G)) < \infty$, then $\sum_{p \in \mathop{\mathrm{Sing}}(\F) \cap D'^\F } \mathop{\mathrm{ind}_{\F}} (p) = \sum_{p \in \mathop{\mathrm{Sing}}(\G) \cap D'^\F } \mathop{\mathrm{ind}_{\G}} (p)$.  
\end{corollary}

The above piecewise $C^1$ properties imply the following invariance of indices on small \nbds of prongs. 

\begin{lemma}\label{lem:inv_index_prong}
Let $\F$ be a singular foliation $\F$ on a surface $M$ and $w \colon [0,\infty] \to [0,1]$ a distance weight. 
For any $\varepsilon > 0$ and any prong $x \in \mathop{\mathrm{int}} M$, there are a positive number $\delta > 0$ and a small open disk $U \subset \mathop{\mathrm{int}} M$ which is an isolated \nbd of $x$ and whose boundary is a finite union of closed transverse arcs $T_0, \ldots , T_k$ and pairwise disjoint $C^1$ closed leaf arcs $C_0, \ldots , C_k$ with $\partial U = \left( \bigsqcup_{i = 0}^k T_i \right) \sqcup \left( \bigsqcup_{i = 0}^k \mathop{\mathrm{int}} C_i \right)$ satisfying the following statements for any line field $\mathcal{G} \in B_{D^\infty_{w}}(\F,\delta)$: 
\\
{\rm(1)} $\mathop{\mathrm{ind}_{\F}} (x) = \sum_{p \in \mathop{\mathrm{Sing}}(\F) \cap U} \mathop{\mathrm{ind}_{\F}} (p) = \sum_{p \in \mathop{\mathrm{Sing}}(\G) \cap U} \mathop{\mathrm{ind}_{\G}} (p)$.  
\\
{\rm(2)} There is an open annular \nbd $A$ of $\partial U$ with 
\[
A \cap \bigcup_{\G \in B_{D^\infty_{w}}(\F,\delta)} \mathop{\mathrm{Sing}}(\G) = \emptyset
\]
such that the union $U \cup A$ is an open disk, and there is a loop $\gamma^\G \subset A$ which is a finite union of closed transverse arcs $T^\G_0, \ldots , T^\G_k$ with respect to $\G$ and pairwise disjoint closed leaf arcs $C^\G_0, \ldots , C^\G_k$ with respect to $\G$ such that $\gamma^\G = \left( \bigsqcup_{i = 0}^k T^\G_i \right) \sqcup \left( \bigsqcup_{i = 0}^k \mathop{\mathrm{int}} C^\G_i \right)$, $d_H (T_i, T^{\G}_i) < \varepsilon$, and $d_H (C_i, C^{\G}_i) < \varepsilon$ for any $i \in \{0, 1, \ldots , k \}$.
\end{lemma}

\begin{proof}
Since any prong is locally $C^1$ with respect to a singular foliation, there is a small open disk $U \subset \mathop{\mathrm{int}} M$ which is an isolated \nbd of $x$ and whose boundary is a finite union of closed transverse arcs $T_0, \ldots , T_k$ and pairwise disjoint $C^1$ closed leaf arcs $C_0, \ldots , C_k$ with $\gamma = \left( \bigsqcup_{i = 0}^k T_i \right) \sqcup \left( \bigsqcup_{i = 0}^k \mathop{\mathrm{int}} C_i \right)$. 
By Corollary~\ref{cor:persist_leaf_arcs}, there is a positive number $\delta >0$ as in Corollary~\ref{cor:persist_leaf_arcs}. 
Then the assertion holds.
\end{proof}

Similarly, we have the following statement. 

\begin{lemma}\label{lem:inv_index_prong_bd}
Let $\F$ be a singular foliation $\F$ tangent to the boundary $\partial M$ on a surface $M$ and $w \colon [0,\infty] \to [0,1]$ a distance weight. 
For any $\varepsilon > 0$ and any prong $x \in \partial M$, there are a positive number $\delta > 0$ and a small open isolated \nbd $U$ of $x$ whose lift to the double of $M$ is an open disk, there are an open arc $I \subset \partial M$, and there is a closed arc $\gamma = \partial U - I$ between distinct two points on $\partial M$ which is a finite union of closed transverse arcs $T_0, \ldots , T_k$ and pairwise disjoint $C^1$ closed leaf arcs $C_0, \ldots , C_k$ with $\gamma = \left( \bigsqcup_{i = 0}^k T_i \right) \sqcup \left( \bigsqcup_{i = 0}^k \mathop{\mathrm{int}} C_i \right)$ satisfying the following statements for any line field $\mathcal{G} \in B_{D^\infty_{w}}(\F,\delta)$: 
\\
{\rm(1)} $\mathop{\mathrm{ind}_{\F}} (x) = \sum_{p \in \mathop{\mathrm{Sing}}(\F) \cap U} \mathop{\mathrm{ind}_{\F}} (p) = \sum_{p \in \mathop{\mathrm{Sing}}(\G) \cap U} \mathop{\mathrm{ind}_{\G}} (p)$.  
\\
{\rm(2)} There is an open \nbd $A$ of $\gamma$ with 
\[
A \cap \bigcup_{\G \in B_{D^\infty_{w}}(\F,\delta)} \mathop{\mathrm{Sing}}(\G) = \emptyset
\]
such that the union $U \cup A$ is a disk which is open, and there is a cllosed arc $\gamma^\G \subset A$ between distinct two points on $\partial M$ which is a finite union of closed transverse arcs $T^\G_0, \ldots , T^\G_k$ with respect to $\G$ and pairwise disjoint closed leaf arcs $C^\G_0, \ldots , C^\G_k$ with respect to $\G$ such that $\gamma^\G = \left( \bigsqcup_{i = 0}^k T^\G_i \right) \sqcup \left( \bigsqcup_{i = 0}^k \mathop{\mathrm{int}} C^\G_i \right)$, $d_H (T_i, T^{\G}_i) < \varepsilon$, and $d_H (C_i, C^{\G}_i) < \varepsilon$ for any $i \in \{0, 1, \ldots , k \}$.
\end{lemma}

Now, the Proof of Proposition~\ref{prop:inv_index} is completed as follows.

\begin{proof}[Proof of Proposition~\ref{prop:inv_index}]
Fix any distance weight $w \colon [0,\infty] \to [0,1]$ and any number $\varepsilon > 0$. 
By Corollary~\ref{cor:nonexistence_sing}, for the number $\varepsilon$, there is a positive number $\delta_0>0$ with $(M - B_\varepsilon(\mathop{\mathrm{Sing}}(\F))) \cap \bigcup_{\G \in B_{D^\infty_{w}}(\F,\delta_0)} \mathop{\mathrm{Sing}}(\G) = \emptyset$.
From Lemma~\ref{lem:inv_index_prong} and Lemma~\ref{lem:inv_index_prong_bd}, the finiteness of singular points of $\F$ implies that there are a positive number $\delta_1 > 0$ and pairwise disjoint disks $D_1, \ldots , D_k$ which is open and whose diameters are less than $\varepsilon$ and each of which contain exactly one singular point as in these lemmas. 
Put $\delta := \min \{ \delta_0, \delta_1 \}$. 
Then $\delta$ and $D_1, \ldots , D_k$ are desired. 
\end{proof}

\section{Non-existence of structural stability of singular foliations on surfaces}

This section shows the non-existence of structural stability of singular foliations on compact surfaces as shown in Theorem~\ref{non-existence_structural_stability}. 
%
%
We recall the following concept. 

\begin{definition}
A point $x'$ is the {\bf first return} of a point $x$ to a transverse $T$ along a leaf arc $C \colon [a,b] \to \F(x)$ if the restriction $C|_{(a,b)}$ is a simple curve in $\F(x) \setminus T$ with $C(a) = x \in T$ and $C(b) = x' \in T$.  
For such first returns along $C$, the existence of fibered charts implies a first return map $T' \to T''$ along $C$ between sub-arcs $T'$ and $T''$ of the transverse $T$, called a {\bf first return map} {\rm(}or holonomy map{\rm)} along $C$ to $T$. 
\end{definition}

To demonstrate Theorem~\ref{non-existence_structural_stability}, we show several lemmas. 

\begin{lemma}\label{lem:no_sing_leaf_int-}
Fix any $p \in \R_{\geq 1} \sqcup \{ \infty \}$ and let $M$ be a surface. 
Suppose that $\mathrm{vol}_g(M) < \infty$ if $p \in \R_{\geq 1}$. 
Then any singular foliation on $M$ whose union of singular leaves has a non-empty interior is not structurally stable in $\mathcal{F}^q(M)$ with respect to the topology $\mathcal{O}^p_{w,\F^q(M)}$. 
\end{lemma}

\begin{proof}
Let $\mathcal{G}$ be a singular foliation on a surface $M$ whose union of singular leaves has a non-empty interior. 
Put $S := \mathrm{int}(\mathop{\mathrm{Sing}}(\mathcal{G})) \neq \emptyset$. 
Then, any leaves of the restriction of $\mathcal{G}$ to the surface $S$ are singletons. 
Fix any distance weight $w \colon [0,\infty] \to [0,1]$ and any positive number $\varepsilon >0$. 
By Lemma~\ref{lem:singular_pt_nbd02}, there is a positive number $\delta > 0$ satisfying the following inequality for any singular foliation $\F$ on $M$:
\[ 
\sup_{x \in B(\mathop{\mathrm{Sing}}(\F),\delta)} d_H(X_{w,\F}(x), X_{w,\G}(x)) < \varepsilon
\]
Replacing the singular foliation $\mathcal{G}\vert_S$ into new one by using a bump function, we can choose a singular foliation $\F \in \mathcal{F}^q(M)$ with $X_\F \vert_{M -S} = X_\mathcal{G}\vert_{M -S}$, $\mathrm{int}(\mathop{\mathrm{Sing}}(\mathcal{\F})) = \emptyset$, and $S \subseteq B(\mathop{\mathrm{Sing}}(\F),\delta)$.
Then $\partial S = \partial \mathop{\mathrm{Sing}}(\G) \subseteq \mathop{\mathrm{Sing}}(\F)$. 
Therefore, we have the following inequalities: 
\[
\begin{split}
D^\infty_{w}(\F, \mathcal{G}) &= \sup_{x \in M} d_H(X_{w,\F}(x), X_{w,\G}(x))\\
&= \sup_{x \in B(\mathop{\mathrm{Sing}}(\F),\delta)}d_H(X_{w,\F}(x), X_{w,\G}(x)) < \varepsilon
\end{split}
\]
\[
\begin{split}
&D^p_{w}(\F, \mathcal{G})\\ 
=& \left(\dfrac{1}{\mathrm{vol}_g(M)}  \int_M d_H(X_{w,\F}(x), X_{w,\G}(x))^p d \mathrm{vol}_g \right)^{1/p}\\
=& \left(\dfrac{1}{\mathrm{vol}_g(M)}  \int_{B(\mathop{\mathrm{Sing}}(\F),\delta)} d_H(X_{w,\F}(x), X_{w,\G}(x))^p d \mathrm{vol}_g \right)^{1/p}\\
< & \left(\dfrac{1}{\mathrm{vol}_g(M)}  \int_{B(\mathop{\mathrm{Sing}}(\F),\delta)} \varepsilon^p d \mathrm{vol}_g \right)^{1/p} = \varepsilon \left(\dfrac{\mathrm{vol}_g(B(\mathop{\mathrm{Sing}}(\F),\delta))}{\mathrm{vol}_g(M)}\right)^{1/p} \leq \varepsilon
\end{split}
\]
By $\mathrm{int}(\mathop{\mathrm{Sing}}(\mathcal{\mathcal{G}})) \neq \emptyset$ and $\mathrm{int}(\mathop{\mathrm{Sing}}(\mathcal{\F})) = \emptyset$, the singular foliations $\mathcal{G}$ and $\F$ are not topologically equivalent. 
This means that $\mathcal{G}$ is not structurally stable in $\mathcal{F}^q(M)$ with respect to the topology $\mathcal{O}^p_{w,\F^q(M)}$. 
\end{proof}

\begin{lemma}\label{lem:no_sing_leaf_int}
Fix any $p \in \R_{\geq 1} \sqcup \{ \infty \}$ and let $M$ be a surface. 
Suppose that $\mathrm{vol}_g(M) < \infty$ if $p \in \R_{\geq 1}$. 
Then any singular foliation with singular points on a surface is not structurally stable in $\mathcal{F}^q(M)$ with respect to the topology $\mathcal{O}^p_{w,\F^q(M)}$. 
\end{lemma}

\begin{proof}
Let $\F$ be a singular foliation with singular points on a surface $M$. 
Lemma~\ref{lem:no_sing_leaf_int-} implies that we may assume that $\mathrm{int}(\mathop{\mathrm{Sing}}(\mathcal{\F})) = \emptyset$. 
Fix any distance weight $w \colon [0,\infty] \to [0,1]$ and any positive number $\varepsilon >0$. 

\begin{claim}\label{claim:07}
There is a singular foliation $\G$ with $D^p_{w}(\F, \mathcal{G}) < \varepsilon$ and $\mathrm{int}(\mathop{\mathrm{Sing}}(\mathcal{\mathcal{G}})) \neq \emptyset$. 
\end{claim}

\begin{proof}[Proof of Claim~\ref{claim:07}]
By Lemma~\ref{lem:singular_pt_nbd02}, there is a positive number $\delta > 0$ satisfying the following inequality for any singular foliation $\G_1$ on $M$ with $B(\mathop{\mathrm{Sing}}(\F),\delta) \subseteq \mathop{\mathrm{Sing}}(\G_1)$ and for any singular foliation $\G_2$ on $M$:
\[
\sup_{x \in B(\mathop{\mathrm{Sing}}(\F),\delta)}d_H(X_{w,\F}(x), X_{w,\G_1}(x)) < \varepsilon/2
\]
\[
\sup_{x \in B(\mathop{\mathrm{Sing}}(\F),\delta)}d_H(X_{w,\G_2}(x), X_{w,\G_1}(x)) < \varepsilon/2
\]
Therefore, we have 
\[
\sup_{x \in B(\mathop{\mathrm{Sing}}(\F),\delta)}d_H(X_{w,\F}(x), X_{w,\G}(x)) < \varepsilon
\]
for any singular foliation $\G$ on $M$. 
Taking $\delta >0$ small if necessary,  by the uniform continuity of $w$, we may assume that $\vert w(r_1) - w(r_2) \vert < \varepsilon$ for any $r_1, r_2 \in [0,\infty)$ with $\vert r_1 - r_2 \vert \leq \delta$. 
Choose a singular foliation $\G$ on $M$ with 
\[
B(\mathop{\mathrm{Sing}}(\F),\delta/2) \subseteq \mathop{\mathrm{Sing}}(\G) \subset B(\mathop{\mathrm{Sing}}(\F), \delta)
\]
and $X_\F \vert_{M - B(\mathop{\mathrm{Sing}}(\F),\delta)} = X_\mathcal{G}\vert_{M - B(\mathop{\mathrm{Sing}}(\F),\delta)}$. 
Then 
\[
\sup_{x \in B(M - \mathop{\mathrm{Sing}}(\F),\delta))} (d_g(x, \mathop{\mathrm{Sing}}(\F)) - d_g(x, \mathop{\mathrm{Sing}}(\G))) \in [\delta/2,\delta]
\]
and 
$\mathrm{int}(\mathop{\mathrm{Sing}}(\mathcal{\mathcal{G}})) \neq \emptyset$. 
The equality $X_\F \vert_{M - B(\mathop{\mathrm{Sing}}(\F),\delta)} = X_\mathcal{G}\vert_{M - B(\mathop{\mathrm{Sing}}(\F),\delta)}$ implies the following inequality: 
\[
\begin{split}
& \sup_{x \in B(M - \mathop{\mathrm{Sing}}(\F),\delta))} d_H(X_{w,\F}(x), X_{w,\G}(x))
\\
= & \sup_{x \in B(M - \mathop{\mathrm{Sing}}(\F),\delta))} (\rho_{w,\F}(x) - \rho_{w,\G}(x))
\\
= & \sup_{x \in B(M - \mathop{\mathrm{Sing}}(\F),\delta))} (w(d_g(x, \mathop{\mathrm{Sing}}(\F))) - w(d_g(x, \mathop{\mathrm{Sing}}(\G))))
\\ < & \varepsilon
\end{split}
\]
This implies $D^\infty_{w}(\F, \mathcal{G}) < \varepsilon$. 
On the other hand, 
\[
\begin{split}
&D^p_{w}(\F, \mathcal{G})\\ 
=& \left(\dfrac{1}{\mathrm{vol}_g(M)}  \int_{B(\mathop{\mathrm{Sing}}(\F),\delta)} d_H(X_{w,\F}(x), X_{w,\G}(x))^p d \mathrm{vol}_g \right)^{1/p}\\
< & \left(\dfrac{1}{\mathrm{vol}_g(M)}  \int_{B(\mathop{\mathrm{Sing}}(\F),\delta)} \varepsilon^p d \mathrm{vol}_g \right)^{1/p} = \varepsilon \left(\dfrac{\mathrm{vol}_g(B(\mathop{\mathrm{Sing}}(\F),\delta))}{\mathrm{vol}_g(M)}\right)^{1/p} \leq \varepsilon
\end{split}
\]
if $p < \infty$. 
\end{proof}

By $\mathrm{int}(\mathop{\mathrm{Sing}}(\mathcal{\mathcal{G}})) \neq \emptyset$ and $\mathrm{int}(\mathop{\mathrm{Sing}}(\mathcal{\F})) = \emptyset$, the singular foliations $\mathcal{G}$ and $\F$ are not topologically equivalent. 
This means that $\mathcal{F}$ is not structurally stable in $\mathcal{F}^q(M)$ with respect to the topology $\mathcal{O}^p_{w,\F^q(M)}$ for any $p \in \R_{\geq 1} \sqcup \{ \infty \}$. 
\end{proof}

\begin{lemma}\label{lem:limit_cycle}
The boundary of any connected component of the interior $\mathrm{int} (\mathop{\mathrm{Per}}(\mathcal{\F}))$ of the union of periodic leaves of a regular foliation $\F$ on a closed surface consists of periodic leaves. 
\end{lemma}

\begin{proof}
Let $C$ be a connected component of the interior $\mathrm{int} (\mathop{\mathrm{Per}}(\mathcal{\F}))$ of the union of periodic leaves of a regular foliation $\F$ on a closed surface $M$. 
Take a boundary component $\partial_C \subseteq \partial C$ of $C$. 

It suffices to show that $\partial_C$ is a periodic leaf. 
Indeed, by the existence of fibered charts, the compactness of $C$ implies that there are finitely many fibered charts $U_1, \ldots , U_k$ centered at points in $C$ such that $C \subset \bigcup_{i=1}^k U_i$. 
By the compactness of $C$, there is an invariant periodic annulus $A \subseteq C \cap \bigcup_{i=1}^k U_i$ such that $\partial_C$ is a boundary component of $A$.  
Then the intersections $A \cap U_i$ consist of closed leaf arcs. 
Fix a transverse closed arc $T \subset \partial_C \sqcup A \subseteq \overline{A}$. 
Then the first return map on the interval $T \cap A = T \setminus \partial_C$ with any direction is identical. 
By the existence of fibered charts that cover $C$, the boundary component $\partial_C$ is a finite union of closed leaf arcs. 
From the existence of the identical first return map on the interval $T \cap A$, the finite union $\partial_C$ is a loop, which is a periodic leaf.  
\end{proof}

\begin{lemma}\label{lem:nss_cpt_leaf}
Any singular foliations on a surface $M$ whose unions of periodic leaves have nonempty interiors are not structurally stable in $\mathcal{F}^q(M)$ with respect to the topology $\mathcal{O}^p_{w,\F^q(M)}$ for any $r \in \Z_{>0}$ and any $p \in \R_{\geq 1} \sqcup \{ \infty \}$. 
\end{lemma}

\begin{proof}
Fix any positive number $\varepsilon > 0$ and any $p \in \R_{\geq 1} \sqcup \{ \infty \}$.  
Put a positive number $N_g \in (0,1]$ with $N_g = \max \{1, (\mathrm{vol}_g(\mathrm{int} (\mathop{\mathrm{Per}}))/\mathrm{int} (\mathop{\mathrm{Per}}(\mathcal{\F})))^{1/p} \}$ if $p < \infty$. 
Let $\F$ be a singular foliation on a surface $M$ whose union of periodic leaves has a nonempty interior. 
Fix any connected component $C$ of $\mathrm{int} (\mathop{\mathrm{Per}}(\mathcal{\F}))$. 
By Lemma~\ref{lem:limit_cycle}, the closure $\overline{C}$ is an invariant closed annulus. 
Fix a fibered chart $B$ in $C$ whose saturation is $C$. 
Replacing the fibered chart $B$ into a fibered chart as in Figure~\ref{flowbox_perturbation02}, we can obtain the resulting foliation $\mathcal{G}_C$ satisfying the following inequality: 
\begin{figure}[t]
\begin{center}
\includegraphics[scale=0.75]{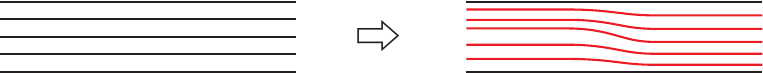}
\end{center} 
\caption{Replacement of a trivial fibered chart.}
\label{flowbox_perturbation02}
\end{figure}
\[
\begin{split}
D^\infty_{w}(\F, \mathcal{G}_C) &= \sup_{x \in M} d_H(X_{w,\F}(x), X_{w,\G_C}(x))\\
&= \sup_{x \in C} d_H(X_{w,\F}(x), X_{w,\G_C}(x)) 
< \dfrac{\varepsilon}{N_g} \leq \varepsilon
\end{split}
\]
Applying the above replacements for all connected components of $\mathrm{int} (\mathop{\mathrm{Per}}(\mathcal{\F}))$, the resulting foliation $\mathcal{G}$ satisfies the following inequality: 
\[
\begin{split}
D^\infty_{w}(\F, \mathcal{G}) &= \sup_{x \in M} d_H(X_{w,\F}(x), C_{w,\mathcal{G}}(x,r,\alpha))\\
&= \sup_{x \in \mathop{\mathrm{Per}}(\mathcal{\F})} d_H(X_{w,\F}(x), C_{w,\mathcal{G}}(x,r,\alpha)) < \dfrac{\varepsilon}{N_g} \leq \varepsilon 
\end{split}
\]
Moreover, we have the following inequalities for any $p < \infty$: 
\[
\begin{split}
&D^p_{w}(\F, \mathcal{G})\\ 
=& \left(\dfrac{1}{\mathrm{vol}_g(M)}  \int_{M} d_H(X_{w,\F}(x), X_{w,\G}(x))^p d \mathrm{vol}_g \right)^{1/p}\\
=& \left(\dfrac{1}{\mathrm{vol}_g(M)}  \int_{\mathrm{int} (\mathop{\mathrm{Per}}(\mathcal{\F}))} d_H(X_{w,\F}(x), X_{w,\G}(x))^p d \mathrm{vol}_g \right)^{1/p}\\
< & \left(\dfrac{1}{\mathrm{vol}_g(M)}  \int_{\mathrm{int} (\mathop{\mathrm{Per}}(\mathcal{\F}))} \left(\dfrac{\varepsilon}{N_g}\right)^p d \mathrm{vol}_g \right)^{1/p}\\
=& \dfrac{\varepsilon}{N_g} \left(\dfrac{1}{\mathrm{vol}_g(M)}  \int_{\mathrm{int} (\mathop{\mathrm{Per}}(\mathcal{\F}))}  d \mathrm{vol}_g \right)^{1/p} = \dfrac{\varepsilon}{N_g} \left(\dfrac{\mathrm{vol}_g(\mathrm{int} (\mathop{\mathrm{Per}}(\mathcal{\F})))}{\mathrm{vol}_g(M)}  \right)^{1/p} \leq \varepsilon
\end{split}
\]
\end{proof}

\begin{lemma}\label{lem:nss_cpt_leaf03}
Fix any $p \in \R_{\geq 1} \sqcup \{ \infty \}$ and let $M$ be a surface. 
Suppose that $\mathrm{vol}_g(M) < \infty$ if $p \in \R_{\geq 1}$. 
Then any singular foliations with periodic leaves on $M$ are not structurally stable in $\mathcal{F}^q(M)$ with respect to the topology $\mathcal{O}^p_{w,\F^q(M)}$ for any $r \in \Z_{>0}$ and any $p \in \R_{\geq 1} \sqcup \{ \infty \}$. 
\end{lemma}

\begin{proof}
Fix any distance weight $w \colon [0,\infty] \to [0,1]$, any $\varepsilon>0$, and any singular foliation $\F$ with periodic leaves on a surface $M$. 
Choose a small positive number $\delta \in (0,1)$. 
Let $L$ be a periodic leaf of $\F$. 
By Lemma~\ref{lem:nss_cpt_leaf}, we may assume that $\mathrm{int}(\mathop{\mathrm{Per}}(\F)) = \emptyset$.  
Choose a fibered chart $U$ centered at a point $p \in L$ such that $U$ (resp. $p$) can be identified with a rectangle $U_\tau \times U_\pitchfork = [-R_1, R_1] \times [-R_2, R_2]$ (resp. $0$), where $R_1, R_2 \in (0, 1)$ are positive numbers with $R_1 > 2 \delta R_2$. 
Take a transverse open arc $T \subseteq \{ 0 \} \times [-R_2, R_2] \subset U$ intersecting $x$ and a small open sub-arc $T' \subset T$  intersecting $x$ such that the first return map $f_{\F} \colon T' \to T$ along the one-dimensional leaf $L$ to $T$ with the negative direction for $U_\tau$ is well-defined. 
Then there is a sequence $(p_n)_{n \in \Z_{\geq 0}}$ with $p = \lim_{n \to \infty} p_n  = \lim_{n \to \infty} f_\F(p_n)$. 
Replacing $[-R_2, R_2]$ with a small sub-arc, we may assume that $R_2 < \delta R_1/6$. 
%
Since $\F$ is of piecewise $C^q$, so is the return map $f_\F$. 
Fix a smooth bump function $\varphi \colon [-R_1, R_1] \to [0,1]$ whose support is $[0,R_1/2]$ such that $\int_{0}^{R_1} \varphi(x_1) dx_1 = R_1/3$. 
Moreover, choose a smooth bump function $\psi \colon [-R_2, R_2] \to [0,1]$ whose support is $[-R_2/2,R_2/2]$ such that $\psi^{-1}(1) = [-R_2/3,R_2/3]$. 
Define a function $f \colon U = [-R_1, R_1] \times [-R_2, R_2] \to [-R_2, R_2]$ as follows: 
\[
f(x_1,x_2) :=  x_2 + \psi(x_2) \dfrac{3}{R_1}(f_\F(x_2)-x_2) \int_{0}^{x_1} \varphi(x_1) dx_1
\]
Since $\varphi^{-1}(0) = [-R_1, 0] \sqcup [R_1/2, R_1]$, we have that $f(x_1,x_2) = x_2$ for any $x_1 \in [-R_1, 0]$, and that 
\[
f(x_1,x_2) = x_2 +  \psi(x_2) (f_\F(x_2)-x_2)
\]
for any $x_1 \in [R_1/2, R_1]$ because $\int_{0}^{R_1} \varphi(x_1) dx_1 = R_1/3$. 
Therefore, for any $x_1 \in [R_1/2, R_1]$, we obtain 
\[
f(x_1,x_2) = x_2 + (f_\F(x_2)-x_2) = f_\F(x_2)
\]
for any $x_2 \in \psi^{-1}(1) = [-R_2/3,R_2/3]$, 
and 
\[
f(x_1,x_2) = x_2 
\]
for any $x_2 \in \psi^{-1}(1) = [-R_2,-R_2/2] \sqcup [R_2/2,R_2]$. 
%
The annulus $A := U - ([0,R_1/2] \times [-R_2/2,R_2/2])$ is a color of the loop $\partial U$ such that the images $f(A \cap ([-R_1, R_1] \times \{ x_2 \} ))$ for $x_2 \in [-R_2, R_2]$ are disjoint union of at most two horizontal arcs (i.e. parallel to the $x_1$-axis). 
By replacing the $\F$-fibered chart $U$ with a fibered chart whose plaques are $l(x_2) := \{(x_1, f(x_1,x_2)) \mid x_1 \in  [-R_1, R_1]\}$ for any $x_2 \in [-R_2, R_2]$, the resulting partition $\mathcal{G}$ is a piecewise $C^q$ singular foliations whose union of periodic leaves has a nonempty interior.

Since the union $\mathop{\mathrm{Per}}(\F)$ of periodic leaves of $\F$ has the empty interior, the foliations $\F$ and $\mathcal{G}$ are not topologically equivalent. 

\begin{claim}\label{claim:08}
We may assume that $D^p_{w}(\F, \mathcal{G}) < \varepsilon$ for any $p \in \R_{\geq 1} \sqcup \{ \infty \}$. 
\end{claim}

\begin{proof}[Proof of Claim~\ref{claim:08}]
By construction, we have $\mathop{\mathrm{Sing}}(\F) = \mathop{\mathrm{Sing}}(\G)$ and so $\rho_{w,\F} = \rho_{w,\G}$.
This implies the following inequalities: 
\[
\begin{split}
D^\infty_{w}(\F, \mathcal{G}) & = \sup_{x \in U} d_H(X_{w,\F}(x), X_{w,\G}(x))
\\
& = \sup_{x \in U} d_H(\rho_{w,\F}(x) X_{\F}(x), \rho_{w,\G}(x) X_{\G}(x))
\\
& = \sup_{x \in U} d_H(\rho_{w,\F}(x) X_{\F}(x), \rho_{w,\F}(x) X_{\G}(x))
\\
& = \sup_{x \in U} \rho_{w,\F}(x) d_H( X_{\F}(x), X_{\G}(x))
\\
& \leq \sup_{x \in U} d_H( X_{\F}(x), X_{\G}(x))
\end{split}
\]

We have the following inequalities: 
\[
\dfrac{\partial f(x_1,x_2)}{\partial x_1} = \psi(x_2) \dfrac{3}{R_1}(f_\F(x_2)-x_2) \varphi(x_1) \leq  \dfrac{3}{R_1}(f_\F(x_2)-x_2) \leq \dfrac{6R_2}{R_1} < \delta
\]
If we make the given number $\delta > 0$ arbitrarily small, then the following maximal angle
\[
\max_{x \in U} \angle (T_{x} \F(x),T_{x} \mathcal{G}(x)) 
\]
is arbitrarily small, where $\angle(v,w)$ is the angle between tangent vectors $v$ and $w$. 
Thus the superior 
\[
\sup_{x \in U} d_H( X_{\F}(x), X_{\G}(x))
\]
becomes arbitrarily small. 
Replacing $\delta > 0$ with a small positive number if necessary, we may assume that the following inequality holds: 
\[
\sup_{x \in U} d_H( X_{\F}(x), X_{\G}(x)) < \varepsilon
\]
This means that $D^\infty_{w}(\F, \mathcal{G}) < \varepsilon$ and so
\[
\begin{split}
D^p_{w}(\F, \mathcal{G}) =& \left(\dfrac{1}{\mathrm{vol}_g(M)}  \int_{U} d_H(X_{w,\F}(x), X_{w,\G}(x))^p d \mathrm{vol}_g \right)^{1/p}\\
\leq & \left(\dfrac{1}{\mathrm{vol}_g(M)}  \int_{U} \varepsilon^p d \mathrm{vol}_g \right)^{1/p} = \varepsilon \left(\dfrac{\mathrm{vol}_g(U)}{\mathrm{vol}_g(M)}\right)^{1/p} < \varepsilon
\end{split}
\]
if $p < \infty$. 
\end{proof}

This completes the proof. 
\end{proof}

\begin{lemma}\label{lem:nss_cpt_leaf02}
Any singular foliations with non-proper leaves on a surface $M$ are not structurally stable in $\mathcal{F}^q(M)$ with respect to the topology $\mathcal{O}^p_{w,\F^q(M)}$ for any $r \in \Z_{>0}$ and any $p \in \R_{\geq 1} \sqcup \{ \infty \}$. 
\end{lemma}

\begin{proof}
Fix any positive number $\varepsilon \in (0,1)$. 
Let $\F$ be a singular foliation with non-proper leaves on a surface $M$. 
By Lemma~\ref{lem:limit_cycle} and Lemma~\ref{lem:nss_cpt_leaf03}, we may assume that $\F$ has no closed leaves. 
This implies that $\F$ is a regular foliation. 
Take a non-proper leaf $L \in \F$ and an arbitrarily small transverse open arc $T$ intersecting a point $p \in L$. 
The non-closedness and non-proper property imply that $\vert T \cap L \vert = \infty$. 
As the proof of the previous lemma, replacing a fibered chart $U$ centered at $p$ with a fibered chart, we can obtain the resulting regular foliation $\mathcal{G}$ with respect to which first return to $T$ of $p$ is the same point $p$ such that 
\[
\max_{x \in U} \angle (T_{x} \F(x),T_{x} \mathcal{G}(x)) < \delta
\]
for arbitrarily small number $\delta >0$. 
By the smallness of $\delta >0$, we may assume that $\sup_{x \in U} d_H( X_{\F}(x), X_{\G}(x)) < \varepsilon$. 
This implies
\[
\begin{split}
D^\infty_{w}(\F, \mathcal{G}) 
&= \sup_{x \in U} d_H( X_{\F}(x), X_{\G}(x))< \varepsilon
\end{split}
\]
if $p = \infty$, and 
\[
\begin{split}
D^p_{w}(\F, \mathcal{G}) =& \left(\dfrac{1}{\mathrm{vol}_g(M)}  \int_{U} d_H(X_{w,\F}(x), X_{w,\G}(x))^p d \mathrm{vol}_g \right)^{1/p}\\
\leq & \left(\dfrac{1}{\mathrm{vol}_g(M)}  \int_{U} \varepsilon^p d \mathrm{vol}_g \right)^{1/p} = \varepsilon \left(\dfrac{\mathrm{vol}_g(U)}{\mathrm{vol}_g(M)}\right)^{1/p} < \varepsilon
\end{split}
\]
if $p < \infty$. 
Since the first return to $T$ of $p$ is the same point $p$ with respect to $\mathcal{G}$, the resulting foliation $\mathcal{G}$ has a periodic leaf. 
Since $\F$ has no periodic leaves, the foliation $\F$ is not topologically equivalent to $\mathcal{G} \in B_{D^p_{w}}(\F,\varepsilon)$, which implies the absence of structural stability.
\end{proof}

\begin{lemma}\label{lem:non_per_proper_leaf}
For any $r \in \Z_{>0}$ and any $p \in \R_{\geq 1} \sqcup \{ \infty \}$, every singular foliation on a surface $M$ which is structurally stable in $\mathcal{F}^q(M)$ with respect to the topology $\mathcal{O}^p_{w,\F^q(M)}$ consists of non-closed proper leaves. 
\end{lemma}

\begin{proof}
Fix any $r \in \Z_{>0}$ and any $p \in \R_{\geq 1} \sqcup \{ \infty \}$. 
Let $\F$ be a singular foliation on a surface $M$. 
Assume that $\F$ is structurally stable in $\mathcal{F}^q(M)$ with respect to the topology $\mathcal{O}^p_{w,\F^q(M)}$. 
Lemma~\ref{lem:no_sing_leaf_int} and Lemma~\ref{lem:nss_cpt_leaf} imply that $\F$ has no closed leaves. 
By Lemma~\ref{lem:nss_cpt_leaf02}, the foliation $\F$ has no non-proper leaves.
Therefore, the regular foliation $\F$ consists of non-closed proper leaves. 
\end{proof}

We demonstrate the main result in this section. 


\begin{proof}{Proof of Theorem~\ref{non-existence_structural_stability}}
Fix any $q \in \Z_{>0}$ and any $p \in \R_{\geq 1} \sqcup \{ \infty \}$. 
Let $\F$ be a singular foliation tangent to the boundary on a compact surface $M$. 
Assume that $\F$ is structurally stable in $\mathcal{F}^q(M)$ with respect to the topology $\mathcal{O}^p_{w,\F^q(M)}$. 
Lemma~\ref{lem:non_per_proper_leaf} implies that $\F$ consists of non-closed proper leaves. 
It is known that the foliation $\F$ contains a minimal set (cf. \cite[Proposition~4.1.8, p.46]{HH1986A}). 
The minimal set is either a periodic leaf or the closure of a non-proper leaf (cf. \cite[Proposition~4.1.7, p.46]{HH1986A}), which contradicts that $\F$ has neither proper leaves nor non-proper leaves.  
\end{proof}

\subsection{Structural stability for foliation transverses to the boundaries on compact surfaces}

We observe that there is a foliation transverse to the boundary on a compact surface $M$ is structurally stable in $\mathcal{F}^q(M)$ with respect to the topology $\mathcal{O}^\infty_{w,\F^q(M)}$ for any $r \in \Z_{>0}$ as follows. 

\begin{lemma}
Let $\mathbb{A} := [0,1] \times (\R/\Z)$ be a closed annulus and $\F := \{ [0,1] \times [\theta] \mid [\theta] \in \R/\Z \}$ a foliation transverse to the boundary on $\mathbb{A}$. 
Then $\F$ is structurally stable in $\mathcal{F}^q(\mathbb{A})$ with respect to the topology $\mathcal{O}^\infty_{w,\F^q(M)}$. 
\end{lemma}

\begin{proof} 
By definition of $\F$, the loops $\{x \} \times (\R/\Z)$ for any $x \in [0,1]$ are closed transversal for $\F$. 
Fix any $\delta \in (0, 1/2)$. 
Since $D^\infty_{w}(\F, \mathcal{G}) \geq 1$ for any singular foliation $\mathcal{G} \in \mathcal{F}^q(\mathbb{A})$ with singular points, the ball $B_{D^\infty_{w}}(\F,\delta)$ contains no singular foliations.
Fix any $\mathcal{G} \in B_{D^\infty_{w}}(\F,\delta)$. 

\begin{claim}\label{claim:09}
The loops $\{x \} \times (\R/\Z)$ for any $x \in [0,1]$ are closed transversal for $\G$
\end{claim}

\begin{proof}[Proof of Claim~\ref{claim:09}]
Notice that $D^\infty_{w}(\F, \mathcal{G}) = \sup_{x \in \A} d_H( X_{\F}(x), X_{\G}(x))< \delta < 1$ and so that $\sup_{x \in \A} \angle (T_x \F(x),T_x \mathcal{G}(x)) < \pi/3 < \pi/2$. 
This means that the loops $\{x \} \times (\R/\Z)$ for any $x \in [0,1]$ are closed transversal for $\G$.
\end{proof}

By the previous claim, the regular foliations $\F$ and $\mathcal{G}$ are topologically equivalent. 
Thus $\F$ is structurally stable in $\mathcal{F}^q(M)$ with respect to the topology $\mathcal{O}^\infty_{w,\F^q(M)}$ for any $r \in \Z_{>0}$. 
\end{proof}

The previous lemma implies that the tangential condition in Theorem~\ref{non-existence_structural_stability} is necessary. 

\section{Divergence-free foliations with pronged singular points}

Liouville’s theorem implies that any smooth vector field on a surface is divergence-free if and only if it is area-preserving (cf. \cite[Lemma~2.2]{yokoyama2021ham}).
Moreover, area-preserving flows on surfaces are also known as locally Hamiltonian flows or equivalently multi-valued Hamiltonian flows. 
This means that any divergence-free flow on a surface has a multi-valued Hamiltonian and so a transversely invariant measure, and so that the set of orbits of any divergence-free non-singular flow on a surface is a Riemannian foliation. 
Thus we introduce a divergence-free singular foliation as below. 

\subsection{Riemannian foliations on surfaces}
To introduce such a divergence-free property for singular foliations, recall the concept of Riemannian foliations on surfaces as follows. 
A codimension $k$ regular foliation $\F$ is {\bf Riemannian} if there is a maximal atlas $\{ (V_\alpha, \varphi_\alpha \colon V_\alpha \to U_{\tau, \alpha} \times U_{\pitchfork, \alpha}) \}_{\alpha \in \Lambda}$ of fibered charts of $\F$ such that the transition maps
\[
\varphi_{\alpha \beta} := \varphi_\alpha \circ \varphi_\beta^{-1} \vert_{\varphi_\beta^{-1}(V_\alpha \cap V_\beta)} \colon \varphi_\beta^{-1}(V_\alpha \cap V_\beta) \to \varphi_\alpha^{-1}(V_\alpha \cap V_\beta)
\]
are of the form
\[
\varphi_{\alpha \beta}(x_\tau, x_\pitchfork) = (\varphi_{\alpha \beta, \tau}(x_\tau, x_\pitchfork), \varphi_{\alpha \beta, \pitchfork}(x_\pitchfork))
\]
with $\varphi_{\alpha \beta, \pitchfork} \in O(k, \R)$, where $O(k, \R)$ is the orthogonal group in dimension $k$. 

\subsection{Divergence-free singular foliations}

We introduce a divergence-free property for singular foliations with finitely many singular points as follows.

\begin{definition}
A singular foliation $\F$ with finitely many singular points on a surface $M$ is {\bf divergence-free} if the regular foliation $\F_{\mathrm{reg}}$ on the surface $M_{\mathrm{reg}} = M - \mathop{\mathrm{Sing}}(\F)$ is a Riemannian regular foliation. 
\end{definition}

We have the following characterization of the divergence-free property. 

\begin{lemma}\label{lem:ch:levelable}
The following statements are equivalent for any singular foliation $\F$ with finitely many singular points on a surface $M$:
\\
{\rm(1)} The singular foliation $\F$ is divergence-free. 
\\
{\rm(2)} The tangent orientation lift $\F^*$ of $\F$ on the tangent orientation branched covering $M^*$ of $M$ is topologically equivalent to the set of orbits of the flow generated by a smooth divergence-free vector field.
\end{lemma}

\begin{proof}
Suppose that $\F$ is divergence-free. 
By construction, the tangent orientation branched covering $\pi \colon M_{\mathrm{reg}}^* \to M_{\mathrm{reg}}$ of $M_{\mathrm{reg}}$ with respect to the Riemannian regular foliation $\F_{\mathrm{reg}}$ is an isometric $2$-fold covering. 
Then the restriction $\F^*|_{M^* - \mathop{\mathrm{Sing}}(\F^*)}$ of the tangent orientation lift $\F^*$ is also a Riemannian regular foliation, because so is $\F_{\mathrm{reg}}$. 
By \cite[Remark~2.3.2]{HH1986A}, there is a flow on the surface $M^* - \mathop{\mathrm{Sing}}(\F^*)$ whose set of orbits is the restriction $\F^*|_{M^* - \mathop{\mathrm{Sing}}(\F^*)}$. 
Therefore, there is a flow $f_{\F^*}$ on the compact surface $M^*$ whose set of orbits is the foliation $\F^*$. 
The existence of the transverse Riemannian structure implies that the complement $M^* - \mathop{\mathrm{Sing}}(\F^*)$ is non-wandering with respect to the flow $f_{\F^*}$. 
Since the set of non-wandering points of $f_{\F^*}$ is closed, the flow $f_{\F^*}$ on $M^*$ is non-wandering.
\cite[Theorem~A]{yokoyama2021ham} implies that $f_{\F^*}$ is generated by a smooth divergence-free vector field up to topological equivalence.

Conversely, suppose that the lift $\F^*$ of $\F$ on $M^*$ is topologically equivalent to the set of orbits of the flow generated by a smooth divergence-free vector field.
By the divergence-free property of the vector field, the restriction $\F^*|_{M^* - \mathop{\mathrm{Sing}}(\F^*)}$ is a Riemannian regular foliation. 
By construction of the tangent orientation branched covering $\pi \colon M_{\mathrm{reg}}^* \to M_{\mathrm{reg}}$, the set $\F_{\mathrm{reg}}$ is the tangent orientation lift. 
\end{proof}


We characterize singular points of divergence-free singular foliations on surfaces.

\begin{lemma}\label{lem:type_singular_points}
Any singular points of a divergence-free singular foliation on a surface are either centers or prongs. 
\end{lemma}

\begin{proof}
Let $\F$ be a divergence-free singular foliation on a surface $M$. 
By Lemma~\ref{lem:ch:levelable}, the tangent orientation lift $\F^*$ of $\F$ on the tangent orientation branched covering $M^*$ of $M$ is generated by a smooth divergence-free vector field up to topological equivalence. 
Since any divergence-free vector field is locally a Hamiltonian vector field, the lift $\widetilde{\F^*}$ of the foliation $\F^*$ on the universal covering $\widetilde{M^*}$ of the surface $M^*$ is generated by a smooth Hamiltonian vector field up to topological equivalence. 
\cite[Lemma~3.1]{yokoyama2024topological} implies that any singular points of $\widetilde{\F^*}$ are either centers or prongs, and so are those of $\F^*$. 
By construction of $\F^*$, any singular point of the singular foliation $\F$ is either a center or a prong.
\end{proof}


%
%
%
%
%

From now on, we assume that all leaves are of piecewise $C^1$, to define the distance on leaves of a foliation induced by a Riemannian metric of a surface. 

\subsection{Divergence singular foliations with pronged singular points}\label{section:div-free}



Notice that there are contractible closed transversals for divergence-free foliations as on the bottom in Figure~\ref{disk_2prongs}. 

\begin{definition}
A leaf $L$ is {\bf transversely non-orientable} if, for any piont $x \in L$ and for any small transverse open arc $T'$ intersecting $x$, there are a transverse open arc $T \subseteq T'$ containing $x$ and a trivial fibered chart $V$ whose transverse boundary is $T - \{ x \}$ such that the first return map with respect to $V$ is orientation-reversing. 
\end{definition}

\begin{definition}
A leaf is {\bf transversely orientable} if it is not transversely non-orientable {\rm(i.e.} there is a small transverse open arc $T'$ intersecting a point $x \in L$ such that, for any trivial fibered chart $V$ whose transverse boundary is $T - \{ x \}$ for some transverse open arc $T \subseteq T'$ containing $x$, the first return with respect to $V$ is orientation-preserving{\rm)}. 
\end{definition}

Recall that a simple leaf arc is a simple arc contained in a leaf. 
For an arc $T$ and points $\alpha,\omega \in T$, denote by $[\alpha, \omega]_{T}$ the closed interval contained in $T$ whose boundary is $\{\alpha, \omega \}$.
We have the following statement for non-orientable first return maps. 

\begin{lemma}\label{lem:leaf_nonori_hol}
Let $\F$ be a divergence-free singular foliation on a surface $M$. 
For any transverse closed arc $T$ and any simple leaf arc $\gamma$ from a point $\alpha$ in $T$ and to the first return $\omega$ of $\alpha$ to $T$ from the same side of $T$, if there is a closed disk $B$ whose boundary is a loop $[\alpha,\omega]_{T} \cup \gamma$, then the following statements hold:
\\
{\rm(1)} The transverse closed arc $T$ intersects a semi-prong separatrix at $[\alpha, \omega]_{T}$.  
\\
{\rm(2)} If $[\alpha, \omega]_{T}$ intersects no non-proper leaves, then $T$ intersects a semi-prong separatrix at $[\alpha, \omega]_{T}$ which is transversely non-orientable on $B$.  
\end{lemma}

\begin{proof}
Lemma~\ref{lem:type_singular_points} implies that any singular points of $\F$ are either centers or prongs.
The union $C := [\alpha,\omega]_{T} \cup \gamma$ is a contractible simple closed curve as on the right in Figure~\ref{nonori_curve_disk}. 
\begin{figure}[t]
\begin{center}
\includegraphics[scale=0.75]{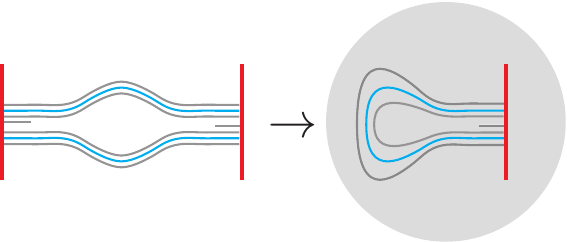}
\end{center} 
\caption{The tangent orientation branched covering of a transverse with a non-orientable first return.}
\label{nonori_curve_disk}
\end{figure}
Denote by $B$ a closed disk whose boundary is $C$. 
For any point $x$ in $[\alpha,\omega]_{T}$, we can take the maximal simple leaf arc $\gamma_x$ in $B$ from it with the initial direction corresponding to $\gamma$. 
Denote by $(\alpha,\omega)_{T,\mathrm{re}} \subsetneq [\alpha,\omega]_{T}$ the set of points $x$ of the closed interval $[\alpha,\omega]_{T}$ such that the boundaries of $\gamma_x$ are contained in $[\alpha,\omega]_{T}$.  
The double $S := B \sqcup_{\partial B} (- B)$ of $B$ is a closed surface with the induced foliation $\F_S$ of $\F\vert_B$ with finitely many singular points, which are either centers or prongs. 
The transverse closed arc $[\alpha,\omega]_{T}$ can be identified with the subset of $S$. 
Then any leaf intersecting $[\alpha,\omega]_{T}$ is a periodic leaf. 
By the existence of a fibered chart, the construction of $S$ implies that the union $\mathop{\mathrm{Per}}(\F_S)$ is open. 
From \cite[Lemma 3.1 and Lemma~4.2]{yokoyama2017decompositions}, any leaf $L \in \F_S$ intersecting $\overline{(\alpha,\omega)_{T,\mathrm{re}}} - (\alpha,\omega)_{T,\mathrm{re}} \subset [\alpha,\omega]_{T}$ is a prong separatrix. 
For any point $x \in \overline{(\alpha,\omega)_{T,\mathrm{re}}} - (\alpha,\omega)_{T,\mathrm{re}}$, the maximal simple leaf arc $\gamma_x$ is a semi-prong separatrix which connects to a prong in $B$.
The finiteness of prongs of $\F$ implies the finiteness of $\overline{(\alpha,\omega)_{T,\mathrm{re}}} - (\alpha,\omega)_{T,\mathrm{re}}$. 

Suppose that $[\alpha, \omega]_{T}$ intersects no non-proper leaves. 
Then $\overline{(\alpha,\omega)_{T,\mathrm{re}}} = [\alpha, \omega]_{T}$. 
Let $D(\F)$ be the prong connection diagram of $\F$. 
Therefore $(\alpha,\omega)_{T,\mathrm{re}}$ is the resulting space from $[\alpha,\omega]_{T}$ by removing the finite points. 

\begin{claim}\label{claim:02}
The closed interval $[\alpha,\omega]_{T}$ intersects leaves which are transversely non-orientable on the closed disk $B$.
\end{claim}
\begin{proof}[Proof of Claim~\ref{claim:02}]
Assume that $[\alpha,\omega]_{T}$ intersects no transversely non-orientable leaf.  
There are transverse orientations on leaves that intersect $T$ continuously. 
Then any connected components of $(\alpha,\omega)_{T,\mathrm{re}}$ have the same orientation. 
On the other hand, since $\partial \gamma = \{ \alpha, \omega \} = \partial [\alpha,\omega]_{T}$ is a part of the boundary $C = [\alpha,\omega]_{T} \cup \gamma$ of the closed disks $B$, the two connected components intersecting $\partial [\alpha,\omega]_{T}$ have the opposite transverse orientation, which contradicts that they have the same transverse orientation. 
\end{proof}

By the previous claim, there is a leaf $L$ intersecting $[\alpha,\omega]_{T}$ which is transversely non-orientable.
Since the intersection $L \cap B$ is contained in a closed disk $B$, the intersection $L \cap B$ is a semi-prong separatrix. 

\begin{claim}\label{claim:03}
The closed interval $[\alpha,\omega]_{T}$ intersects leaves that are non-orientable on the closed disk $B$.
\end{claim}
\begin{proof}[Proof of Claim~\ref{claim:03}]
Assume that $[\alpha,\omega]_{T}$ intersects no non-orientable leaf.  
Fix the orientation of leaves of the connected component of $(\alpha,\omega)_{T,\mathrm{re}}$ containing $\omega$ continuously. 
The orientation implies one of the connected components of $(\alpha,\omega)_{T,\mathrm{re}}$ containing $\alpha$ continuously because $\gamma$ is a part of the boundary of the closed disk $B$. 
Then the orientations of the intervals induce those of the boundaries of the intervals. 
Moreover, the orientations of the boundaries induce those of intervals whose boundaries intersect them. 
Since $(\alpha,\omega)_{T,\mathrm{re}}$ is the resulting space from $[\alpha,\omega]_{T}$ by removing the finite points, by finite iterations, we have the same orientation on $[\alpha,\omega]_{T}$ continuously. 
On the other hand, the existence of the leaf arc $\gamma$ which is a part of the boundary of the closed disk $B$ implies that the finite disjoint union $[\alpha,\omega]_{T}$ of intervals has the opposite orientation on $[\alpha,\omega]_{T}$, which contradicts the existence of same orientation on it. 
\end{proof}
By the previous claim, there is a leaf $L'$ intersecting $[\alpha,\omega]_{T}$ which is non-orientable.
Since the intersection $L' \cap B$ is contained in a closed disk $B$, the intersection $L' \cap B$ is a semi-prong separatrix. 
\end{proof}

\begin{figure}[t]
\begin{center}
\includegraphics[scale=0.375]{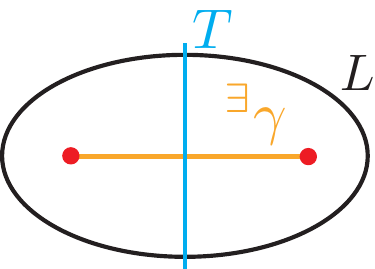}
\end{center} 
\caption{A transverse closed arc $T$ which intersects a contractible periodic leaf $L$ exactly twice must intersect a prong separatrix $\gamma$.}
\label{quasi-circle_01}
\end{figure}

Applying the previous lemma to the half disk as in Figure~\ref{quasi-circle_01}, we have the following property as in Figure~\ref{quasi-circle_01}. 


\begin{corollary}\label{cor:leaf_separatrix}
Let $\F$ be a divergence-free singular foliation on a surface $M$. 
For any transverse closed arc $T$ and any contractible periodic leaf $L$ bounding an open disk $D$ such that $L$ intersects exactly twice to the interior of $T$ and the set difference $D \setminus T$ consists of two disjoint open disks, the interior of the arc $T$ intersects a prong separatrix. 
\end{corollary}

\subsection{Divergence-free singular foliations on compact surfaces}

We have the following statements for compact surface cases. 

\begin{lemma}\label{lem:leaf_hol}
Any divergence-free singular foliation without non-proper leaves on a compact surface consists of periodic leaves except finitely many leaves, which are prongs and prong separatrices. 
\end{lemma}

\begin{proof}
Let $\pi \colon M^* \to M$ be the tangent orientation branched covering and $\F^*$ the lift of $\F$ on $M^*$ which is the set of orbits of a smooth divergence-free vector field $X$.
Then $X$ has no non-proper orbits.
Since any divergence-free vector field is non-wandering, by Lemma~\ref{lem:type_singular_points}, any singular points of $X$ are either centers or multi-saddles. 
By \cite[Lemma~2.4 and Proposition~2.6]{yokoyama2016topological}, since there are at most finitely many singular points, the periodic point set of $X$ is open and dense, and any non-closed orbits are multi-saddles separatrices. 
This means that the complement of the periodic point set consists of finitely many centers, multi-saddles, and multi-saddle separatrices. 
Therefore, the assertion holds. 
\end{proof}

\section{Levelable singular foliations}

To introduce an analogous concept of a Hamiltonian vector field, we recall the following concepts. 

\subsection{Levelable foliations}

We now introduce an analogous concept of a Hamiltonian vector field with finitely many singular points, as follows. 

\begin{definition}\label{def:lebelable}
A singular foliation $\F$
on a surface $M$ is {\bf levelable} if $\F$ is tangent to the boundary and there is a continuous function $H \colon M \to \R$ such that the following conditions hold:
\\
{\rm(1)} Each connected component of the inverse image $H|_{M - \mathop{\mathrm{Sing}}(\F)}^{-1}(c)$ for any $c \in H(M - \mathop{\mathrm{Sing}}(\F))$ is a leaf of $\F$. 
\\
{\rm(2)} For any non-singular transversely orientable leaf $L$ and any small transverse arc $T$ from a point in $L$ whose saturation $\F(T)$ is a \nbd of $L$, the restriciton $H \vert_T$ is strictly monotonic. 
\end{definition}

Then, the function $H$ is called the {\bf height function} of $\F$ and is denoted by $H_\F$. 
By definition, the set of orbits of a Hamiltonian is a levelable foliation. 
On the other hand, a levelable foliation need not be the set of orbits of a Hamiltonian flow even if the foliation is the set of orbits of a non-singular flow (see Example~\ref{ex:nonHam}). 
In addition, notice that the singular foliation on a periodic torus is not levelable. 
Denote by $\bm{\mathcal{L}(M)} \subset \F(M)$ the set of levelable foliations on a surface $M$,  and by $\bm{\mathcal{L}^\star}(M) \subset \mathcal{L}(M)$ the subspace of levelable foliations with finitely many pronged singular points on $M$, and by $\bm{\mathcal{L}^*}(M) \subset \mathcal{L}^\star(M)$ the subspace of levelable foliations without fake prongs on $M$.
We have the following observation. 

\begin{lemma}\label{lem:regularity_proper}
Let $\F$ be a levelable foliation on a surface $M$. 
For any leaf $L \in \F$, the set difference $\overline{L} - L$ consists of singular points. 
\end{lemma}

\begin{proof}
Denote by $H_\F$ the height function of $\F$. 
Fix any leaf $L \in \F$. 
By the definition of a levelable foliation, there is a value $c \in H(M - \mathop{\mathrm{Sing}}(\F))$ such that the leaf $L$ is a connected component of the inverse image $H|_{M - \mathop{\mathrm{Sing}}(\F)}^{-1}(c)$. 
Therefore, the leaf $L$ is closed in $M - \mathop{\mathrm{Sing}}(\F)$ and so $\overline{L} - L \subseteq \mathop{\mathrm{Sing}}(\F)$. 
\end{proof}

%

%
%
%

\subsubsection{Topological graphs and Reeb graphs}

Recall several topological concepts as follows. 

\begin{definition}
A {\bf topological graph} is a cell complex whose dimension is at most one and which is a geometric realization of a finite abstract multi-graph.
\end{definition}

A topological graph which is a geometric realization of a finite abstract directed multi-graph equipped with the directed structure is called a {\bf (topological) directed graph}.

\begin{definition}
A {\rm(}topological{\rm)} {\bf directed graph $\bm{G}$ with possibly semi-open edges} is homeomorphic to the resulting space from a topological directed graph $H$ by removing vertices. 
Then the topological directed graph $H$ is called the {\bf end-completion} of $G$. 
\end{definition}

\begin{definition}
A directed graph (with possibly semi-open edges) is {\bf acyclic} if it contains no directed cycles. 
\end{definition}

\begin{definition}
A directed graph (with possibly semi-open edges) is {\bf finite} if it consists of finitely many vertices and edges. 
\end{definition}

A topological directed graph equipped with lengths of edges is called a {\bf directed graph with edge lengths}, and the lengths of edges are also called the {\bf edge lengths}. 
Denote by $\bm{\mathcal{G}}$ the set of directed graphs equipped with edge lengths. 
Recall the Reeb graphs as follows. 

\begin{definition}
For a function $f \colon  X \to \R$ on a topological space $X$, the {\bf Reeb graph} $G_f$ of a function $f \colon  X \to \R$ on a topological space $X$ is a quotient space $X/\mathop{\sim}_{\mathrm{Reeb}}$ defined by $x \sim_{\mathrm{Reeb}} y$ if there are a number $c \in \R$ and a connected component of $f^{-1}(c)$ which contains $x$ and $y$.
\end{definition}

Note that Reeb graphs are not topological directed graphs with possibly semi-open edges in general. 
We have the following statement. 

\begin{lemma}\label{lem:periodic_noncompact}
The following statements hold for any levelable foliation with pronged singular points on a connected surface:
\\
{\rm(1)} Every leaves are proper. 
\\
{\rm(2) } The prong connection diagram $D(\F)$ is a closed subset which consists of fintely many leaves. 
\\
{\rm(3)} The complement $M - D(\F)$ consists of finitely many connected components. 
\\
{\rm(4)} Any connected component of  $M - D(\F)$ is either a periodic annulus, a periodic M{\"o}bius band, a periodic Klein bottle, or a trivial fibered chart which consists of non-closed leaves which are closed subsets. 
\\
{\rm(5)} The Reeb graph of the height function of a levelable foliation on a surface is an acyclic finite directed graph with possibly semi-open edges.
\\
{\rm(6)} The union $M - D(\F)$ of non-singular leaves which are closed subsets is open and dense.  
\end{lemma}

\begin{proof}
Let $\F$ be a levelable foliation with pronged singular points on a surface $M$ and $H$ the height function of $\F$. 
Lemma~\ref{lem:regularity_proper} implies that all leaves are proper. 
Therefore, assertion (1) holds. 

From Lemma~\ref{lem:regularity_proper}, since all singular points are prongs, all non-closed leaves are semi-prong separatrices. 
By definition of the prong, any singular foliation with pronged singular points has at most finitely many singular points and semi-prong separatrices.  
Then the prong connection diagram $D(\F)$ consists of finitely many leaves and is a closed subset. 
Therefore, assertion (2) holds. 

Since all non-closed leaves are prong separatrices, any leaf is either a semi-prong separatrix, a closed leaf, or a non-closed leaf which is a closed subset.
Then the complement $M - D(\F)$ consists of finitely many connected components such that any connected component is contained in either the union of periodic leaves or the union of non-closed leaves which are closed subsets. 
By the existence of the height function of $\F$ and of fibered chart, any connected component of $M - D(\F)$ consisting of non-closed leaves which are closed subsets is a trivial flow flox, and any connected component of $M - D(\F)$ consisting of periodic leaves. 

\begin{claim}\label{claim:10}
Any connected component of $M - D(\F)$ consisting of periodic leaves is either a periodic annulus, a periodic M{\"o}bius band, or a periodic Klein bottle. 
\end{claim}

\begin{proof}[Proof of Claim~\ref{claim:10}]
Since the first return of any point in a periodic leaf $L$ to a transverse open interval along $L$ is either the identical homeomorphism or non-orientable involutive homeomorphism, each connected component of $\mathop{\mathrm{Per}}(\F)$ is either an annulus, a M{\"o}bius band, a torus, or a Klein bottle.
Since any periodic torus is not levelable, the assertion holds. 
\end{proof}

Claim~\ref{claim:10} implies assertion (4). 
Notice that, for any connected component $C$ of $M - D(\F)$ consisting of periodic leaves which is non-orientable, the connected component $C$ contains non-orientable periodic leaves and the complement $C - C_N$ of the union $C_N$ of non-orientable periodic leaves in $C$ is a periodic annulus. 
By assertions (2)--(4), the Reeb graph of $H$ is a finite directed graph whose vertices are prong connections and boundary components of $M$ and transversely non-orientable periodic leaves, and whose edges are periodic annuli and trivial fibered chartes consisting of non-closed leaves which are closed subsets. 
The existence of a height function implies the non-existence of directed cycle, and so assertion (5) holds. 
In addition, the complement $M - D(\F)$ is an open dense subset which is the union of non-singular leaves which are closed subsets, which implies assertion (6). 
\end{proof}

The previous lemma implies the following observations in the case on a spherical case. 

\begin{corollary}\label{cor:tree_rep-}
The Reeb graph of the height function of a levelable foliation $\F$ with pronged singular points on a surface $S$ contained in a sphere is a finite directed tree with possibly semi-open edges. 
\end{corollary}

\begin{theorem}\label{cor:tree_rep}
The Reeb graph of the height function of a levelable foliation $\F$ with pronged singular points on a sphere $\mathbb{S}^2$ is a finite directed tree. 
Moreover, the complement $\mathbb{S}^2 - D(\F)$ is a finite disjoint union of open periodic annuli and is open and dense. 
In particular, for any levelable foliation $\F_1, \F_2$ on $\mathbb{S}^2$, the foliation $\F_1$ and $\F_2$ are topologically equivalent if and only if their prong connection diagrams $D(\F_1)$ and $D(\F_2)$ are isomorphic as a surface graph. 
\end{theorem}


We have the following statement. 

\begin{lemma}\label{lem:leaf_hol_02}
Let $\F$ be a levelable singular foliation with pronged singular points on a surface $M$. 
For any transverse closed arc $T$ and any simple leaf curve $\gamma$ from a point $\alpha \in T$ and to the first return $\omega \in T$, if $T$ contains no transversely non-orientable leaves, then $\alpha = \omega$ and so $\gamma$ is a periodic leaf.  
\end{lemma}

\begin{proof}
Let $H$ be the height function of $\F$.
Fix a transverse closed arc $T$ which intersects only transversely orientable leaves, and fix a simple leaf curve $\gamma$ from a point $\alpha \in T$ and to the first return $\omega \in T$. 
The transverse orientability on $T$ implies the strict monotonicity of $H \vert_T$. 
Since $\gamma$ is contained in $H^{-1}(c)$ for some value $c \in \R$, we have $H(\alpha) = H(\omega)$, which contracts $H(\alpha) \neq H(\omega)$. 
\end{proof}

\subsection{Properties for levelable foliations on compact surfaces}
We have the following statement. 

\begin{lemma}\label{lem:periodic}
The following statements hold for any levelable foliations with pronged singular points on compact surfaces:
\\
{\rm(1)} Any leaves except finite exceptions are periodic. 
\\
{\rm(2) } The Reeb graph of a levelable foliation on a compact surface is an acyclic finite directed graph.
\\
{\rm(3)} The union of non-periodic leaves consists of finitely many singular leaves and prong separatrices. 
\\
{\rm(4)} The union of periodic leaves is open and dense.  
\end{lemma}

\begin{proof}
Let $\F$ be a levelable foliation with pronged singular points on a compact surface $M$. 
Since any leaves that are closed subsets are closed leaves, there are no non-closed leaves that are closed subsets. 
By Lemma~\ref{lem:periodic_noncompact}, the assertions (1), (3), and (4) hold. 
The compactness of $M$ implies the non-existence of semi-open edges, and so the assertion (2) holds because of Lemma~\ref{lem:periodic_noncompact}(5). 
\end{proof}

\begin{figure}[t]
\begin{center}
\includegraphics[scale=0.375]{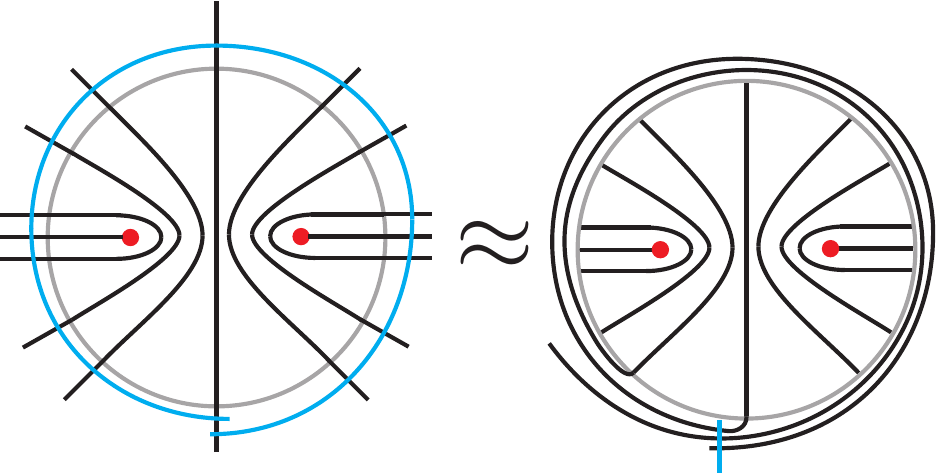}
\end{center} 
\caption{A contractible loop consisting of a transverse closed arc and a simple leaf arc along which the first return between distinct points in the transversal is orientable.}
\label{disk_2prongs}
\end{figure}

Notice that the transversely non-orientable condition in the previous lemma is necessary. 
In fact, there is a contractible loop consisting of a transverse closed arc and a simple leaf arc along which the first return between distinct points in the transversal is orientable as in Figure~\ref{disk_2prongs} such that intersecting no proper leaves, because the disk can be embedded in a sphere with a singular foliation with four thorns and minimal (non-proper) leaves as in Example~\ref{ex:rot_disk}. 
%
%
%
We have the following implication. 
 
\begin{lemma}
Any levelable foliation with pronged singular points on a compact surface is divergence-free. 
\end{lemma}


\begin{proof}
Let $F$ be a levelable foliation on a compact surface $M$ and $\F^*$ the tangent orientation lift of $\F$ on the tangent orientation branched covering $M^*$ of $M$. 
From Lemma~\ref{lem:periodic}, any leaves except finite exceptions of $\F$ are periodic. 
Since a regular foliation on a surface is orientable if and only if there is a flow on the surface whose set of orbits is the foliation (cf. \cite[Remark~2.3.2]{HH1986A}), there is a flow $f \colon \R \times M^* \to M^*$ whose set of orbits is the foliation $\F^*$. 
Since any leaves except finite exceptions of $\F$ are periodic, the tangent orientation lift $\F^*$ of $\F$ consists of periodic leaves except finite exceptions. 
This means that the periodic point set $\mathop{\mathrm{Per}}(f)$ of $f$ is dense, and so the flow $f$ is non-wandering. 
By \cite[Theorem~A]{yokoyama2021ham}, the non-wandering flow $f$ is topologically equivalent to the flow generated by a smooth divergence-free vector field.
From Lemma~\ref{lem:ch:levelable}, the assertion holds. 
\end{proof}

\section{Topologies of the set of levelable foliations on surfaces}

Fix a Riemannian metric $g$ on a surface $M$. 
Let $\F$ be a levelable foliation on $M$ whose height function $H$ is of piecewise $C^1$. 

\subsection{Balls on the set of levelable foliations on surfaces}

To define distances on the set of levelable foliations on surfaces, we define weight for non-periodic leaves. 

%
%

\subsubsection{Weight for non-periodic leaves}

We define the weight for non-periodic leaves.

\begin{definition}
For a distance weight $W \colon [0,\infty] \to [0,1]$, 
the {\bf weight} $\bm{\sigma_{W}} \colon M \times \mathcal{L}(M) \to [0,1]$ with respect to $W$ and $\mathop{\mathrm{Per}}(\F)$ is defined by 
\[
\sigma_{W}(x, \F) := W(d_g(x, M - \mathop{\mathrm{Per}}(\F)))
\]
\end{definition}

Put $\sigma_{W,\F} := \sigma_{W}(\cdot, \F) \colon M \to [0,1]$. 
%
For the weight $\sigma_{W,\F} \colon M \to [0,1]$ with respect to a distance weight $W \colon [0,\infty] \to [0,1]$ and $\F$, put 
\[
X_{W_{\mathrm{per}},\F_-}(x) :=  \sigma_{W,\F}(x)X_{\F_-}(x)
\hspace{20pt}
X_{W_{\mathrm{per}},\F_+}(x) :=  \sigma_{W,\F}(x)X_{\F_+}(x)
\]
for any $x \in M$. 
We define the continuous associated directed field of $\F$ with respect to the weight $\sigma_{W,\F} \colon M \to [0,1]$ as follows. 

\begin{definition}
Define the {\bf continous associated directed field} $X_{W_{\mathrm{per}},\F} \colon M \to TM_\pm$ of $\F$ with respect to the weight $\sigma_{W,\F} \colon M \to [0,1]$ by 
\[
\begin{split}
X_{W_{\mathrm{per}},\F}(x) := \sigma_{W,\F}(x) X_{\F}(x) &= [0_x, \sigma_{W,\F}(x)X_{\F_-}(x)] \cup [0_x, \sigma_{W,\F}(x)X_{\F_+}(x)]
\\
&= [0_x, X_{W_{\mathrm{per}},\F_-}(x)] \cup [0_x, X_{w,\F_+}(x)]
\end{split}
\]
for any $x \in M$.
\end{definition}

Notice the continous associated directed field $X_{W_{\mathrm{per}},\F} \colon M \to TM_\pm$ of the levelable foliation $\F$ is $\{ \{ X_{W_{\mathrm{per}},\F_-}(x), X_{W_{\mathrm{per}},\F_+}(x) \} \mid x \in M \}$, which consists of two/one-point sets, such that $X_{W_{\mathrm{per}},\F}\vert_{M - \mathop{\mathrm{Per}}(\F)} \equiv \{ 0 \}$.

\subsubsection{Distance between singular foliations and their topology}

Fix a distance weight $w,W \colon [0,\infty] \to [0,1]$. 

\begin{definition}
Define the wighted metric $\bm{D_{W_{\mathrm{per}}}}$, called the {\bf locally wighted levelable Hausdorff superior distance} for periodic leaves on the set $\mathcal{L}(M)$ of singular foliations on a  manifold $M$ for an ordered field $[0, \infty]$: 
\[
\begin{split}
\bm{D_{W_{\mathrm{per}}}(\F, \mathcal{G})} = D^{\infty}_{W_{\mathrm{per}}}(\F, \mathcal{G}) &:= \sup_{x \in M} d_H(X_{W_{\mathrm{per}},\F}(x), X_{W_{\mathrm{per}},\G}(x))
\end{split}
\]
where $d_H$ is the Hausdorff distance on $M$ with respect to the distance $d$. 
\end{definition}

We define the following distance. 

\begin{definition}
Define the wighted metric $\bm{D_{w,W}}$, called the {\bf locally wighted levelable Hausdorff superior distance} on the set $\mathcal{L}(M)$ of singular foliations on a  manifold $M$: 
\[
\begin{split}
\bm{D_{w,W}(\F, \mathcal{G})} = D^{\infty}_{w,W}(\F, \mathcal{G}) &:= D(\F, \mathcal{G}) + D_{W_{\mathrm{per}}}(\F, \mathcal{G})
\end{split}
\]
\end{definition}

Similarly, we define the following distance. 

\begin{definition}
For any $p \in \R_{\geq 1}$, any $r \in \R_{> 0}$, and any $\delta \in \R_{> 0}$ and a regular foliation $\F$ on a  manifold $M$, define the wighted metrics $\bm{D^{p}_{W_{\mathrm{per}}}}$ and $\bm{D^{p}_{w,W}}$, which is called the {\bf locally wighted levelable Hausdorff $L^p$-distance}, on the set $\mathcal{L}(M)$ of singular foliations on a  manifold $M$ for an ordered field $[0, \infty]$: 
\[
\begin{split}
\bm{D^{p}_{W_{\mathrm{per}}}(\F, \mathcal{G})} := & \left(\dfrac{1}{\mathrm{vol}_g(M)}  \int_M (d^+_{w,W}(\F, \mathcal{G})(x))^p d \mathrm{vol}_g \right)^{1/p}
\end{split}
\]
\[
\begin{split}
\bm{D^{p}_{w,W}(\F, \mathcal{G})} := & D^p_{w}(\F, \mathcal{G}) + D^{p}_{W_{\mathrm{per}}}(\F, \mathcal{G})
\end{split}
\]
\end{definition}

\subsubsection{Open $\delta$-balls for levelable singular foliations on manifolds and their topology}
Fix distance weights $w, W \colon [0,\infty] \to [0,1]$. 

\begin{definition}
For any $p \in \R_{\geq 1} \sqcup \{ \infty \}$, any $\delta \in \R_{> 0}$, and a levelable singular foliation $\F$ on the manifold $M$, define the {\bf open $\delta$-ball} $\bm{B_{D^{p}_{w,W}}(\F,\delta)}$ centered at $\F$ with respect to $(w, W)$ as follows:  
\[
B_{D^{p}_{w,W}}(\F,\delta) := \{ \mathcal{G} \in \mathcal{L}(M) \mid D^{p}_{w,W}(\F, \mathcal{G}) < \delta \}
\]
\end{definition}

By definitions, we have $B_{D^{p}_{w,W}}(\F,\delta) \subseteq B_{D^{p}_{W_{\mathrm{per}}}}(\F,\delta)$. 
Moreover, we set the following notation. 

\begin{definition}
Define the set $\mathcal{B}^{p}_{w,W,\mathcal{L}(M)}$ of open balls with respect to the distance weights $(w, W)$ as follows: 
\[
\mathcal{B}^{p}_{w,W,\mathcal{L}(M)} := \{ B_{D^{p}_{w,W}}(\F,\delta) \mid \delta \in \R_{> 0}, \F \in \mathcal{L}(M) \}
\]
\end{definition}

\subsubsection{Topologies for levelable singular foliations on surfaces}

\begin{definition}\label{dec:levelable_top}
For any $p \in \R_{\geq 1} \sqcup \{ \infty \}$ and any distance weight $w, W \colon [0,\infty] \to [0,1]$, the smallest topology on $\mathcal{L}(M)$ which contains the set $\mathcal{B}^{p}_{w,W,\mathcal{L}(M)} $ of open balls  with respect to $w$ is called the {\bf topology induced by the local Hausdorff distance} with respect to $w$ and $W$, and is denoted by $\mathcal{O}^p_{w,W,\mathcal{L}(M)}$. 
\end{definition}

Then the set $\mathcal{B}^{p}_{w,W,\mathcal{L}(M)} $ is a subbase of the topology $\mathcal{O}^p_{w,W,\mathcal{L}(M)}$. 

\section{Generic coheight-zero levelable foliationswith finitely many pronged singular points on surfaces}

Let $\F$ be a levelable foliation with finitely many pronged singular points on a surface $M$.  

Recall that a {\bf bouquet} of $k$ circles for any $k \in \Z_{\geq 0}$ is a topological graph that consists of one vertex and $k$ loops. 
We have the following observation.

\begin{lemma}\label{lem:betti}
The first betti number $b_1(G)$ of a connected finite graph $G = (V,E)$ is $1 - \vert V \vert + \vert E \vert$.   
\end{lemma}

\subsection{Coheight of levelable foliations}

To define the coheight of levelable foliations, we define coheights with respect to heteroclinic separatrices and with respect to the degeneracy of singular points.

\subsubsection{Coheights of prong connections of levelable foliations}

We have the following observation. 

\begin{lemma}
Let 
$D(\F)$ be the prong connection diagram.  
Then $D(\F)$ is a topological graph whose degrees of vertices are one or three such that every connected component $C$ of $D(\F)$ is homotopic to a bouquet of $b_1(C)$ circles such that $b_1(C) = 1 - n_{V}(C) + n_{E}(C)$, where $n_{V}(C)$ is the number of prongs in $C$ and $n_{E}(C)$ is the number of prong separatrices in $C$. 
\end{lemma}

We define coheights with respect to heteroclinic separatrices as follows. 

\begin{definition}
Define the {\bf coheight} $\bm{\mathrm{coheight}_{\mathrm{h},\F}(C)}$ of a prong connection $C$ with respect to heteroclinic separatrices as follows: 
\[
\mathop{\mathrm{coheight}_{\mathrm{h},\F}}(C) := n_{V_-}(C) - b_1(C)
\]
where $n_{V_-}(C) := \vert \mathop{\mathrm{Sing}_{-}}(\F) \cap C \vert$ is the number of prongs in $C$ whose indices are negative. 
\end{definition}

\begin{definition}
The {\bf coheight} $\bm{\mathrm{coheight}_{\mathrm{h}}(\F)}$ of $\F$ with respect to heteroclinic separatrices is defined as follows: 
\[
\begin{split}
\mathrm{coheight}_{\mathrm{h}}(\F) := \, &  \sum_{C:\text{prong connection}} \hspace{-5pt} \mathop{\mathrm{coheight}_{\mathrm{h},\F}}(C) 
\\
= \, &  n_{V_-}(\F)  - b_1(D(\F))
\\
= \, & n_{V_-}(\F) +  n_{V}(\F) - n_{E}(\F) -n_{CC}(\F) 
\end{split}
\]
where $n_{V_-}(\F) := \vert \mathop{\mathrm{Sing}_{-}}(\F) \vert$ is the number of prongs of $\F$ whose indices are negative, $n_{V}(\F) := \vert \mathop{\mathrm{Sing}}(\F) \vert$ is the number of prongs of $\F$, $n_{E}(\F)$ is the number of prong separatrices of $\F$, and $n_{CC}(\F)$ is the number of connected components of the prong connection diagram $D(\F)$ 
\end{definition}

\subsubsection{Coheights of singular points of levelable foliations}

We define coheights with respect to the degeneracy of singular points as follows. 

\begin{definition}
The {\bf coheight} $\bm{\mathrm{coheight}_{\mathrm{p},\F}(y)}$ of a center $y$ with respect to degeneracy of singular points is one. 
\end{definition}

\begin{definition}
The {\bf coheight} $\bm{\mathrm{coheight}_{\mathrm{p},\F}(y)}$ of a fake prong $y$ outside of the boundary with respect to the degeneracy of singular points is zero. 
\end{definition}

\begin{definition}
The {\bf coheight} $\bm{\mathrm{coheight}_{\mathrm{p},\F}(y)}$ of a $1$-prong $y$ with respect to degeneracy of singular points is zero. 
\end{definition}

\begin{definition}
The {\bf coheight} of $k$-prong $y$ ($k \geq 3$) outside of the boundary with respect to degeneracy of singular points is $\bm{\mathrm{coheight}_{\mathrm{p},\F}(y)} := 2(k-3) = -2(1 + 2 \mathop{\mathrm{ind}_{\F}}(y)) \geq 0$, because $\mathop{\mathrm{ind}_{\F}}(y) = -(k-2)/2$.
\end{definition}

\begin{definition}
The {\bf coheight} $\bm{\mathrm{coheight}_{\mathrm{p},\F}(y)}$ of a fake prong $y$ on the boundary with respect to the degeneracy of singular points is two. 
\end{definition}

\begin{definition}
The {\bf coheight} of $k$-prong $y$ ($k \geq 3$) on the boundary with respect to degeneracy of singular points is $\bm{\mathrm{coheight}_{\mathrm{p},\F}(y)} := 2(k-3) = -2(1 + 2 \mathop{\mathrm{ind}_{\F}}(y)) \geq 0$. 
\end{definition}

\begin{definition}
The {\bf coheight} $\bm{\mathrm{coheight}_{\mathrm{p,\F}}(C)}$ of a prong connection $C$ of $\F$ with respect to degeneracy of singular points is defined as follows: 
\[
\begin{split}
\mathrm{coheight}_{\mathrm{p,\F}}(C) := & \sum_{y \in \mathop{\mathrm{Sing}}(\F) \cap C} \mathrm{coheight}_{\mathrm{p},\F}(y)
\end{split}
\]
\end{definition}

\begin{definition}
The {\bf coheight} $\bm{\mathrm{coheight}_{\mathrm{p}}(\F)}$ of $\F$ with respect to degeneracy of singular points is defined as follows: 
\[
\begin{split}
\mathrm{coheight}_{\mathrm{p}}(\F) := & \sum_{y \in \mathop{\mathrm{Sing}}(\F)} \mathrm{coheight}_{\mathrm{p},\F}(y)
\end{split}
\]
\end{definition}

Notice that $\mathrm{coheight}_{\mathrm{p,\F}}$ is non-negative for any of a prong, a prong connection, or a line field. 

\subsubsection{Coheights of prong connections of levelable foliations}

We define coheights of prong connections of levelable foliations as follows. 

\begin{definition}
The {\bf coheight} $\bm{\mathrm{coheight}_{\F}(C)}$ of a prong connection $C$ of $\F$ 
is defined as follows: 
\[
\begin{split}
& \mathrm{coheight}_{\F}(C) \\
:= \, &  n_{\partial}(C) + n_{\mathrm{fake}}(C) + \mathop{\mathrm{coheight}_{\mathrm{h},\F}}(C) +  \mathrm{coheight}_{\mathrm{p,\F}}(C)
\\
= \, & n_{\partial}(C) + n_{\mathrm{fake}}(C) + n_{V_-}(C)  - b_1(C) + \sum_{y \in \mathop{\mathrm{Sing}}(\F) \cap C} \mathrm{coheight}_{\mathrm{p},\F}(y)
\\
= \, & n_{\partial}(C) + n_{\mathrm{fake}}(C) + n_{V_-}(C)  +  n_{V}(C) - n_{E}(C) -1 + \sum_{y \in \mathop{\mathrm{Sing}}(\F) \cap C} \mathrm{coheight}_{\mathrm{p},\F}(y)
\\
= \, & n_{\partial}(C) + n_{\mathrm{fake}}(C) + n_{V_-}(C)  - n_{E}(C) -1 + \sum_{y \in \mathop{\mathrm{Sing}}(\F) \cap C} (1 + \mathrm{coheight}_{\mathrm{p},\F}(y))
\end{split}
\]
where $n_{\partial}(C)$ is the number of boundary components of $M$ intersecting $C$, and $n_{\mathrm{fake}}(C)$ is the number of fake prongs in $C$. 
\end{definition}

\subsubsection{Coheights of levelable foliations}

We define coheights of levelable foliations as follows. 

\begin{definition}
The {\bf coheight} $\bm{\mathrm{coheight}(\F)}$ of $\F$ 
is defined as follows: 
\[
\begin{split}
& \mathrm{coheight}(\F)\\
 :=  & \sum_{C:\text{prong connection}} \mathrm{coheight}_{\F}(C) 
\\
 =  \, & n_{\partial}(\F) + n_{\mathrm{fake}}(\F) +  \mathop{\mathrm{coheight}_{\mathrm{h}}}(\F) + \mathrm{coheight}_{\mathrm{p}}(\F)
\\
= \, & n_{\partial}(\F) + n_{\mathrm{fake}}(\F) + n_{V_-}(\F) - b_1(D(\F)) + \sum_{y \in \mathop{\mathrm{Sing}}(\F)} \mathrm{coheight}_{\mathrm{p},\F}(y)
\\
= \, & n_{\partial}(\F) + n_{\mathrm{fake}}(\F) + n_{V_-}(\F) +  n_{V}(\F) - n_{E}(\F) -n_{CC}(\F) + \sum_{y \in \mathop{\mathrm{Sing}}(\F)} \hspace{-7pt} \mathrm{coheight}_{\mathrm{p},\F}(y)
\\
= \, & n_{\partial}(\F) + n_{\mathrm{fake}}(\F) + n_{V_-}(\F) - n_{E}(\F) -n_{CC}(\F) + \sum_{y \in \mathop{\mathrm{Sing}}(\F)} (1 + \mathrm{coheight}_{\mathrm{p},\F}(y))
\end{split}
\]
where $n_{\partial}(\F)$ is the number of boundary components of $M$ intersecting the prong connection diagram $D(\F)$ and $n_{\mathrm{fake}}(\F)$ is the number of fake prongs of $\F$. 
\end{definition}

List the numbers related to the coheight of a levelable foliation in Table~\ref{no_codim}. 

\begin{table}[t]\label{no_codim}
\centering
\begin{tabular}{|l|l|c|}
\hline
\textbf{Notation} & \textbf{Definition} & \textbf{codim} \\
\hline
\begin{tabular}{l} $n_{V_-}(\mathcal{F})$ \\ $= \left| \mathop{\mathrm{Sing}_{-}}(\mathcal{F}\right|$ \end{tabular} \hspace{-10pt} & the number of prongs of $\mathcal{F}$ with negative indices  & + \\
\hline
\begin{tabular}{l}$n_{V}(\mathcal{F})$ \\ $= \left| \mathop{\mathrm{Sing}}(\mathcal{F})\right|  $ \end{tabular} \hspace{-10pt} & the number of prongs of $\mathcal{F}$ & + \\
\hline
$n_{E}(\mathcal{F})$ &  the number of prong separatrices of $\mathcal{F}$ & - \\
\hline
$n_{CC}(\mathcal{F})$ & \begin{tabular}{c} the number of connected components of \\ the prong connection diagram $D(\mathcal{F})$ \end{tabular} & - \\
\hline
$n_{\partial}(\mathcal{F})$ &  
\begin{tabular}{c}
the number of boundary components of $M$ \\ intersecting the prong connection diagram $D(\mathcal{F})$ 
\end{tabular} & +
\\
\hline
$n_{\mathrm{fake}}(\mathcal{F})$ &  the number of fake prongs of $\mathcal{F}$ & + \\
\hline
\end{tabular}
\caption{Numbers related to the coheight of a levelable foliation $\mathcal{F}$ with finitely many pronged singular points}
\end{table}

%

\subsubsection{Properties of $\mathop{\mathrm{coheight}_{\mathrm{h},\F}}$ and $\mathop{\mathrm{coheight}_{\mathrm{h},\F}} + \mathop{\mathrm{coheight}_{\mathrm{p},\F}}$}

We have the following positivities. 

\begin{lemma}
For a prong connection $C$ each of whose prongs is either a $1$-prong and a $3$-prong, we have $\mathop{\mathrm{coheight}_{\mathrm{h},\F}}(C) \geq 0$. 
\end{lemma}

\begin{proof}
Fix any prong connection $C$ each of whose prongs is either a $1$-prong or a $3$-prong. 
Then $C$ contains at least two prongs. 
Since any prong in $C$ is either a $1$-prong and a $3$-prong, we have $n_{E}(C) = (n_{V_+}(C) + 3n_{V_-}(C))/2 = n_{V}(C)/2 + n_{V_2}(C)$. 
Lemma~\ref{lem:betti} implies that 
\[
\begin{split}
b_1(C) &= 1 - n_{V}(C) + n_{E}(C) 
\\
&= 1 - n_{V}(C) + n_{V}(C)/2 + n_{V_-}(C) = 1 - n_{V}(C)/2 + n_{V_-}(C)
\end{split}
\]
where $n_{V}(C)$ is the number of prongs in $C$ and $n_{E}(C)$ is the number of prong separatrices in $C$. 
By definition of the coheight $\mathop{\mathrm{coheight}_{\mathrm{h},\F}}$, we obtain the following equality: 
\[
\begin{split}
\mathop{\mathrm{coheight}_{\mathrm{h},\F}}(C) &= n_{V_-}(C) - b_1(C) 
\\
&= n_{V_-}(C) - \left(1 - \dfrac{n_{V}(C)}{2} + n_{V_-}(C) \right) = \dfrac{n_{V}(C)}{2} - 1
\end{split}
\]
Since $C$ contains at least two prongs, we have $\mathop{\mathrm{coheight}_{\mathrm{h},\F}}(C) = n_{V}(C)/2 - 1 \geq 0$. 
\end{proof}

\begin{lemma}\label{lem:1_3_prongs}
For a prong connection $C$, we have 
\[
\mathop{\mathrm{coheight}_{\mathrm{h},\F}}(C) + \mathop{\mathrm{coheight}_{\mathrm{p},\F}}(C) \geq n_0(C) +  \sum_{k \geq 3} (k-3) n_k(C) \geq 0
\]
where $n_k(C)$ is the number of $k$-prongs in $C$ for any $k \in \Z_{\geq 0}$. 
Moreover, if $C$ has no centers, then there is a prong connection $C'$ with $n_{\partial}(C) = n_{\partial}(C')$ and $n_{\mathrm{fake}}(C) = n_{\mathrm{fake}}(C')$ and 
\[
n_V(C') - n_V(C) = n_{V_-}(C') - n_{V_-}(C) = \sum_{k \geq 3} (k-3) n_k(C)
\]
such that 
\[
\begin{split}
\mathop{\mathrm{coheight}_{\mathrm{h},\F}}(C) + \mathop{\mathrm{coheight}_{\mathrm{p},\F}}(C) =  \mathop{\mathrm{coheight}_{\mathrm{h},\F}}(C') &+ \mathop{\mathrm{coheight}_{\mathrm{p},\F}}(C') 
\\
& + \sum_{k \geq 3} (k-3) n_k(C)
\end{split}
\]
where $n_{\partial}(C)$ is the number of boundary components of the surface $M$ intersecting $C$, and $n_{\mathrm{fake}}(C)$ is the number of fake prongs in $C$. 
\end{lemma}

\begin{proof}
If the prong connection $C$ consists of one center, then $\mathop{\mathrm{coheight}_{\mathrm{h},\F}}(C) = 0 - 0 =0$ and $\mathop{\mathrm{coheight}_{\mathrm{p},\F}}(C) = 1$, which implies the assertion. 
Thus we may assume that $C$ consists of non-central (i.e. non-zero) prongs. 
Applying the Whitehead moves to $k$-prongs for any $k >0$, we can deform $C$ into a prong connection $C'$ which is homotopy equivalent to $C$ and each of whose prongs is either a $1$-prong and a $3$-prong such that $n_{\partial}(C) = n_{\partial}(C')$ and $n_{\mathrm{fake}}(C) = n_{\mathrm{fake}}(C')$.
Then $b_1(C) = b_1(C')$. 
Put $l := \sum_{k \geq 3} (k-3) n_k(C) = n_{V_-}(C') - n_{V_-}(C) \geq 0$.  
Therefore, we have $\mathop{\mathrm{coheight}_{\mathrm{h},\F}}(C) - \mathop{\mathrm{coheight}_{\mathrm{h},\F}}(C') = (n_{V_-}(C) - b_1(C)) - (n_{V_-}(C') - b_1(C')) = -l$.  
By definition of $\mathop{\mathrm{coheight}_{\mathrm{p},\F}}$, we obtain $\mathop{\mathrm{coheight}_{\mathrm{p},\F}}(C) - \mathop{\mathrm{coheight}_{\mathrm{p},\F}}(C') = 2l$.   
From Lemma~\ref{lem:1_3_prongs}, we have $\mathop{\mathrm{coheight}_{\mathrm{h},\F}}(C') \geq 0$. 
Therefore, we have the following inequality: 
\[
\begin{split}
\mathop{\mathrm{coheight}_{\mathrm{h},\F}}(C) + \mathop{\mathrm{coheight}_{\mathrm{p},\F}}(C) &= \mathop{\mathrm{coheight}_{\mathrm{h},\F}}(C') -l + \mathop{\mathrm{coheight}_{\mathrm{p},\F}}(C') + 2l
\\
&=  \mathop{\mathrm{coheight}_{\mathrm{h},\F}}(C') + \mathop{\mathrm{coheight}_{\mathrm{p},\F}}(C') + l
\\
&\geq 0 + l = \sum_{k \geq 3} (k-3) n_k(C)  \geq 0
\end{split}
\]
\end{proof}

\subsection{Coheight-zero prong connections}\label{sec:sscpc}

We define the types of prong connections as follows. 

\begin{definition}
A prong connection $C$ is $p_o$ if it contains one $1$-prong, one $3$-prong $x$, one heteroclinic separatrix, and one homoclinic separatrix $\gamma$ such that the prongs are $3$-prongs and that there is a closed disk $B$ which does not contain $C$ and that $\partial B =  \{x \} \sqcup \gamma$.
\end{definition}

\begin{definition}
A prong connection $C$ is $p_i$ if it contains one $1$-prong, one $3$-prong $x$, one heteroclinic separatrix, and one homoclinic separatrix $\gamma$ such that the prongs are $3$-prongs and that there is a closed disk $B$ containing $C$ with $\partial B =  \{x \} \sqcup \gamma$.
\end{definition}

\begin{definition}
A prong connection $C$ is $b_\theta$ if it contains two prongs and three heteroclinic separatrices.
\end{definition}

\begin{definition}
A prong connection $C$ is $b_o$ if it contains two prongs, one heteroclinic separatrix $\gamma$, and two homoclinic separatrices such that the prongs are $3$-prongs, the set difference $C - \gamma$ is disconnected,  and there is no closed annulus $A$ containing $C$ with $\partial A =  C - \gamma$.
\end{definition}

\begin{definition}
A prong connection $C$ is $b_i$ if it contains two prongs, one heteroclinic separatrix $\gamma$, and two homoclinic separatrices such that the prongs are $3$-prongs and that there is a closed annulus $A$ containing $C$ with $\partial A =  C - \gamma$.
\end{definition}

\begin{definition}
A prong connection is {\bf semi-self-connected} if it of coheight-zero (i.e. isomorphic to one of Figure~\ref{self_connected_prong} as a topological planar graph, because of Lemma~\ref{lem:ch_codim_zero_01} below). 
\end{definition}

\begin{figure}[t]
\begin{center}
\includegraphics[scale=1.]{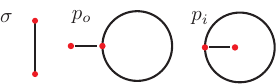}
\includegraphics[scale=1.]{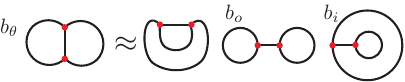}
\end{center} 
\caption{The list of semi-self-connected prong connections.}
\label{self_connected_prong}
\end{figure}



\subsubsection{Characterization of coheight-zero prong connections of levelable foliations}

The definition of coheight of prong connections implies the following observation. 

\begin{lemma}\label{lem:ch_codim_zero}
Let $\F$ be a levelable foliation with finitely many pronged singular points on a surface $M$ and $D(\F)$ the prong connection diagram.  
Then each connected component $C$ of $D(\F)$ is coheight-zero if and only if there are no centers and any prongs in $C$ are $1$-prongs and $3$-prongs outside of the boundary $\partial M$ such that $n_{V_-}(C) = b_1(C)$.
%
\end{lemma}

We have the following characterization of coheight-zero prong connections. 

\begin{proposition}\label{lem:ch_codim_zero_01}
Let $\F$ be a levelable foliation with finitely many pronged singular points on a surface $M$ and $D(\F)$ the prong connection diagram.  
Then each connected component $C$ of $D(\F)$ is of coheight-zero if and only if $C$ is one of Figure~\ref{self_connected_prong} as a topological graph. 
\end{proposition}

\begin{proof}
If a connected component $C$ of $D(\F)$ is one of Figure~\ref{self_connected_prong} as a topological graph, then Lemma~\ref{lem:ch_codim_zero} implies that $C$ is of coheight-zero. 

Conversely, suppose that a connected component $C$ of $D(\F)$ is of coheight-zero. 
By Lemma~\ref{lem:ch_codim_zero}, any prongs in $C$ are $1$-prongs and $3$-prongs outside of the boundary $\partial M$ such that $n_{V_-}(C) = b_1(C)$ because of $\mathop{\mathrm{coheight}_{\mathrm{h},\F}}(C) = 0$. 
If $C$ has no $3$-prongs, then $C$ is of type $\sigma$. 
Thus we may assume that $C$ has $3$-prongs (i.e. $n_{V_-}(C) > 0$). 

\begin{claim}\label{claim:12}
$b_1(C) = 1 + (\vert \{ 3\text{-prong in } C \} \vert - \vert \{ 1\text{-prong in } C \} \vert) / 2 $. 
\end{claim}

\begin{proof}[Proof of Claim~\ref{claim:12}]
By $n_{E}(C) = \sum_{v \in \mathop{\mathrm{Sing}}(\F) \cap C} \deg (v)/2$, we have the following equality: 
\[
\begin{split}
b_1(C) = & 1 - n_{V}(C) + n_{E}(C) 
\\
= & 1 - n_{V}(C) + \sum_{v \in \mathop{\mathrm{Sing}}(\F) \cap C} \deg (v)/2 
\\
= & 1 + \sum_{v \in \mathop{\mathrm{Sing}}(\F) \cap C} \left( \dfrac{\deg (v)}{2} -1 \right) 
\\
= & 1 + (\vert \{ 3\text{-prong in } C \} \vert - \vert \{ 1\text{-prong in } C \} \vert) / 2 
\end{split}
\]
\end{proof}

\begin{claim}\label{claim:11}
$C$ has at most two $3$-prongs 
\end{claim}

\begin{proof}[Proof of Claim~\ref{claim:11}]
Assume $n_{V_-}(C) \geq 3$. 
Then 
\[
\begin{split}
b_1(C) = & 1 + (\vert \{ 3\text{-prong in } C \} \vert - \vert \{ 1\text{-prong in } C \} \vert) / 2 
\\
\leq & 1 + \vert \{ 3\text{-prong in } C \} \vert  / 2  = 1 + n_{V_-}(C)/2 < n_{V_-}(C) = b_1(C)
\end{split}
\]
which is a contradiction. 
\end{proof}

Suppose that $C$ has exactly one $3$-prong $x$. 
Then $b_1(C) = n_{V_-}(C) = 1$. 
From $b_1(C) = 1$, Claim~\ref{claim:12} implies that the prong connection $C$ has exactly one $1$-prong. 
Since all prongs except $x$ are $1$-prongs, each prong separatrix of $C$ either is a separatrix from and to $x$ or connects to a $1$-prong. 
From $b_1(C) = n_{V_-}(C) = 1$, the prong connection $C$ is of type either $p_o$ or $p_i$. 

Suppose that $C$ has at least two $3$-prong $x$. 
By Claim~\ref{claim:11}, we may assume that $C$ has exactly two $3$-prongs $x$ and $y$. 
Then $b_1(C) = n_{V_-}(C) = 2$. 
From $b_1(C) = 2$, Claim~\ref{claim:12} implies that the prong connection $C$ has no $1$-prongs. 
Then each prong separatrix of $C$ is a separatrix from and to the two-point set $\{ x, y \}$. 
If $C$ contains a separatrix from and to a $3$-prong, then $C$ is of type either $b_o$ or $b_i$. 
Thus we may assume that any separatrix connects $x$ and $y$. 
Then $C$ is of type $b_\theta$. 
\end{proof}

The previous proposition implies Theorem~\ref{main:char_0}. 
Moreover, the previous proposition and Theorem~\ref{cor:tree_rep} imply the following statement. 

\begin{theorem}\label{th:tree_rep_codim_zero}
Let $S$ be a compact surface contained in a sphere $\mathbb{S}^2$ and $\F$ a coheight-zero levelable foliation. 
Then the Reeb graph of the height function of $\F$ on $S$ is a finite directed tree whose vertices of degree one are either a prong connection of type $\sigma$ or a periodic leaf on the boundary $\partial \mathbb{S}^2$, and whose vertrice of degree non-one are either $p_o$, $p_i$, $b_\theta$, $b_o$, or $b_i$. 
\end{theorem}

\section{Properties with respect to the topology $\mathcal{O}^p_{w,W,\mathcal{L}(M)}$}

We describe the behaviors near non-periodic leaves to characterize the instability of some kinds of non-periodic points of levelable foliations. 

Recall that the union of periodic leaves of any levelable foliation on a surface is open and the union of non-periodic leaves is closed because of Lemma~\ref{lem:periodic_noncompact} 
Therefore we have the following similar statements to Lemma~\ref{lem:open_reg_pt}, Lemma~\ref{lem:open_reg_pt_03}, and Lemma~\ref{lem:singular_pt_nbd03}. 

\subsection{Properties for levelable singular foliations}

Fix any distance weights $w,W \colon [0,\infty] \to [0,1]$. 
Let $\F$ be a levelable singular foliation with finitely many pronged singular points on a surface $M$ with or without boundary in this subsection. 

\subsubsection{Persistence of the local non-existence of non-periodic leaves}

By replacing $\mathop{\mathrm{Sing}}(\F)$ with $M - \mathop{\mathrm{Per}}(\F)$ in the proofs of Lemma~\ref{lem:open_reg_pt} and Lemma~\ref{lem:open_reg_pt_03}, the same argument in the proofs implies that the local non-existence of non-periodic leaves persists as follows. 

\begin{lemma}\label{lem:open_periodic_pt}
For any $p \in \R_{\geq 1} \sqcup \{ \infty \}$, there are continuous functions $\varepsilon, \delta \colon M \to \R_{\geq 0}$ with $\varepsilon^{-1}(0) = \delta^{-1}(0) = M - \mathop{\mathrm{Per}}(\F)$ satisfying the following equations for any $x_0 \in \mathop{\mathrm{Per}}(\F)$: 
\[
B(x_0,\delta(x_0)) \cap \bigcup_{\mathcal{G} \in B_{D^{p}_{W_{\mathrm{per}}}}(\F,\varepsilon(x_0))}  M - \mathop{\mathrm{Per}}(\mathcal{G}) = \emptyset
\] 
and so 
\[
B(x_0,\delta(x_0)) \cap \bigcup_{\mathcal{G} \in B_{D^{p}_{w,W}}(\F,\varepsilon(x_0))}  M - \mathop{\mathrm{Per}}(\mathcal{G}) = \emptyset
\] 
\end{lemma}

\subsubsection{Descriptions near non-periodic points}

By replacing $\mathop{\mathrm{Sing}}(\F)$ with $M - \mathop{\mathrm{Per}}(\F)$ in the proof of Lemma~\ref{lem:singular_pt_nbd03}, the same argument in the proof implies the following statements. 

\begin{lemma}\label{lem:non-periodic_pt_nbd}
For any positive number $\varepsilon >0$, there is a positive number $\delta > 0$ satisfying the following inequality: 
\[ 
\sup_{x_0 \in M - \mathop{\mathrm{Per}}(\F)} \sup_{x \in B(x_0,\delta)}d_H(X_{w,\F}(x), \{ 0_{x} \}) \leq \varepsilon
\]
\end{lemma}

\subsubsection{Persistence of periodic leaves of a levelable singular foliation on a surface}

Lemma~\ref{lem:open_periodic_pt} implies the following statements. 

\begin{lemma}\label{lem:per_persistence}
For any $p \in \R_{\geq 1} \sqcup \{ \infty \}$, any periodic leaf $L \in \F$, an open transverse arc $T \subset M$ intersecting $L$ at a point $x$, and any positive number $\varepsilon>0$, there are a positive number $\delta > 0$ and closed sub-arcs $T',T'',T''' \subset T$ with $x \in T''' \subseteq T'' \subseteq T' \subseteq T \subset B_\varepsilon (T')$ and $\F(T'') \subseteq B_\varepsilon(L)$  satisfying the following statements for any foliation $\mathcal{G} \in B_{D^{p}_{w,W}}(\F,\delta)$: 
\\
{\rm(1)} $T' \subset \mathop{\mathrm{Per}}(\mathcal{\mathcal{G}})$. 
\\
{\rm(2)} If $L$ is transversely orientable {\rm(resp.} transversely non-orientable{\rm)}, then the saturation $\F(T'')$ is an periodic annulus {\rm(resp.} M{\"o}bius band{\rm)} with respect to $\F$ and the saturation $\G(T''')$ is an periodic annulus {\rm(resp.} M{\"o}bius band{\rm)} with respect to $\G$ such that $\F(x) \subset \G(T''') \subseteq \F(T'')$. 
\\
{\rm(3)} There is a point $y \in T'''$ such that the leaf $\G(y)$ is isotopic to $\F(x)$ in $\G(T''')$. 
\end{lemma}

We recall the following concepts. 
\begin{definition}
A leaf $L$ of $\F$ is {\bf two-sided} if it has an open small neighborhood $U$ which is an annulus such that $U - L$ is a disjoint union of two open annuli. 
\end{definition}

\begin{definition}
A leaf of $\F$ is {\bf one-sided} if it is not two-sided. 
\end{definition}

\begin{definition}
The {\bf leaf space} $M/\F$ is the quotient space $M/\mathop{\sim_{\F}}$ by $x \sim_{\F} y$ if $\F(x) = \F(y)$. 
\end{definition}

Notice that any periodic leaves on the boundary are one-sided. 
Moreover, for any periodic leaf $L$ such that any small neighborhood of $L$ contains a neighborhood which is a M{\"o}bius band, the leaf $L$ is one-sided. 
We have the following properties. 

\begin{lemma}\label{lem:per_persistence_03}
For any $p \in \R_{\geq 1} \sqcup \{ \infty \}$ and any periodic leaf $L \in \F$, there are a positive number $\delta > 0$ satisfying the following property for any foliation $\mathcal{G} \in B_{D^{p}_{w,W}}(\F,\delta)$: 
\\
{\rm(1)} If $L$ is transversely orientable, then there is a closed periodic annulus with respect to $\G$ whose interior contains $L$ and which is contained in a closed periodic annulus with respect to $\F$ in $B_\varepsilon(L)$. 
\\
{\rm(2)} If $L$ is transversely non-orientable, then there is a closed periodic M{\"o}bius band with respect to $\G$ whose interior as a subset of $M$ contains $L$ and which is contained in a closed periodic M{\"o}bius band with respect to $\F$ in $B_\varepsilon(L)$. 
\end{lemma}

\begin{proof}
Let $L \in \F$  be a periodic leaf. 
By the compactness of $L$, Lemma~\ref{lem:open_periodic_pt} implies that the constants $\varepsilon := \min_{x \in L} \varepsilon(x)$ and $\delta := \min_{x \in L} \delta(x)$ are positive, and that there are finite points $x_1, \ldots, x_k \in L$ with $L \subset \bigcup_{i=1}^k B(x_i, \delta)$ satisfying the following equality: 
\[
\left( \bigcup_{i=1}^k B(x_i,\delta) \right) \cap \bigcup_{\mathcal{G} \in B_{D^{p}_{W_{\mathrm{per}}}}(\F,\varepsilon)}  M - \mathop{\mathrm{Per}}(\mathcal{G}) = \emptyset
\] 
Then there is a positve number $\delta' \in (0, \delta)$ such that $B(L,\delta') \subset \bigcup_{i=1}^k B(x_i, \delta)$. 
Therefore we have 
\[
B(L,\delta') \cap \bigcup_{\mathcal{G} \in B_{D^{p}_{W_{\mathrm{per}}}}(\F,\varepsilon)}  M - \mathop{\mathrm{Per}}(\mathcal{G}) = \emptyset
\] 
and so the following relation holds: 
\[
B(L,\delta') \subseteq  \bigcap_{\mathcal{G} \in B_{D^{p}_{W_{\mathrm{per}}}}(\F,\varepsilon)}   \mathop{\mathrm{Per}}(\mathcal{G})
\] 
Fix an open transverse arc $T$ intersecting $x$ such that $\F(T) \subset B_\varepsilon(L)$ is a periodic annulus with respect to $\F$. 
Lemma~\ref{lem:per_persistence} implies that there is a closed sub-arcs $T''',T'' \subset T$ with $x \in T''' \subseteq T''$ such that 
$\F(x) \subset \G(T''') \subset \F(T'') \subset B_\varepsilon(L)$. 
Moreover, if $L$ is transversely orientable (resp. non-orientable), then we may assume that $\G(T''')$ is a closed periodic annulus (resp. M{\"o}bius band) with respect to $\G$ whose interior contains $L$ and which is contained in the closed periodic annulus  (resp. M{\"o}bius band) $\F(T'')$ with respect to $\F$ in $B_\varepsilon(L)$. 
Therefore, the assertion holds. 
\end{proof}

\begin{lemma}\label{lem:per_persistence_04}
Let $G_{H_{\F}}$ be the Reeb graph of the height function $H_{\F}$ of the levelable singular foliation $\F$, which is an acyclic finite directed graph. 
There is a set $\{ L_1, \ldots, L_k \}$ of distinct periodic leaves contained in the interiors of edges of $G_{H_{\F}}$ satisfying the following properties: 
\\
{\rm(1)} Every connected component of $G_{H_{\F}} - \{ L_1, \ldots, L_k \}$ contains exactly one vertex of the end-completion of $G_{H_{\F}}$. 
\\
{\rm(2)} There is a positive number $\delta > 0$ such that, for any foliation $\mathcal{G} \in B_{D^{p}_{w,W}}(\F,\delta)$, there are $\mathcal{G}$-invariant closed periodic annuli $\mathbb{A}_i$ which contains $L_i$ for any $i \in \{ 1, \ldots, k\}$ such that $V_{\G,j} \subseteq V_{\F,j}$ for any $j \in \{ 1, \ldots, l \}$, where $\{ V_{\F,1}, \ldots, V_{\F,l} \}$ is the set of connected components of $M - \bigsqcup_{i=1}^k L_i$ and $V_{\G,j}$ is the connected component of $M - \bigsqcup_{i=1}^k \mathbb{A}_i$ intersecting $V_{\F,j}$.
\end{lemma}

\begin{proof}
For any edge of the Reeb graph $G_{H_{\F}}$, choose exactly one periodic leaf in the interior of the edge. 
The set of such periodic leaves is denoted by $\{ L_1, \ldots, L_k \}$, which implies assertion (1). 
By Lemma~\ref{lem:per_persistence_03}, there are $\mathcal{G}$-invariant closed periodic annuli $\mathbb{A}_i$ which contains $L_i$ for any $i \in \{ 1, \ldots, k\}$. 
Let $V_{\G,j}$ be the connected component of $M - \bigsqcup_{i=1}^k \mathbb{A}_ i$ intersecting $V_{\F,j}$. 
By construction, we have $V_{\G,j} \subseteq V_{\F,j}$ for any $j \in \{ 1, \ldots, l \}$. 
\end{proof}

\subsection{Local maximality of coheights}

We have the following observation.

\begin{lemma}\label{lem:codim_max}
Let $M$ be a surface. 
For any $\F \in \mathcal{L}^*_{i}(M)$, there is a positive number $\delta > 0$, the following statements hold for any $\G \in \mathcal{L}^*_{i}(M)$ near $\F$:
\\
{\rm(1)} $n_{\partial}(\G) \leq n_{\partial}(\F)$. 
\\
{\rm(2)} $n_{V_-}(\G) - n_{E}(\G) + n_{0}(\G) = n_{V_-}(\F) - n_{E}(\F) + n_{0}(\G)$, where $n_{0}(\F')$ is the number of centers of $\F'$. 
\\
{\rm(3)} $- n_{CC}(\G) \leq - n_{CC}(\F)$. 
\\
{\rm(4)} $\sum_{y \in \mathop{\mathrm{Sing}}(\G)} (1 + \mathrm{coheight}_{\mathrm{p},\G}(y)) \leq \sum_{y \in \mathop{\mathrm{Sing}}(\F)} (1 + \mathrm{coheight}_{\mathrm{p},\F}(y))$. 
\\
{\rm(5)} $\mathrm{coheight}(\G) \leq \mathrm{coheight}(\F)$.  
\end{lemma}

\begin{proof}
Fix a line field $\G \in \mathcal{L}^*_{i}(M)$ near $\F$. 
By Proposition~\ref{prop:inv_index}, assertion (1) holds and $n_{V_-}(\G) \geq n_{V_-}(\F)$. 
Lemma~\ref{lem:per_persistence_04} implies assertion (3).

From the invariance of the sum of positive indices of singular leaves, the sum of negative indices of singular leaves is invariant in $\mathcal{L}^*_{i}(M)$. 
This implies that any $1$-prong is invariant under small perturbations and that any center either is preserved or becomes two $1$-prongs under small perturbations. 
The invariance of the sum of positive indices of singular leaves implies assertion (2). 

Notice that $\sum_{y \in \mathop{\mathrm{Sing}}(\F)} (1 + \mathrm{coheight}_{\mathrm{p},\F}(y))$ is invariant under  any small perturbations near $p$, which either preserves $p$ or split $p$ into two $1$-prongs.  
From the definition of $\mathrm{coheight}_{\mathrm{p},\F}$, by the invariance of the sum of negative indices of singular leaves, Lemma~\ref{lem:inv_index_prong} implies 
 \[
 \begin{split}
 \sum_{y \in \mathop{\mathrm{Sing}}(\G)} (1 + \mathrm{coheight}_{\mathrm{p},\G}(y)) &\leq 
n_{V_-}(\G) - n_{V_-}(\F) + \sum_{y \in \mathop{\mathrm{Sing}}(\G)} (1 + \mathrm{coheight}_{\mathrm{p},\G}(y)) 
\\
&= \sum_{y \in \mathop{\mathrm{Sing}}(\F)} (1 + \mathrm{coheight}_{\mathrm{p},\F}(y))
 \end{split}
\]
because of $0 \leq n_{V_-}(\G) - n_{V_-}(\F)$. 
By assertion (1)--(4), the non-existence of fake prongs implies assertion (5). 
\end{proof}

\subsection{Correspondence between structural stabilities for Hamiltonian flows and levelable foliations on surfaces}

%

Fix any distance weight $W \colon [0,\infty] \to [0,1]$.

\subsubsection{Whitehead moves}

\begin{figure}
\begin{center}
\includegraphics[scale=0.6]{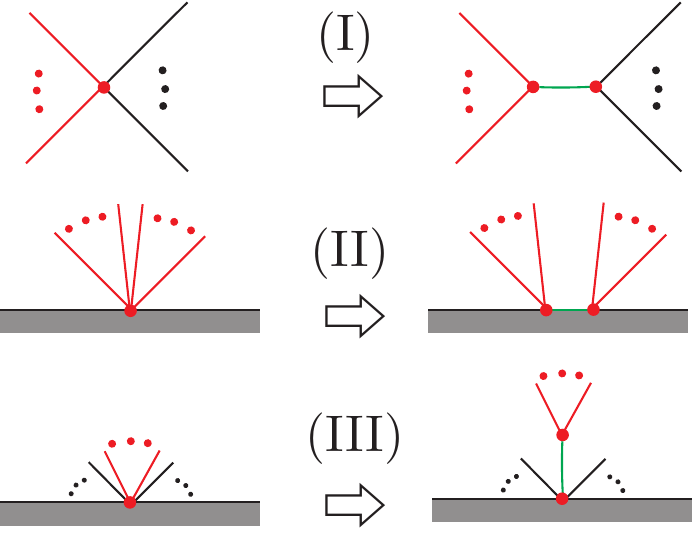}
\end{center}
\caption{Whitehead moves}
\label{splitting}
\end{figure}

Each collapsing of a heteroclinic separatrix and the inverse operation as in Figure~\ref{splitting} is called the {\bf Whitehead move}. 
Roughly speaking, the Whitehead moves are ``splitting of prongs'' and the inverse operation of the Whitehead moves are ``merges of prongs''.

\subsubsection{Density of levelable foliations}

Denote by $\bm{\mathcal{L}^*_{i}(M)} \subset \mathcal{L}^*(M)$ the subspace of levelable foliations whose sums of indices of pronged singular points with positive indices are $i$. 
The subspace $\mathcal{L}^*_{i}(M)$ is considered as the space without creations and annihilations. 

Denote by $\bm{\mathcal{L}^{0}(M)} \subset \mathcal{L}^*(M)$ the subspace of levelable foliations whose prong connection diagrams are semi-self-connected.  
Put $\mathcal{L}^0_{i}(M) := \mathcal{L}^0(M) \cap \mathcal{L}^*_{i}(M)$. 
We will show that $\mathcal{L}^0_{i}(M)$ for any compact surface $M$ are of coheight-zero and consist of structurally stable ones (see Theorem~\ref{th:stability03}). 
%
%
The Whitehead move implies the following observation. 

\begin{lemma}\label{lem:stability04}
For  any $p \in \R_{\geq 1} \sqcup \{ \infty \}$ and any $i \in \Z_{>1}$, the following statements hold with respect to the topology $\mathcal{O}^p_{w,W,\mathcal{L}(M)}$ and the relative topologies:   
\\
{\rm(1)} The subspace in $\mathcal{L}^\star(M)$ which consists of levelable foliations each of whose prongs is either a $1$-prong {\rm(}thorn{\rm)}, a $2$-prong {\rm(}fake prong{\rm)}, or a $3$-prong {\rm(}tripod{\rm)} is dense in $\mathcal{L}^\star(M)$. 
\\
{\rm(2)} The subspace in $\mathcal{L}^*(M)$ which consists of levelable foliations each of whose prongs is either a $1$-prong or a $3$-prong is dense in $\mathcal{L}^*(M)$. 
\end{lemma}

\begin{proof}
Replacing a small \nbd of singular points with the resulting \nbd of Whitehead moves, any foliations in $\mathcal{L}^\star(M)$ (resp. $\mathcal{L}^*(M)$) can be obtained from levelable foliations whose prongs are $1$-prongs, $2$-prongs, and $3$-prongs (resp. $1$-prongs and $3$-prongs) by arbitrarily small perturbations. 
This implies the assertion. 
\end{proof}

\begin{lemma}\label{lem:stability03}
For any $i \in \Z_{>1}$, the subspace $\mathcal{L}^{0}(M)$ is dense in $\mathcal{L}^*(M)$ with respect to the relative topology of the topology $\mathcal{O}^p_{w,W,\mathcal{L}(M)}$. 
\end{lemma}

\begin{proof}
Fix a lelvalable foliation $\F \in \mathcal{L}^{*}_{i}(M)$.  
By Lemma~\ref{lem:stability04}, we may assume that any prongs of $\F$ are $1$-prongs and $3$-prongs. 
Let $C$ be a prong connection. 

\begin{claim}\label{claim:04}
We may assume there is no $1$-prong in $C$. 
\end{claim}
\begin{proof}[Proof of Claim~\ref{claim:04}]
Suppose that there is a $1$-prong $x$ in $C$.
If $C$ is $\sigma$, then $C$ is semi-self-connected. 
Thus we may assume that $C$ is not $\sigma$. 
Then $C$ contains $3$-prongs and so a separatrix between a $1$-prong and a $3$-prong. 
A replacement of two fibered charts near $x$ into a fibered chart as in the upper on Figure~\ref{prong_perturbation02} preserves the outside of an arbitrarily small \nbd of $C$ and so is an arbitrarily small perturbation with respect to any locally weighted Hausdorff distance with respect to $w$ and $W$. 
\begin{figure}[t]
\begin{center}
\includegraphics[scale=1.]{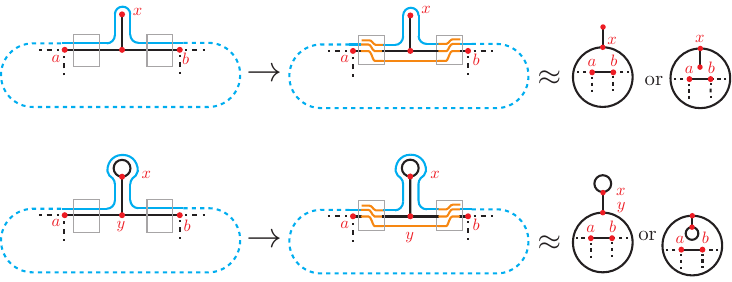}
\end{center} 
\caption{Replacement of two fibered charts near $3$-prongs with three heteroclinic separatrices. (Note that $a = b$ possibly.)}
\label{prong_perturbation02}
\end{figure}
This completes the claim since the resulting prong connection containing $x$ is either $p_o$ or $p_i$. 
\end{proof}

\begin{claim}\label{claim:05}
We may assume there are no self-connected separatrices in $C$. 
\end{claim}
\begin{proof}[Proof of Claim~\ref{claim:05}]
Suppose that there is a prong $x$ with a self-connected separatrix in $C$. 
Then there is a separatrix between $x$ and another prong $y$. 
A replacement of two fibered charts near $y$ into a fibered chart as in the lower on Figure~\ref{prong_perturbation02} preserves the outside of an arbitrarily small \nbd of $C$ and so is an arbitrarily small perturbation with respect to any locally weighted Hausdorff distance with respect to $w$ and $W$. 
This completes the claim since the resulting prong connection containing $x$ is either $b_o$ or $b_i$. 
\end{proof}

\begin{claim}\label{claim:06}
We may assume there are no heteroclinic separatrices in prong connections that are not $b_\theta$. 
\end{claim}
\begin{proof}[Proof of Claim~\ref{claim:06}]
Assume that there is a heteroclinic separatrix $\gamma$ in a prong connection that is not $b_\theta$. 
A replacement of four fibered charts near $\gamma$ into a fibered chart as in Figure~\ref{prong_perturbation03} preserves the outside of an arbitrarily small \nbd of $C$ and so is an arbitrarily small perturbation with respect to any locally weighted Hausdorff distance with respect to $w$ and $W$. 
\begin{figure}[t]
\begin{center}
\includegraphics[scale=1.]{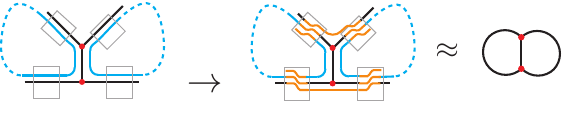}
\end{center} 
\caption{Replacement of four fibered charts near a heteroclinic separatrix.}
\label{prong_perturbation03}
\end{figure}
Since the resulting prong connection containing $\gamma$ is $b_\theta$, this completes the claim. 
\end{proof}

The claims imply that $C$ is $b_\theta$.
This completes the proof. 
\end{proof}

\begin{lemma}\label{lem:stability04+}
Suppose that $M$ is compact. 
For any $i \in \Z_{>1}$,  the subspace $\mathcal{L}^{0}_{i}(M)$ is open dense in the space $\mathcal{L}^*_{i}(M)$ with respect to the relative topology of the topology $\mathcal{O}^p_{w,W,\mathcal{L}(M)}$. 
\end{lemma}

\begin{proof}
Lemma~\ref{lem:stability03} implies the density. 
Without loss of generality, we may assume that $M$ is connected.
Let $\F \in \mathcal{L}^{0}_{i}(M)$ be a levelable foliation on $M$. 
Then, any prong connections of $\F$ consist of prong separatrices and two prongs, which are $1$-prongs and $3$-prongs. 
By Lemma~\ref{lem:periodic}, the complement of the prong connection diagram of $\F$ is a finite union of annuli which is an open subset of $M$ and is the union of periodic leaves. 
From Lemma~\ref{lem:per_persistence_04}, there are periodic leaves $L_1 , \ldots, L_k$ and $V_{\F,1}, \ldots, V_{\F,l}$ connected components of $M - \bigsqcup_{j=1}^k L_j$ as in Lemma~\ref{lem:per_persistence_04}. 
Then all connected components $V_{\F,j}$ are disks, annuli, and two-punctured disks. 
By Poincar{\'e}-Hopf theorem for one-dimensional singular foliations with finitely many singular points on a surface \cite[p.113 Thoerem~2.2]{Hopf1989diff}, if such a connected component $V_{\F,j}$ is a disk (resp. annulus, two-punctured disk) then it contains exactly two singular leaves which consists of two $1$-prongs (resp, one $1$-prong and one $3$-prong, two $3$-prongs).
Let $i_{\F,j'} \in \{ 0,1\}$ be the sum of positive indices of singular leaves of $\F$ in $V_{\F,j'}$. 
Then $i = \sum_{j'=1}^l i_{\F,j'}$. 
Lemma~\ref{lem:per_persistence_04} implies that there is a positive number $\delta > 0$ such that, for any foliation $\mathcal{G} \in B_{D^{p}_{w,W}}(\F,\delta)$, there are $\mathcal{G}$-invariant closed periodic annuli $\mathbb{A}_j$ which contains $L_j$ for any $j \in \{ 1, \ldots, k\}$ such that $V_{\G,j'} \subseteq V_{\F,j'}$ for any $j' \in \{ 1, \ldots, l \}$, where $V_{\G,j'}$ is the connected component of $M - \bigsqcup_{j=1}^k \mathbb{A}_{j'}$ intersecting $V_{\F,j}$. 
This implies that all connected components $V_{\G,j'}$ are disks, annuli, and two-punctured disks. 

Let $i_{\G,j'} \in \Z_{\geq 0}$ be the sum of positive indices of singular leaves of $\G$ in $V_{\G,j'} \subset V_{\F,j'}$. 
Lemma~\ref{lem:singular_pt_nbd} implies that $i_{\G,j'} \geq i_{\F,j'}$ for any $j' \in \{ 1, \ldots, l \}$, and that the numbers $n_{\G,j'}$ contained in singular leaves of $\G$ in $V_{\G,j'}$ is more than or equal to those $n_{\F,j'}$ of $\F$ in $V_{\G,j'}$.
By $\sum_{j'=1}^l i_{\F,j'} = i = \sum_{j'=1}^l i_{\G,j'}$, we obtain $i_{\G,j'} = i_{\F,j'}$ and $n_{\F,j'} = n_{\G,j'}$ for any $j' \in \{ 1, \ldots, l \}$. 
Since the indices of $k$-prongs are $1 - k/2$, Poincar{\'e}-Hopf theorem implies that the indices of any singular leaves of $\G$ are $1/2$ and $-1/2$. 
This means that all singular leaves of $\G$ are $1$-prongs and $3$-prongs. 
Therefore each connected component $V_{\G,j}$ which is a disk (resp. annulus, two-punctured disk) contains exactly two singular leaves which consist of two $1$-prongs (resp, one $1$-prong and one $3$-prong, two $3$-prongs).
All possible prong connections with exactly two prongs which are $1$-prongs and $3$-prongs are semi-self-connected as shown in Figure~\ref{self_connected_prong}. 
This means that $\mathcal{L}^{0}_{i}(M)$ is open. 
\end{proof}

\subsubsection{Structural stability of levelable foliations}

To demonstrate the structural stability, we define the canonical local structures as follows.

\begin{figure}[t]
\begin{center}
\includegraphics[scale=.5]{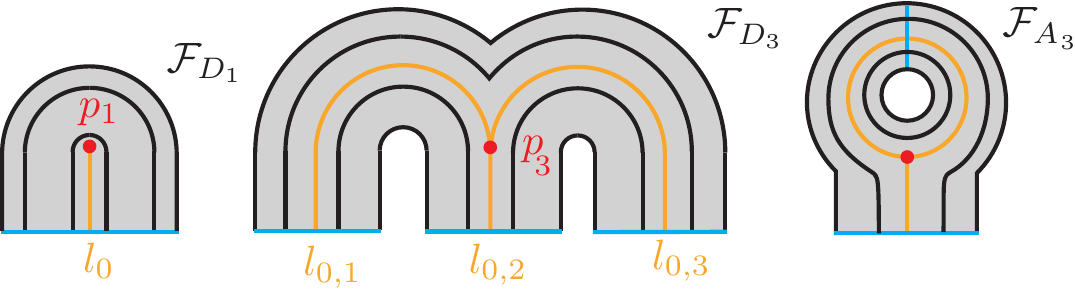}
\end{center} 
\caption{The partition $\F_{D_1}$ on the disk $D_1$, the partition $\F_{D_3}$ on the disk $D_3$, and the partition $\F_{D_A}$ on the annulus $A_3$.}
\label{fig:stdd_loc_str}
\end{figure} 

\begin{definition}
Set $D_1 := \{ (x,y) \mid x \in [-1,1], y \in [0, \sqrt{1 - x^2} + 1] \}$ $p_1 := (0,1)$, and $l_0 := \{0\} \times [0,1) \subset D_1$. 
Define 
\[
l_r := \{ p \in D_1 \mid d(p, \{p_1 \} \sqcup l_0) = r \}
\]
for any $r \in (0,1]$. 
Then $\{ \{ p_1 \} \} \sqcup \{ l_r \mid r \in [0,1] \}$ is a partition and denote by $\F_{D_1}$ as in the left in Figure~\ref{fig:stdd_loc_str}. 
Moreover, the transverse arc $[-1,1] \times \{0 \} \subset D_1$ is called the {\bf core transverse arc} of $D_1$, and the periodic orbit intersecting the boundary of the core transverse arc is called the {\bf core periodic orbit} for $D_1$. 
\end{definition}

\begin{definition}
Set $p_3 := (0,1)$, $l_{0,2} := \{0 \} \times [0,1)$, 
\[
\begin{split}
l_{0,1} &:= (\{-1\} \times [0,1)) \sqcup \{(x,y) \mid x \in [-1,0), y = \sqrt{1/4 - (x+1/2)^2} + 1 \}
\\
l_{0,3} &:= (\{1\} \times [0,1)) \sqcup \{(x,y) \mid x \in (0,1], y = \sqrt{1/4 - (x-1/2)^2} + 1 \}
\end{split}
\]
and $D_3 := (\R \times \R_{\geq 0}) \cap B_{1/3}(l_0)$. 
Define 
\[
l_r := \{ p \in D_3 \mid d(p, \{p_3 \} \sqcup l_{0,1} \sqcup l_{0,2} \sqcup l_{0,3}) = r \}
\]
for any $r \in (0,1]$. 
Then $l_r$ for any $r \in (0,1]$ consists of three closed arcs $l_{r,1},l_{r,2}$, and $l_{r,3}$. 
Then $\{\{p_3\} \} \sqcup \{ l_{r,i} \mid r \in [0,1/3], i \in \{1,2,3\} \}$ is a partition and denote by $\F_{D_3}$ as in the middle in Figure~\ref{fig:stdd_loc_str}. 
Moreover, the transverse arcs $[-4/3,-2/3] \times \{0 \}, [-1/3,1/3] \times \{0 \}, [2/3,4/3] \times \{0 \} \subset D_3$ are called the {\bf core transverse arcs} of $D_3$ and the periodic orbits intersecting the boundaries of the core transverse arcs are called the {\bf core periodic orbit} for $D_3$. 
\end{definition}

\begin{definition}\label{def:disk2anuulus}
Consider the partition $\F_{D_3}$ on $D_3$. 
Set $T_- := [-4/3,-2/3] \times \{0 \}$, $T_0 := [-1/3,1/3] \times \{0 \}$, and $T_+ := [2/3,4/3] \times \{0 \}$.
Then $D_3 \cap (\R \times \{ 0 \} = T_- \sqcup T_0 \sqcup T_+$. 
Define an equivalence relation $\sim_{A_3}$ on $D_3$ as follows: 
\[
\begin{split}
(x,y) \sim_{A_3} (x',y') \text{ if } & \text{ either } (x,y) = (x',y') 
\\
& \text{ or } \left( x+x' = y = y' = 0 \text{ and } (x,y) \in T_- \sqcup T_+ \right)
\end{split}
\]
The quotient space $A_3 := D_3/\mathop{\sim_{A_3}}$ is a closed annulus, and the induced partition is denoted by $\F_{A_3}$ as in the right in Figure~\ref{fig:stdd_loc_str}. 
Moreover, the transverse arcs $\pi([-4/3,-2/3] \times \{0 \}), \pi([-1/3,1/3] \times \{0 \}) \subset A_3$ is called the {\bf core transverse arcs} of $A_3$ and the periodic orbits intersecting the boundaries of the core transverse arcs are called the {\bf core periodic orbit} for $A_3$, where $\pi \colon D_3 \to A_3$ is the quotient map. 
\end{definition}

Using these partitions, we define the following types. 

\begin{definition}
A partition on a disk $D$ is of {\bf type $\bm{\sigma/2}$} if it is locally topologically equivalent to the partition $\F_{D_1}$.
Then $D$ is called a disk of  {\bf type $\bm{\sigma/2}$}. 
\end{definition}

\begin{definition}
A partition on a disk $D$ is of {\bf type $\bm{b_{o/i}/2}$} if it is locally topologicaly equivalent to the partition $\F_{D_1}$. 
Then the disk $D$ is called a domain of  {\bf type $\bm{b_{o/i}/2}$}. 
\end{definition}

\begin{definition}
A partition on an annulus $A$ is of  {\bf type $\bm{b_\theta/2}$} if it is locally topologically equivalent to the partition $\F_{A_1}$. 
Then the annulus $A$ is called a domain of  {\bf type $\bm{b_\theta/2}$}. 
\end{definition}

We characterize the local structures as follows.

\begin{lemma}\label{lem:char_local_prongs}
Let $\F$ be a levelable foliation with finitely many pronged singular points on a surface $M$ and $D \subset M$ a disk which contains no prong except one $1$-prong $q_1$ and whose boundary consists of one closed leaf arc $C$ and one closed transverse arc $T$ such that the intersection $C \cap T = \partial C = \partial T$ consists of two points. 
Then the restriction $\F\vert_{D}$ is of type $\sigma/2$.
\end{lemma}

\begin{figure}[t]
\begin{center}
\includegraphics[scale=.5]{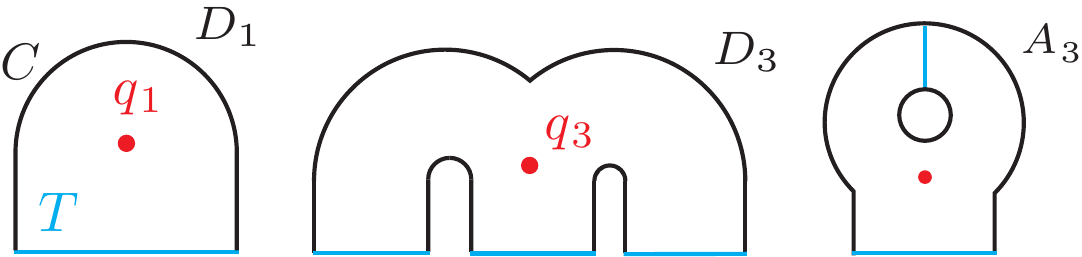}
\end{center} 
\caption{A disk $D_1$ which contains no prong except one $1$-prong and whose boundary consists of one closed leaf arc $C$ and one closed transverse arc $T$, a disk $D_3$ which contains no prong except one $3$-prong and whose boundary consists of three closed leaf arcs and three closed transverse arcs, and an anulus $A_3$ which contains no prong except one $3$-prong and whose boundary consists of two boundary components such that one of them consists of one closed leaf arc and one closed transverse arc and another is a periodic leaf.}
\label{fig:semi_disk}
\end{figure} 

\begin{proof}
Since any leaf of the levelable foliation $\F$ is closed in $M - \mathop{\mathrm{Sing}}(\F)$, the semi-prong separatrix $L$ of $q_1$ intersects the interior of $T$. 
Let $l \subset D$ be the closed leaf arc between $q_1$ and the point $q \in L \cap T$, which is the first hitting point to $T$ from $q_1$ by $L$ as on the left in Figure~\ref{fig:semi_disk}. 
Then the complement $T - \{ q\}$ is the disjoint union of two intervals $T_-$ and $T_+$. 

\begin{claim}\label{claim:17}
For any point $p \in D - (\{q_1\} \sqcup l)$, the element $\F\vert_D(p)$ of the restriction $\F\vert_D$ is closed. 
\end{claim}

\begin{proof}[Proof of Claim~\ref{claim:17}]
Fix any point $p \in D - (\{q_1\} \sqcup l)$. 
By Lemma~\ref{lem:regularity_proper}, for any point $p \in D - (\{q_1\} \sqcup l)$, the set difference $\overline{\F(p)} - \F(p)$ consists of singular points and so $(\overline{\F(p)} - \F(p)) \cap D \subset \{ q_1 \}$.  
Assume that $q_1 \in \overline{\F\vert_D(p)}$. 
Since $D$ is closed, by $\overline{\F\vert_D(p)} \subset \overline{\F(p)} \cap D$, we have $\overline{\F\vert_D(p)} = \F\vert_D(p) \sqcup \{ q_1 \}$. 
Then there is a connected component $\gamma$ of $\F\vert_D(p) - \{ p \} \subset D - l$ which is an open arc whose closure is a closed arc $\gamma \sqcup \{p, q_1\}$. 
This means that $q_1$ is not a $1$-prong because $l \cap \gamma = \emptyset$ and $\overline{l} \cap \overline{\gamma} = q_1$. 
Thus $q_1 \notin \overline{\F\vert_D(p)}$ and so $(\overline{\F\vert_D(p)} - \F\vert_D(p)) \cap D = \emptyset$, which implies the assertion.  
\end{proof}
If there is a periodic leaf in $D$, then Poincar{\'e}-Hopf theorem implies that $D$ contains a singular point which is not $q_1$, which contradicts the hypothesis. 
Thus any elements of $\F\vert_D(p)$ for any point $p \in D - (\{q_1\} \sqcup l)$ are closed intervals between $T$. 
Since $T$ has no tangencies, each element of $\F\vert_D(p)$ for any point $p \in D - (\{q_1\} \sqcup l)$ connect a point in $T_-$ and a point in $T_+$. 
This means that the assertion holds. 
\end{proof}

\begin{lemma}\label{lem:char_local_prongs_02}
Let $\F$ be a levelable foliation with finitely many pronged singular points on a surface $M$ and $D \subset M$ a disk which contains no prong except one $3$-prong $q_3$ and whose boundary consists of pairwise disjoint three closed leaf arc $C_1,C_2,C_3$ and pairwise disjoint three closed transverse arcs $T_1,T_2,T_3$ such that $(C_i \sqcup C_{i+1}) \cap T_i = \partial T_i$ and  $(T_i \sqcup T_{i+1}) \cap T_{i+1} = \partial T_i$ for any $i \in \{1,2,3\}$, where $C_4 := C_1$ and $T_4 := T_1$. 
Then the restriction $\F\vert_{D}$ is of type $b_\theta/2$. 
\end{lemma}

\begin{proof}
Let $D$ be a disk as above, as on the right in Figure~\ref{fig:semi_disk}. 
Assume there is a prong separatrix of $q_3$ in $D$. 
Then the disk in $D$ whose boundary consists of the prong separatrix and $q_3$ contains a disk whose boundary is a periodic leaf. 
Therefore, Poincar{\'e}-Hopf theorem implies that $D$ contains a singular point which is not $q_1$, which contradicts the hypothesis. 

Thus there are no prong separatrix of $q_3$ in $D$.
Let $l_1,l_2,l_3$ be the elements of the restriction $\F\vert_D$ each of whose closures contains $q_3$. 
Since any leaf of the levelable foliation $\F$ is closed in $M - \mathop{\mathrm{Sing}}(\F)$, each semi-prong separatrix of $q_3$ intersects $\mathrm{int} T_i$ for some $i \in \{ 1,2,3 \}$. 
Because $T_1 \sqcup T_2 \sqcup T_3$ has no tangencies, the three elements $l_1,l_2$, and $l_3$ intersects distinct connected components of $T_1 \sqcup T_2 \sqcup T_3$. 
By renaming $T_i$ if neccesary, we may assume that $T_i$ intersects $l_i$ for any $i \in \{ 1,2,3 \}$. 
As the same argument in the proof of Claim~\ref{claim:17} and the subsequent explanation,  for any point $p \in D - (\{q_3\} \sqcup l_1 \sqcup l_2 \sqcup l_3)$, the element $\F\vert_D(p)$ of the restriction $\F\vert_D$ is a closed intervals between $T_1 \sqcup T_2 \sqcup T_3$. 
Since each $T_i$ has no tangencies, each element of $\F\vert_D(p)$ for any point $p \in D - (\{q_3\} \sqcup l_1 \sqcup l_2 \sqcup l_3)$ connect points in different connected components of $T_1 \sqcup T_2 \sqcup T_3$. 
This means that the assertion holds. 
\end{proof}

\begin{lemma}\label{lem:char_local_prongs_03}
Let $\F$ be a levelable foliation with finitely many pronged singular points on a surface $M$ and $A \subset M$ an annulus which contains no prong except one $3$-prong $q_3$ and whose boundary consists of two boundary components $\partial_1$ and $\partial_2$ such that $\partial_1$ consists of one closed leaf arc $C$ and one closed transverse arc $T$ with $C \cap T = \partial C = \partial T$ and that $\partial_2$ is a periodic leaf. 
Suppose that there is a closed transverse between $\partial_1$ and $\partial_2$. 
Then the restriction $\F\vert_{A}$ is type $b_{o/i}/2$.
\end{lemma}

\begin{proof}
Cutting the closed transverse arc $T$ into two closed intervals $T_-$ and $T_+$, the resulting space $D$ is a disk as in Lemma~\ref{lem:char_local_prongs_02}. 
By Lemma~\ref{lem:char_local_prongs_02}, the resulting partition $\F_D$ of $\F$ is locally topologicaly equivalent to the partition $\F_{D_3}$. 
The existence of the height function $H_\F$ of $\F$ implies that the restriction $H_\F\vert_T$ on the interval is strictly monotonic. 
Since $\F$ can be considered with the resulting partition on the quotient space as in Definition~\ref{def:disk2anuulus}, the restriction $\F\vert_{A}$ is locally topologically equivalent to the partition $\F_{A_3}$.
\end{proof}

For any element $\F$ of $\mathcal{L}^{0}(M)$ and for any prong connection $C$ of $\F$, by Lemmas~\ref{lem:char_local_prongs}--\ref{lem:char_local_prongs_03}, we can choose an invariant closed \nbd $D_- \sqcup_T D_+$ of $C$ which is a union of two domains $D_-$ and $D_+$ of type $\sigma/2$, $b_\theta/2$, or $b_{o/i}/2$ such that $T$ is both the union of core transverse arcs of $D_-$ intersecting $D_+$ and one of $D_+$ intersecting $D_-$, and that the core periodic orbits intersecting $D_-$ coincide the periodic orbits intersecting $D_+$. 
Then the union of the core transverse arcs in $D_-$ and $D_+$ and core periodic orbits intersecting $D_-$ is called the {\bf core component} of $D_- \sqcup_T D_+$. 
Choosing one core component for any prong connection, the union of such core components is called a {\bf core} of $\F$. 
We have the following completeness. 

\begin{proposition}\label{prop:com_inv_02}
For any elements $\F$ and $\G'$ of the subspace $\mathcal{L}^{0}(M)$, if there are cores of elements $\F$ and $\G'$ which are isomorphic as a surface graph, then the prong connection diagrams of $\F$ and $\G$ are isomorphic as a surface graph. 
\end{proposition}

\begin{proof}
Let $D(\F)$ be the prong connection diagram of $\F$ and $H_\F$ the height function of $\F$. 
Fix a core component $C(D_- \sqcup_T D_+)$ of $\F$. 
Then $D_-$ (resp. $D_+$) is a domain of type $\sigma/2$, $b_\theta/2$, or $b_{o/i}/2$, and so they contain exactly one prong $q_-$ (resp. $q_+$), which is either a $1$-prong or a $3$-prong because of $\F \in \mathcal{L}^{0}(M)$. 
Since the union $D_- \sqcup_T D_+$ is a closed invariant subset, the intersection $D(\F) \cap (D_- \sqcup_T D_+)$ is a finite graph which intersects $T$ finitely many times. 
Because that the union $D_- \sqcup_T D_+$ is a closed invariant subset, by Lemma~\ref{lem:regularity_proper}, for any non-closed leaf $L$ in $D_- \sqcup_T D_+$, the leaf $L$ is a prong connection whose closure $\overline{L}$ is a closed curve between prongs with $\overline{L} - L \subseteq \{ q_- , q_+ \}$, and so any semi-prong separatrix intersecting $D_- \sqcup_T D_+$ is a prong separatrix among $\{ q_- , q_+ \}$. 
Since both prongs $q_-$ and $q_+$ have an odd number of separatrices, there is a prong separatrix connecting $q_-$ and $q_+$. 
Therefore, we obtain $H_\F(q_-) = H_\F(q_+)$. 
Then $C := D(\F) \cap (D_- \sqcup_T D_+)$ is the prong connection containing $q_-$ such that $D - C \subseteq \mathop{\mathrm{Per}}(\F)$. 
By the existence of a prong separatrix connecting $q_-$ and $q_+$, the prong connection $C$ is type $\sigma$ (resp. $b_o$ or $b_i$) if $\F\vert_{D_-}$ and $\F\vert_{D_+}$ are of type $\sigma/2$ (resp. $b_{o/i}$). 
Moreover, the prong connection $C$ is type $p_o$ or $p_i$ if $\F\vert_{D_-}$ and $\F\vert_{D_+}$ are of type $\sigma/2$ and $b_{o/i}$. 

\begin{claim}\label{claim:19}
The prong connection $C$ is type $b_\theta$ if $\F\vert_{D_-}$ and $\F\vert_{D_+}$ are of type $b_{\theta}/2$.
\end{claim}

\begin{proof}[Proof of Claim~\ref{claim:19}]
Suppose that $\F\vert_{D_-}$ and $\F\vert_{D_+}$ are of type $b_{\theta}/2$.
Assume that there are elements $l_- \in \F\vert_{D_-}$ with $q_- \in \overline{l_-} - l_-$ and $l_+ \in \F\vert_{D_+}$ with $q_+ \in \overline{l_+} - l_+$ such that $l_-$ and $l_+$ intersect a same connected component $T'$ of $T$ and that $l_- \cap l_+ = \emptyset$. 
Denote by $t'_- \in T'$ the point with $\{ t'_- \} = l_- \cap T'$ and by $t'_+ \in T'$ the point with $\{ t'_+ \} = l_+ \cap T'$. 
Let $[t'_-, t'_+] \subset T'$ be the sub-arc whose bounary is $\{ t'_-, t'_+ \}$. 
Since there is a small \nbd $I_-$ (resp. $I_+$) of $t'_-$ (resp. $t'_+$) in $[t'_-, t'_+] \subset T'$ consists of periodic points except $t'_-$ (resp. $t'_+$) in $T'$, the restriction $H_\F \vert_{I_-}$ (resp. $H_\F \vert_{I_+}$) is strictly monotonic. 
This implies the extrema $t'_0$ in an open interval $[t'_-, t'_+] - \{ t'_-, t'_+ \}$ of $H_\F$ with $H_\F(t'_0) \neq H_\F(t'_-) = H_\F(t'_+)$. 
Therefore, the point $t'_0 \in D_- \sqcup_T D_+$ is contained in a prong connection $C' \neq C$, which contradicts $C = D(\F) \cap (D_- \sqcup_T D_+)$. 
Thus $C$ is type $b_\theta$. 
\end{proof}
Therefore, any core component determines the unique prong connection, and so any core determines the unique prong connection diagram. 
In addition, if the cores of elements of $\mathcal{L}^{0}(M)$ are isomorphic as a surface graph, then their prong connection diagrams are isomorphic as a surface graph. 
\end{proof}

\begin{lemma}\label{prop:com_inv}
The following are equivalent for any $\F, \G \in \mathcal{L}^{0}(M)$:
\\
{\rm(1)} The prong connection diagrams of elements $\F$ and $\G$ are isomorphic as a surface graph.  
\\
{\rm(2)} The elements $\F$ and $\G$ are topologically equivalent. 
\end{lemma}

\begin{proof}
By definition of topologically equivalence, assertion (2) implies assertion (1). 
Conversely, by definition of prong connection diagram, each connected component of the complement of the prong connection diagram of a foliation in $\mathcal{L}^*_{i}(M)$ is either an annulus, a M{\"o}bius band, or a Klein bottle. 
On the other hand, because the leaf spaces of any foliations on such an invariant annulus (resp. M{\"o}bius band, Klein bottle) which is tangent to the boundary and consists of periodic leaves are topologicaly equivalent are intervals, such foliations on the annulus (resp. M{\"o}bius band, Klein bottle) are topologically equivalent respectively. 
Therefore, $\F$ is topologically equivalent to $\G$. 
\end{proof}

We have the following structural stability for the compact surface cases. 

\begin{proposition}\label{prop:stability_infty}
Suppose that $M$ is compact. 
For any $i \in \Z_{>1}$, the subspace $\mathcal{L}^{0}_{i}(M)$ consists of structurally stable ones in $\mathcal{L}^*_{i}(M)$ with respect to the relative topology of the topology $\mathcal{O}^{\infty}_{w,W,\mathcal{L}(M)}$. 
\end{proposition}

\begin{proof}
Choose periodic leaves $L_1 , \ldots, L_k$ and $V_{\F,1}, \ldots, V_{\F,l}$ connected components of $M - \bigsqcup_{j=1}^k L_j$ as in the proof of Lemma~\ref{lem:stability04+}. 
Fix any core $C(\F)$ of $\F$. 
From the existence of fibered charts, any core transverse arc $T$ of $\F$ can be extended to a closed transverse arc whose interior contains $T$. 
By Lemma~\ref{lem:trans}, the core transverse arcs of $\F$ is transverse to any $\G \in \mathcal{L}^*_{i}(M)$ near $\F$. 
Since the sums of the indices of singular points of $\F$ and $\G$ whose indices are positive are $i$, Proposition~\ref{prop:inv_index} implies that the numbers of $1$-prongs (resp. $3$-prongs) of $\F$ equals to the numbers of $1$-prongs (resp. $3$-prongs) of $\G$, and that the Hausdorff distance of $\mathop{\mathrm{Sing}}(\F)$ and $\mathop{\mathrm{Sing}}(\G)$ is small for any $\G \in \mathcal{L}^*_{i}(M)$ near $\F$. 

Fix any $\G \in \mathcal{L}^*_{i}(M)$ near $\F$. 
Therefore, from Lemma~\ref{lem:per_persistence_03}, the existence of the above extension of any core transverse arc of $\F$ implies that there is a core $C(\G)$ of $\G$ which is isotopic to $C(\F)$ for any $\G \in \mathcal{L}^*_{i}(M)$ near $\F$. 
In particular, the cores $C(\mathcal {F})$ and $C(\mathcal {G})$ are isomorphic as surface graphs. 
Lemma~\ref{prop:com_inv} implies that $\F$ and $\G$ are topologically equivalent, which implies the assertion. 
\end{proof}

The previous proposition does not hold for $p \neq \infty$ as follows. 

\begin{lemma}\label{lem:non_stability_p}
Suppose that $M$ is compact. 
For any $p \in \R_{\geq 1}$ and any $i \in \Z_{>1}$, the subspace $\mathcal{L}^{0}_{i}(M)$ does not consist of structurally stable ones in $\mathcal{L}^*_{i}(M)$ with respect to the relative topology of the topology $\mathcal{O}^{p}_{w,W,\mathcal{L}(M)}$. 
\end{lemma}

\begin{proof}
Replacing a \nbd of a heteroclinic prong separatrix in a prong connection whose type is $p_o$, we can deform the prong connection whose type is $p_o$ into a prong connection whose type is $p_i$ by arbitrarily small perturations whose supports are a small \nbd of the heteroclinic prong separtrix between the two prongs with respect to the relative topology of the topology $\mathcal{O}^p_{w,W,\mathcal{L}(M)}$. 
\end{proof}

Lemma~\ref{lem:stability04+} and Proposition~\ref{prop:stability_infty} imply Theorem~\ref{th:stability03-} and Theorem~\ref{th:stability03}. 
Note that one can demonstrate a similar statement without the hypothesis of the non-existence of fake prongs using similar arguments but the statement becomes more complicated.

\section{Coheight-one elements in the space of levelable foliations on surfaces}\label{sec:codim_one}

To characterize coheight-one elements in the space of levelable foliations on surfaces, we introduce the following instabilities, which are analogical concepts for flows (cf. \cite[Definitions~2.5 and 2.6]{sakajo2015transitions} and \cite[Definitions~1.10 and 1.11]{kibkalo2021topological}).

\subsection{Generic unstable elments}

Let $\F$ be a levelable foliation with finitely many pronged singular points on a surface $M$ and $D(\F)$ the prong connection diagram.  
We have the following instabilities.

\begin{figure}[t]
\begin{center}
\includegraphics[scale=1.]{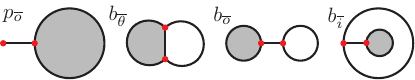}
\end{center} 
\caption{The list of $b$-unstable prong connections.}
\label{b_prong}
\end{figure}

\begin{definition}
A prong connection $C$ is {\bf $b$-unstable} if any prongs in $C$ are $1$-prongs and $3$-prongs, $n_{V_-}(C) -1 = b_1(C)$, and $n_{\partial}(C) = 1$, see Figure~\ref{b_prong}. 
\end{definition}

\begin{figure}[t]
\begin{center}
\includegraphics[scale=1.]{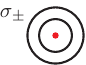}
\end{center} 
\caption{The list of $c$-unstable prong connection.}
\label{c_prong}
\end{figure}

\begin{definition}
A prong connection $C$ is {\bf $c$-unstable} if it consists of one center (i.e. $0$-prong), see Figure~\ref{c_prong}.
\end{definition}

\begin{definition}
A prong connection $C$ is {\bf $f$-unstable} if it contains exactly one fake prong, any prongs except the fake prong in $C$ are $1$-prongs and $3$-prongs, $n_{V_-}(C) = b_1(C)$, and $D(\F) \cap \partial M = \emptyset$. 
\end{definition}

\begin{figure}[t]
\begin{center}
\includegraphics[scale=1.25]{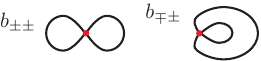}
\end{center} 
\caption{The list of $p$-unstable prong connections.}
\label{p_unstable_prong_connection}
\end{figure}

\begin{definition}
A prong connection $C$ is {\bf $p$-unstable} if $\mathrm{coheight}_{\mathrm{p,\F}}(C) = 2$, $n_{V_-}(C) + 1 = b_1(C)$, and $D(\F) \cap \partial M = \emptyset$, see Figure~\ref{p_unstable_prong_connection}.
\end{definition}

\begin{figure}[t]
\begin{center}
\includegraphics[scale=.8]{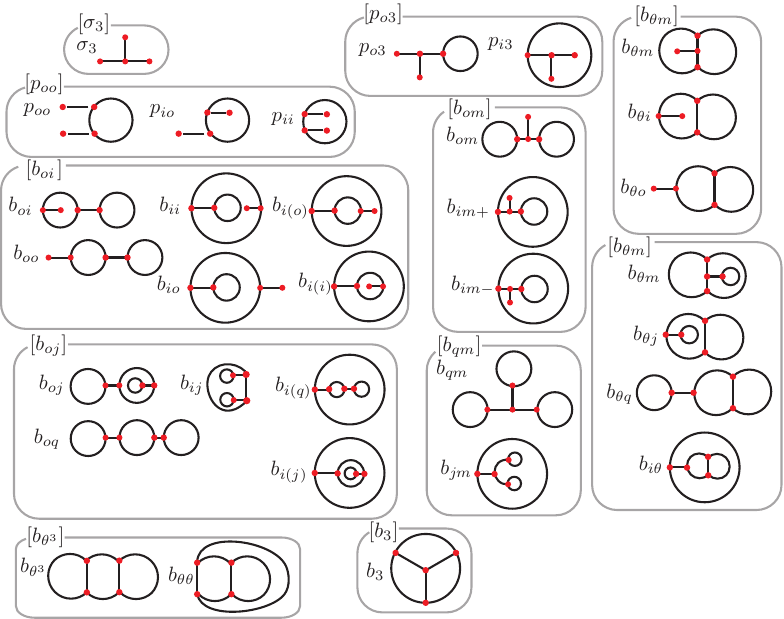}
\end{center} 
\caption{The list of $h$-unstable prong connections as topological graphs, which are equivalence classes of planar graphs.}
\label{h_unstable_prong_connection}
\end{figure}

\begin{definition}
A prong connection $C$ is {\bf $h$-unstable} if any prongs in $C$ are $1$-prongs and $3$-prongs, $n_{V_-}(C) -1 = b_1(C)$, and $D(\F) \cap \partial M = \emptyset$, see Figure~\ref{h_unstable_prong_connection}. 
\end{definition}

Notice that every prong connection $C$ is $b$-unstable (resp. $f$-unstable, $h$-unstable) if and only if $n_{\partial}(\F)  = \mathop{\mathrm{coheight}_{\F}}(C) = 1$ (resp. $n_{\mathrm{fake}}(\F) = \mathop{\mathrm{coheight}_{\F}}(C) = 1$, $\mathop{\mathrm{coheight}_{\mathrm{h},\F}}(C)  = \mathop{\mathrm{coheight}_{\F}}(C) = 1$). 


\begin{definition}
The prong connection diagram of a levelable foliation on a surface is $b$-unstable (resp. $c$-unstable, $f$-unstable, $h$-unstable) if it contains exactly one $b$-unstable (resp. $c$-unstable, $f$-unstable, $h$-unstable) prong connection and any prong connections except the coheight-one prong connection are coheight-zero. 
\end{definition}

\subsubsection{Characterization of $p$-unstable levelable foliations}

Let $n_{V_+}(C)$ be the number of prongs in $C$ whose indices are positive. 
Then $n_{V_+}(C) = n_{V}(C) - n_{V_-}(C)$. 
We have the following observation. 

\begin{lemma}\label{lem:ch_p-unstable}
Each $p$-unstable prong connection is isomorphic to one of Figure~\ref{p_unstable_prong_connection} as a topological graph. 
\end{lemma}

\begin{proof}
Let $C$ be a $p$-unstable prong connection. 
From $D(\F) \cap \partial M = \emptyset$, any prongs in $C$ are outside of the boundary $\partial M$. 
By $\mathrm{coheight}_{\mathrm{p,\F}}(C) = 2$, since the coheight of $k$-prong $y$ ($k \geq 3$) with respect to degeneracy of singular points is $\mathrm{coheight}_{\mathrm{p},\F}(y) = 2(k-3)$, the prong connection $C$ contains exactly one $4$-prong and any prongs except the $4$-prong are $1$-prongs and $3$-prongs. 
Let $l$ be the number of $3$-prongs in $C$. 
Then $n_{V_-}(C) = 1 + l$ and so $n_{V_+}(C) = n_{V}(C) - n_{V_-}(C) = n_{V}(C) - (1+l)$. 

\begin{claim}\label{claim:14}
$l = 0$. 
\end{claim}

\begin{proof}[Proof of Claim~\ref{claim:14}]
By $n_{E}(C) = \sum_{v \in \mathop{\mathrm{Sing}}(\F) \cap C} \deg (v)/2$, from Lemma~\ref{lem:betti}, we have the following equality: 
\[
\begin{split}
n_{V_-}(C) + 1 = b_1(C) = & 1 - n_{V}(C) + n_{E}(C) 
\\
= & 1 - n_{V}(C) + \sum_{v \in \mathop{\mathrm{Sing}}(\F) \cap C} \deg (v)/2 
\\
= & 1 + \sum_{v \in \mathop{\mathrm{Sing}}(\F) \cap C} \left( \dfrac{\deg (v)}{2} -1 \right) 
\\
= & 2 + (\vert \{ 3\text{-prong in } C \} \vert - \vert \{ 1\text{-prong in } C \} \vert) / 2 
\\
= & 2 + (\vert \{ 3\text{-prong in } C \} \vert - n_{V_+}(C)) / 2 
\\
= & 2 + (l - (n_{V_-}(C) - (1 + l))) / 2
\\
= & l + (5 - n_{V_-}(C))/2 
\end{split}
\]
Then $1 + l = n_{V_-}(C)  = 1 + 2l/3$ and so $l = 0$. 
\end{proof}

Thus $C$ consists of one $4$-prong and $1$-prongs. 
By the definition of $p$-instability, we have $b_1(C) = n_{V_-}(C) + 1 = 1 + 1 = 2$. 
From $2 = b_1(C) = 2 - \vert \{ 1\text{-prong in } C \} \vert / 2$, the prong connection  $C$ consists of the $4$-prong and two homoclinic separatrices. 
\end{proof}

\subsubsection{Characterization of $h$-unstable levelable foliations}

We have the following statement. 

\begin{lemma}\label{lem:ch_h-unstable}
Each $h$-unstable prong connection is isomorphic to one of Figure~\ref{h_unstable_prong_connection} as a topological graph. 
\end{lemma}

%

\begin{proof}
Let $C$ be a $p$-unstable prong connection. 
By definition of $h$-unstable, any prongs in $C$ are $1$-prongs and $3$-prongs, $n_{V_-}(C) -1 = b_1(C)$, and $D(\F) \cap \partial M = \emptyset$. 

Since the difference $n_{V_+}(C) = n_{V}(C) - n_{V_-}(C)$ is the number of $1$-prongs in $C$, by $n_{E}(C) = \sum_{v \in \mathop{\mathrm{Sing}}(\F) \cap C} \deg (v)/2$, from Lemma~\ref{lem:betti}, we have the following equality: 
\[
\begin{split}
n_{V_-}(C) -1 = b_1(C) = & 1 - n_{V}(C) + n_{E}(C) 
\\
= & 1 - n_{V}(C) + \sum_{v \in \mathop{\mathrm{Sing}}(\F) \cap C} \deg (v)/2 
\\
= & 1 + \sum_{v \in \mathop{\mathrm{Sing}}(\F) \cap C} \left( \dfrac{\deg (v)}{2} -1 \right) 
\\
= & 1 + (\vert \{ 3\text{-prong in } C \} \vert - \vert \{ 1\text{-prong in } C \} \vert) / 2 
\\
= & 1 + (n_{V_-}(C) - n_{V_+}(C)) / 2 
\\
= & 1 + (n_{V_-}(C) - (n_{V}(C) - n_{V_-}(C))) / 2 
\\
= & 1 + n_{V_-}(C) - n_{V}(C)/2
\end{split}
\]
Then $n_{V}(C) = 4$. 
If $n_{V_-}(C) = 4$ (i.e. any prongs in $C$ are $1$-prong), then $C$ is the disjoint union of two $\sigma$ and so is not connected, because of Figure~\ref{h_unstable_prf}, whcich contradicts to the connectivity of $C$ as in the left in Figure~\ref{h_unstable_prf}. 
\begin{figure}[t]
\begin{center}
\includegraphics[scale=.8]{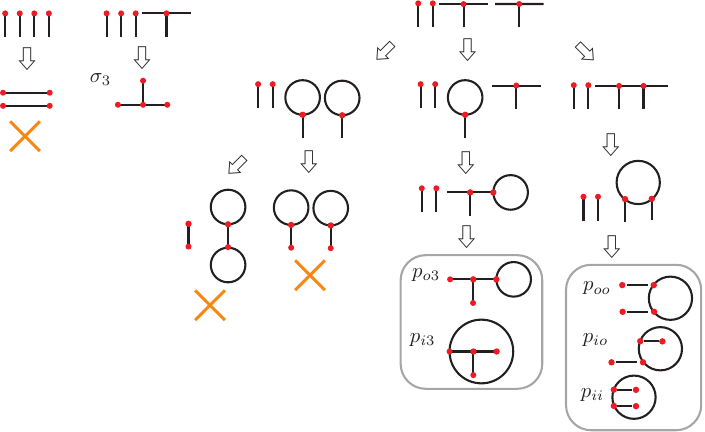}
\end{center} 
\caption{List of $h$-unstable prong connections for the case of zero, one, or two $3$-prongs, which is proved by the induction of the number of homoclinic separatrices.}
\label{h_unstable_prf}
\end{figure}
Thus $n_{V_-}(C) \in \{1,2,3\}$. 

\begin{figure}[t]
\begin{center}
\includegraphics[scale=.8]{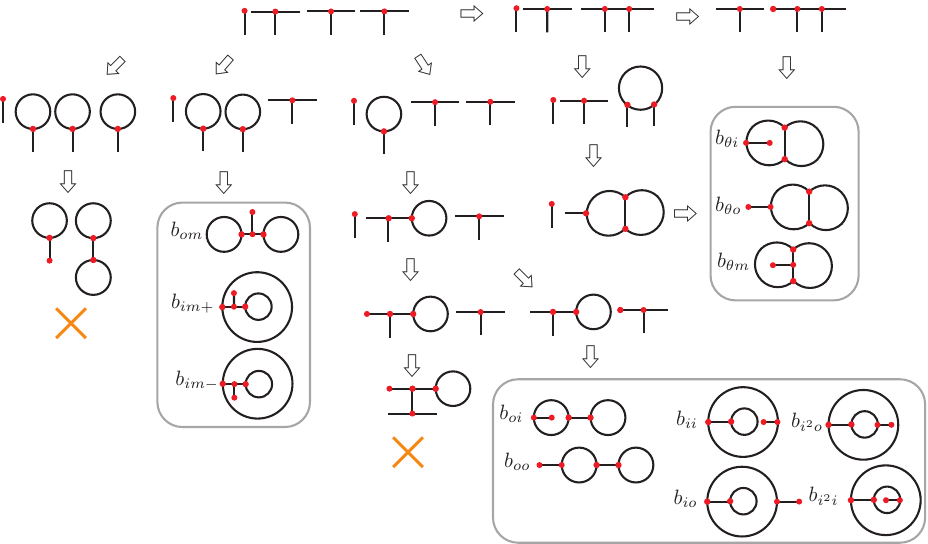}
\end{center} 
\caption{List of $h$-unstable prong connections for the case of three $3$-prongs, which is proved by the induction of the number of homoclinic separatrices.}
\label{h_unstable_prf_02}
\end{figure}

Suppose that $n_{V_-}(C) = 1$. 
Then $C$ consists of three $1$-prongs, a $3$-prong, and three separatrices between them. 
Therefore $C$ is isomorphic to $[\sigma_3]$ as in Figure~\ref{h_unstable_prong_connection} as a topological graph, because of Figure~\ref{h_unstable_prf}. 

Suppose that $n_{V_-}(C) = 2$. 
Assume that there are exactly two homoclinic separatrices in $C$. 
Then $C$ consists of two connected components as in Figure~\ref{h_unstable_prf}, which contradicts the connectivity of $C$. 
Thus there is at most one homoclinic separatrix in $C$ as in the middle and right in the case of the right in Figure~\ref{h_unstable_prf}. 
If there is exactly one homoclinic separatrix in $C$, then $C$ is $[p_{o3}]$ as in Figure~\ref{h_unstable_prong_connection} as a topological graph (i.e. $p_{o3}$ or $p_{i3}$). 
If there is no homoclinic separatrix in $C$, then $C$ is $[p_{oo}]$ as in Figure~\ref{h_unstable_prong_connection} as a topological graph (i.e. $p_{oo}$, $p_{io}$, or $p_{ii}$). 

Suppose that $n_{V_-}(C) =3$. 
Assume that there are exactly three homoclinic separatrices in $C$. 
Then $C$ consists of two connected components as in the left in Figure~\ref{h_unstable_prf_02}, because of Figure~\ref{h_unstable_prf}, which contradicts the connectivity of $C$. 
Thus, there are at most two homoclinic separatrices in $C$. 
If $C$ has two homoclinic separatrices in $C$, then the connectivity of $C$ implies that $C$ is $[b_{om}]$ as a topological graph as in Figure~\ref{h_unstable_prong_connection}, because of Figure~\ref{h_unstable_prf_02}. 
If $C$ has one homoclinic separatrix in $C$, then then the connectivity of $C$ implies that $C$ is $[b_{oi}]$ as a topological graph as in Figure~\ref{h_unstable_prong_connection}, because of Figure~\ref{h_unstable_prf_02}. 
Thus we may assume that there are no homoclinic separatrices in $C$. 
Then $C$ is $[b_{\theta m}]$ as in Figure~\ref{h_unstable_prong_connection} as a topological graph, because of Figure~\ref{h_unstable_prf_02}. 

\begin{figure}[t]
\begin{center}
\includegraphics[scale=.8]{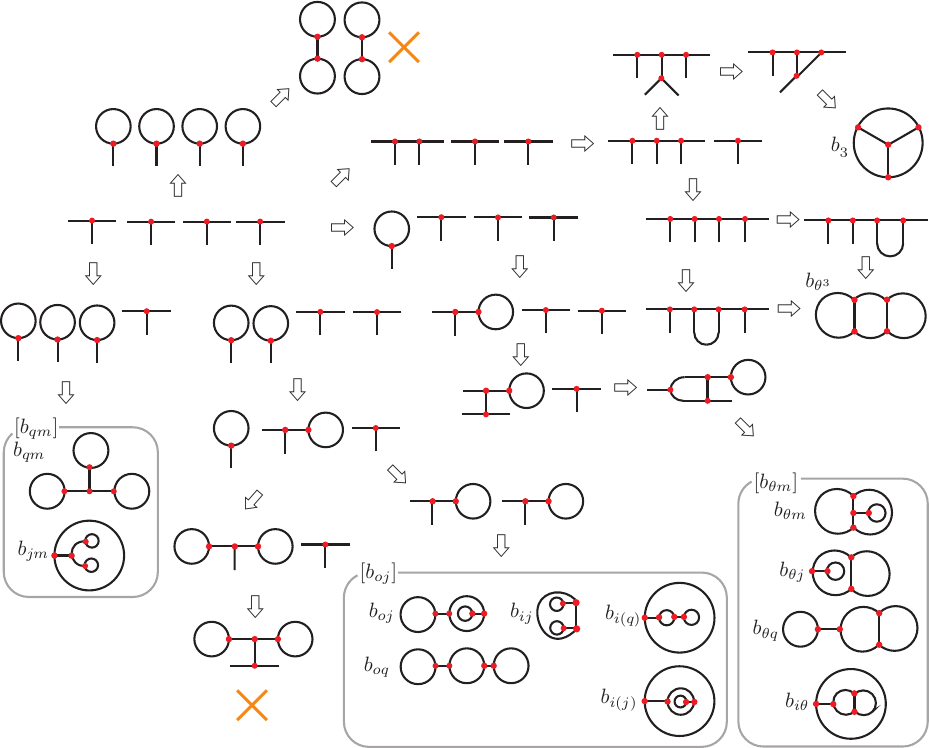}
\end{center} 
\caption{List of $h$-unstable prong connections for the case of four $3$-prongs, which is proved by the induction of the number of homoclinic separatrices.}
\label{h_unstable_prf_03}
\end{figure}

Suppose that $n_{V_-}(C) =4$. 
Assume that there are exactly four homoclinic separatrices in $C$. 
Then $C$ consists of two connected components as in the top of Figure~\ref{h_unstable_prf_03}, which contradicts the connectivity of $C$. 
Thus there are at most three homoclinic separatrices in $C$. 
If $C$ has exactly three homoclinic separatrices in $C$, then the connectivity of $C$ implies that $C$ is $[b_{qm}]$ as a topological graph as in Figure~\ref{h_unstable_prong_connection} and \ref{h_unstable_prf_03}. 
If $C$ has exactly two homoclinic separatrices in $C$, then the connectivity of $C$ implies that $C$ is $[b_{oj}]$ as a topological graph as in Figure~\ref{h_unstable_prong_connection} and \ref{h_unstable_prf_03}. 
If $C$ has exactly one homoclinic separatrix in $C$, then the connectivity of $C$ implies that $C$ is either $b_{\theta^3}$ or $[b_{\theta m}]$ as a topological graph as in Figure~\ref{h_unstable_prong_connection} and \ref{h_unstable_prf_03}. 
Thus, we may assume that there is exactly one homoclinic separatrix in $C$. 
The connectivity of $C$ implies that $C$ is $[b_{3}]$ as a topological graph as in Figure~\ref{h_unstable_prong_connection} and \ref{h_unstable_prf_03}. 
\end{proof}


We have the following characterization of coheight-one prong connections. 

\begin{figure}[t]
\begin{center}
\includegraphics[scale=1.2]{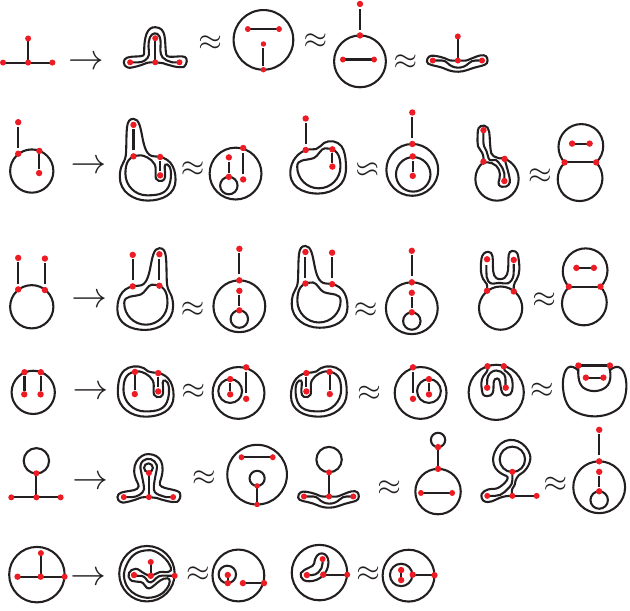}
\end{center} 
\caption{The list of transitions via $h$-unstable prong connections.}
\label{h_unstable_prong_transition}
\end{figure}

\begin{lemma}\label{lem:plus_zero}
Notice that every prong connection $C$ is either $c$-unstable or $p$-unstable if and only if $\mathrm{coheight}_{\mathrm{h},\F}(C) + \mathrm{coheight}_{\mathrm{p},\F}(C) = \mathop{\mathrm{coheight}_{\F}}(C) = 1$. 
In any case, we have $n_{\partial}(C) = n_{\mathrm{fake}}(C) = 0$. 
\end{lemma}

\begin{proof}
By the definition of $c$-unstable prong connection, we have the following equality: 
\[
\mathrm{coheight}_{\mathrm{h},\F}(C) + \mathrm{coheight}_{\mathrm{p},\F}(C) = 0+ 1 = 1 = \mathop{\mathrm{coheight}_{\F}}(C)
\] 
By the definition of $p$-unstable prong connection, we have the following equality: 
\[
\mathrm{coheight}_{\mathrm{h},\F}(C) + \mathrm{coheight}_{\mathrm{p},\F}(C) = 0+ 1 = 1 = \mathop{\mathrm{coheight}_{\F}}(C)
\]

Conversely, suppose that $\mathrm{coheight}_{\mathrm{h},\F}(C) + \mathrm{coheight}_{\mathrm{p},\F}(C) = \mathop{\mathrm{coheight}_{\F}}(C) = 1$. 
By definition of coheight, we have $n_{\partial}(C) = n_{\mathrm{fake}}(C) = 0$. 
If $C$ is $c$-unstable, then the assertion holds. 
Thus we may assume that $C$ is not $c$-unstable and so has centers. 
Then $n_{V_-}(C) \geq 1$. 
From Lemma~\ref{lem:1_3_prongs}, the prong connection $C$ has a $4$-prong and there is a coheight-zero prong connection $C'$ with $n_V(C') - n_V(C) = n_{V_-}(C') - n_{V_-}(C) = 1$. 
Then $n_{V_-}(C') \geq 2$.
By Proposition~\ref{lem:ch_codim_zero_01}, the prong connection is of type either $b_o$, $b_i$, or $b_\theta$ and so $2 = n_{V_-}(C') = n_{V_-}(C')$.  
Therefore, we have $1 = n_{V_-}(C) = n_{V_-}(C)$. 
Since $C$ has a $4$-prong, the prong connection $C$ consists of one $4$-prong and two self-connected prong separatrices and so is $p$-unstable.
\end{proof}

\begin{proposition}\label{lem:ch_codim_one_00}
Each connected component $C$ of $D(\F)$ is coheight-one if and only if it satisfies exactly one ofthe following conditions:
\\
{\rm(1)} The prong connection $C$ is $b$-unstable. 
\\
{\rm(2)} The prong connection $C$ is $c$-unstable. 
\\
{\rm(3)} The prong connection $C$ is $f$-unstable. 
\\
{\rm(4)} The prong connection $C$ is $h$-unstable. 
\end{proposition}

\begin{proof}
Let $ C$ be a prong connection of $\F$. 
Notice that every prong connection $C$ is $b$-unstable (resp. $f$-unstable, $h$-unstable) if and only if $n_{\partial}(\F)  = \mathop{\mathrm{coheight}_{\F}}(C) = 1$ (resp. $n_{\mathrm{fake}}(\F) = \mathop{\mathrm{coheight}_{\F}}(C) = 1$, $\mathop{\mathrm{coheight}_{\mathrm{h},\F}}(C)  = \mathop{\mathrm{coheight}_{\F}}(C) = 1$). 
By Lemma~\ref{lem:plus_zero}, every prong connection $C$ is either $c$-unstable or $p$-unstable if and only if $\mathrm{coheight}_{\mathrm{h},\F}(C) + \mathrm{coheight}_{\mathrm{p},\F}(C) = \mathop{\mathrm{coheight}_{\F}}(C) = 1$. 
Thus if $C$ is either $b$-unstable, $c$-unstable, $f$-unstable, and $h$-unstable, then $\mathop{\mathrm{coheight}_{\F}}(C) = 1$. 

Conversely, suppose that $\mathop{\mathrm{coheight}_{\F}}(C) = 1$. 
Then $1 = \mathrm{coheight}_{\F}(C) =  n_{\partial}(C) + n_{\mathrm{fake}}(C) + \mathop{\mathrm{coheight}_{\mathrm{h},\F}}(C) +  \mathrm{coheight}_{\mathrm{p,\F}}(C)$. 
Since the right-hand side of the previous equality consists of non-negative integers $n_{\partial}(C)$, $n_{\mathrm{fake}}(C)$, and $\mathop{\mathrm{coheight}_{\mathrm{h},\F}}(C) +  \mathrm{coheight}_{\mathrm{p,\F}}(C)$, the prong connection $C$ is either $b$-unstable, $c$-unstable, $f$-unstable, and $h$-unstable. 
\end{proof}

The previous proposition implies Theorem~\ref{main:char}. 
Moreover, Proposition~\ref{lem:ch_codim_one_00} and Theorem~\ref{cor:tree_rep} imply the following statement. 

\begin{theorem}\label{th:tree_rep_codim_one}
Let $S$ be a compact surface contained in a sphere $\mathbb{S}^2$ and $\F$ a coheight-one levelable foliation. 
Then the Reeb graph of the height function of $\F$ on $S$ is a finite directed tree whose vertices of degree one are either a prong connection of type $\sigma$, a prong connection of type $\sigma_3$, or a periodic leaf on the boundary $\partial \mathbb{S}^2$, and whose vertrice of degree non-one are either $p_o$, $p_i$, $b_\theta$, $b_o$, $b_i$, or one of prong connections in Figure~\ref{fig:list_self_connected_prong}. 
\end{theorem}

\subsection{Open dense property and structural stability in the space of levelable foliations on  compact surfaces}

Denote by $\mathcal{L}^1_{i}(M)$ the subspace of levelable foliation in $\mathcal{L}^*_{i}(M) - \mathcal{L}^0_{i}(M)$ whose prong connection diagram is of coheight-one. 
Since any elements of $\mathcal{L}^*_{i}(M)$ contain no fake prongs, the prong connection diagrams of any elements of $\mathcal{L}^*_{i}(M)$ are either $b$-unstable, $c$-unstable, $p$-unstable, or $h$-unstable. 
We state generic intermediate graphs as follows.

\subsubsection{}

The same argument of the proof of Lemma~\ref{lem:char_local_prongs_02} implies the following statement.

\begin{figure}[t]
\begin{center}
\includegraphics[scale=.5]{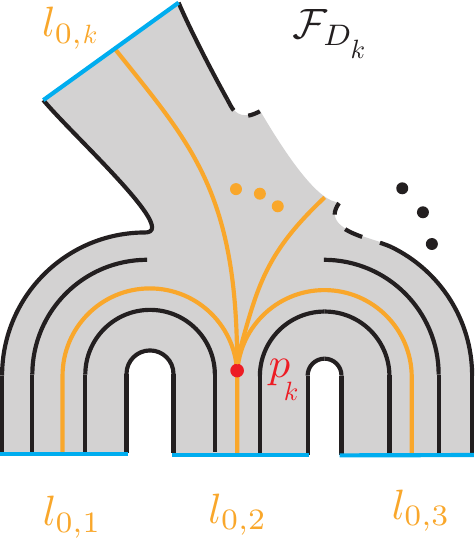}
\end{center} 
\caption{The partition on a disk.}
\label{fig:stdd_loc_str_k}
\end{figure}

\begin{lemma}\label{lem:char_local_prongs_k}
Let $\F$ be a levelable foliation with finitely many pronged singular points on a surface $M$ and $D \subset M$ a disk which contains no prong except one $k$-prong $q_k$ for any $k \in \Z_{\geq 3}$ and whose boundary consists of pairwise disjoint three closed leaf arc $C_1, \ldots ,C_k$ and pairwise disjoint three closed transverse arcs $T_1, \ldots ,T_k$ such that $(C_i \sqcup C_{i+1}) \cap T_i = \partial T_i$ and  $(T_i \sqcup T_{i+1}) \cap T_{i+1} = \partial T_i$ for any $i \in \{1, \ldots ,k\}$, where $C_{k+1} := C_1$ and $T_{k+1} := T_1$. 
Then the restriction $\F\vert_{D}$ is as in Figure~\ref{fig:stdd_loc_str_k}
\end{lemma}

\begin{theorem}\label{th:stability04+}
The following statements hold for any compact surface and any positive integer $i \geq 2$: 
\\
{\rm(1)} The subspace $\mathcal{L}^1_{i}(M)$ is an open dense subset of $\mathcal{L}^*_{i} - \mathcal{L}^0_{i}(M)$. 
\\
{\rm(2)} Any levelable foliation in $\mathcal{L}^1_{i}(M)$ are structurally stable in $\mathcal{L}^*_{i} - \mathcal{L}^0_{i}(M)$. 
\end{theorem}

\begin{proof}
The open denseness follows from Lemma~\ref{lem:codim_max}.
By a similar argument to the proof of Lemma~\ref{lem:char_local_prongs_03}, we can construct the core components of any prong connections such that the ``generalized core'' of line fields in $\mathcal{L}^1_{i}(M)$ determines the prong connection diagram. 
From the same argument of the proof of Proposition~\ref{prop:stability_infty}, using periodic leaves $L_1 , \ldots, L_k$ and $V_{\F,1}, \ldots, V_{\F,l}$ connected components of $M - \bigsqcup_{j=1}^k L_j$ as in the proof of Lemma~\ref{lem:stability04+}, we can obtain the structural stability. 
\end{proof}

The previous theorem implies Theorem~\ref{th:stability04}. 

\section{COT representations for generic levelable foliations on a sphere and a disk}\label{sec:cot_zero}

We introduce all local leaf structures contained in the prong connection diagram, and provide a unique symbol, called a {\bf COT symbol}, for each one, with these symbols being analogous to concepts from Hamiltonian flow.
The COT symbols express local leaf structures and the composition of their surrounding leaf structures.
\begin{figure}[t]
\begin{center}
\includegraphics[scale=1.]{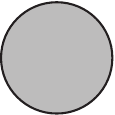}
\end{center} 
\caption{Boundary component.}
\label{fig:boudnary}
\end{figure}

Recall that, by Theorem~\ref{th:tree_rep_codim_zero}, for any coheight-zero levelable foliation $\F$ on a compact surface contained in a sphere $\mathbb{S}^2$, the Reeb graph of the height function of $\F$ is a finite directed tree whose vertices of degree one are either a prong connection of type $\sigma$ or a periodic leaf on the boundary $\partial \mathbb{S}^2$, and whose vertrice of degree non-one are either $p_o$, $p_i$, $b_\theta$, $b_o$, or $b_i$. 
Therefore, we can fix a root of the Reeb graph as a vertex corresponding to either a prong connection of type $\sigma$ or a periodic leaf on the boundary $\partial \mathbb{S}^2$. 
Then the line field equipped with the root is called a line field with a \textbf{root component}. 

First, we introduce coheight-zero local topological structures and their corresponding letters, named COT symbols, as follows.

\subsection{Candidates of roots, which are vertices of degree one}

To encode finite directed trees, we choose the roots of the finite directed trees. 
Therefore, we introduce the following two candidates for the roots of the trees. 

\bigskip
\noindent
\textbf{($\sigma_{\emptyset}$ structure)} The simplest one-dimensional local orbit structure is a closed interval $I$ at the point at infinity $\infty \in S_\sigma^2$ shown in Figure~\ref{fig:boudnary}, where $S_\sigma^2 := \R^2 \sqcup I$ is the sphere which is a compactification of the plain $\R^2$. 
We assign the COT symbol $\sigma_{\emptyset}(\B_{b})$, where $\B_{b}$ denotes the inner structure chosen as in Table~\ref{tbl:COT_Structures}.

\bigskip
\noindent
\textbf{($\beta_{\emptyset}$ structure)} The simplest two-dimensional local orbit structure is a closed disk $D$ at infinity $\infty \in S_\beta^2$ shown in Figure~\ref{self_connected_prong}, where $S_\beta^2 := \R^2 \sqcup D$ is the sphere. 
We assign the COT symbol $\beta_{\emptyset}(\B_{b})$.

\subsection{Candidates of terminals}

We introduce the following two candidates for the terminals of the trees. 

\bigskip
\noindent
\textbf{($\sigma$ structure)} The simplest one-dimensional local orbit structure is a closed interval shown in Figure~\ref{self_connected_prong}. 
We assign the COT symbol $\sigma$.

\bigskip
\noindent
\textbf{($\beta$ structure)} The simplest two-dimensional local orbit structure is a closed disk $D$ shown in Figure~\ref{fig:boudnary}. 
We assign the COT symbol $\beta$.

\subsection{Candidates of intermediate vertices}

We introduce the following two candidates of intermediate vertices of the trees. 

\bigskip
\noindent
\textbf{($p_{o}$, $p_{i}$ structures)} 
They are prong connections which consist of one $1$-prong, one $3$-prong, one prong separatrix between them, and one homoclinic prong separatrix of the $3$-prong shown in Figure~\ref{self_connected_prong}. 
When there is a disk whose boundary consists of the $3$-prong and its homoclinic prong separatrix which contains (resp. does not contain) the $1$-prong, the COT symbol is given by $p_i(\B_{b})$ (resp. $p_o(\B_{b})$).

\bigskip
\noindent
\textbf{($b_\theta$ structure)} The one-dimensional local orbit structure is a prong connection which consists of two $3$-prongs and three heteroclinic prong separatrices between them shown in Figure~\ref{self_connected_prong}.
We assign the COT symbol $b_\theta\{ \B_{b}, \B_{b} \}$.
Note that the $b_\theta$ structure has a rotational symmetry, and so its first inner structure and the second inner structure are arranged not in a line but in a circle. 
In other words, the two inner structures are cyclically ordered. 
Therefore, we enclosed the two inner structures in braces $\{ \}$.

\bigskip
\noindent
\textbf{($b_o$, $b_i$ structures)} They are prong connections which consist of two $3$-prongs, one heteroclinic prong separatrix between them, and two homoclinic separatrices of them respectively shown in Figure~\ref{self_connected_prong}.
When the prong connection has (reps. does not have) a rotational symmetry,  the COT symbol is given by $b_o\{\B_{b}, \B_{b}\}$ (resp. $b_i(\B^1_{b},\B^2_b)$), where $\B^2_{b}$ is the innermost structures and $\B^1_{b}$ is the intermediate structure between the innermost disk and the outer disk contains the structure at $\infty$.

\begin{figure}[t]
\begin{center}
\includegraphics[scale=1.]{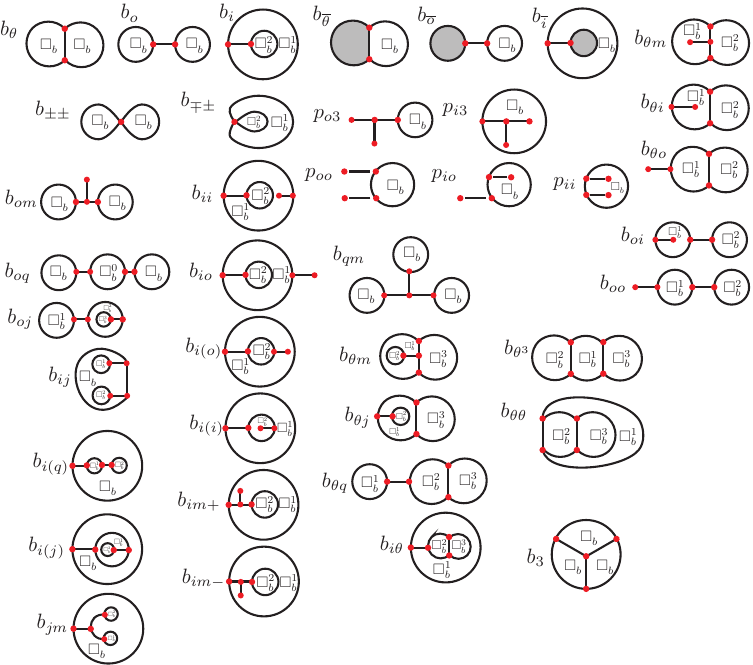}
\end{center} 
\caption{The list of COT symbols of at most coheight-one prong connections with nonempty inner structures.}
\label{fig:list_self_connected_prong}
\end{figure}

\subsection{COT representations of coheight-zero line fields on the sphere}

We have the nest structures for coheight-zero structures as shown in Table~\ref{tbl:All_structures} and Table~\ref{tbl:COT_Structures}. 

\begin{table}
\begin{center}
\begingroup
\begin{tabular}{|l|c|l|l|}\hline
\multicolumn{2}{|c|}{Root structure} & COT symbol & Figure\\ \hline
closed interval & $\sigma_{\emptyset}$ &$\sigma_{\emptyset}(\B_{b})$  & \ref{self_connected_prong} \\ \hline
boundary component & $\beta_{\emptyset}$ &$\beta_{\emptyset}(\B_{b})$  & \ref{fig:boudnary} \\  \hline \hline
\multicolumn{2}{|c|}{Terminal structure} & COT symbol & Figure \\ \hline
boundary component& $\beta$ & $\beta$ & \ref{fig:boudnary}  \\   \hline
prong connection & $\sigma$ & $\sigma$ &\ref{self_connected_prong} \\  \hline \hline
\multicolumn{2}{|c|}{Intermedidate structure} & COT symbol & Figure \\ \hline
& $p_o$ & $p_o(\B_b)$ &  \ref{fig:boudnary} \\ \cline{2-3}
prong connection & $p_i$ & $p_i(\B_b)$ & \\ \cline{2-3}
& $b_\theta$ & $b_\theta\{ \B_b, \B_b \}$ & \\ \cline{2-3}
& $b_o$ & $b_o\{ \B_b, \B_b \}$ & \\ \cline{2-3}
& $b_i$ & $b_i(\B_b,\B_b)$ & \\ \hline
\end{tabular}
\endgroup
\caption{All leaf structures generated by the levelable line field and their COT symbols of coheight-zero elements.}
\label{tbl:All_structures}
\end{center}
\end{table}

\begin{table}
\begin{center}
\begingroup
\begin{tabular}{|l|c|c|}\hline
Class & Boxes  & Orbit structures  \\ \hline
class-$b$ & $\B_{b}$ & $\{ \sigma, \beta, p_o, p_i, b_\theta, b_o, b_i \}$  \\ \hline
\end{tabular}
\endgroup
\caption{Classes of leaf structures used in COT symbols of coheight-zero elements.}
\label{tbl:COT_Structures}
\end{center}
\end{table}

As above, for any line field $\F$ of $\mathcal{L}^0(\mathbb{S}^2)$ with a root component, we  can assign the COT representation, which is a sequence of symbols, of $\F$ as follows: 
\\
(1) If the root component is of type $\sigma$ (resp. a periodic leaf on the boundary), then the intial structure of the sequence of symbols is $\sigma_{\emptyset}(\B_{b})$ (resp. $\beta_{\emptyset}(\B_{b})$). 
\\
(2) Since the Reeb graph has a root given by the vertex of degree one which corresponds to the root component, tracing adjacent vertices from the root uniquely determines the structure of the associated prong connections or boundary components. 
Therefore, for every vertex, we recursively substitute either $\sigma$, $\beta$, $p_o(\B_{b})$, $p_i(\B_{b})$, $b_\theta(\B_{b})$, $b_o(\B_{b})$, or $b_i(\B_{b})$, into $\B_b$ according to the corresponding prong connection. 
The substitution proceeds in order from the vertices closest to the root, where possible.

The resulting sequence of symbols is called the partially cyclically ordered rooted tree ({\bf COT}) {\bf representation} of the line field $\F$ of $\mathcal{L}^0(\mathbb{S}^2)$ with the root component. 

By construction of the COT representation, Proposition~\ref{lem:ch_codim_zero_01} implies the following completeness. 

\begin{theorem}\label{th:cot_01}
For any $p \in \R_{\geq 1} \sqcup \{ \infty \}$, the COT representation 
for line fields of $\mathcal{L}^0(\mathbb{S}^2)$ with root components is a finite complete invariant. 
\end{theorem}

\subsection{Examples of COT representations of levelable line fields}

\begin{figure}[t]
\begin{center}
\includegraphics[scale=0.95]{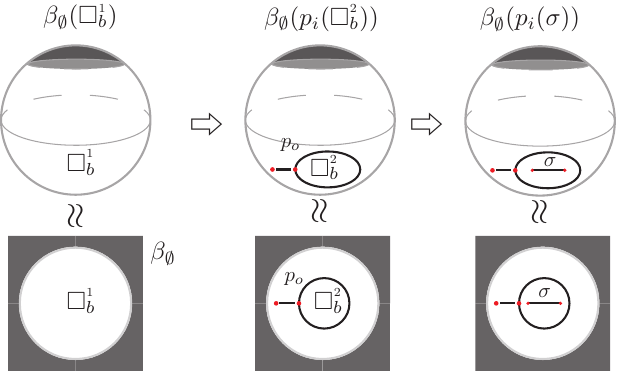}
\end{center} 
\caption{A construction of a line field whose COT representation is $\beta_{\emptyset}(p_i(\sigma))$.}
\label{fig:ex_cot_line_01}
\end{figure}

We have the COT representation $\beta_{\emptyset}(p_i(\sigma))$ of the example of a levelable line field on a sphere as shown on the right in Figure~\ref{fig:ex_cot_line_01}. 
In fact, the line field can be obtained by inserting local structures as follows: 
First, we consider the closed disk $D$ at infinity $\infty \in S_\beta^2$ shown in Figure~\ref{fig:boudnary} whose COT symbol is $\beta_{\emptyset}(\B^1_{b})$, where $S_\beta^2 := \R^2 \sqcup D$ is the sphere. 
Second, we insert a prong connection shown in Figure~\ref{self_connected_prong} whose COT symbol is $p_i(\B^2_{b})$ into the open disk corresponding to $\B^1$. 
Then the resulting COT representation is $\beta_{\emptyset}(p_i(\B_{b}))$. 
Finally, we insert a prong connection shown in Figure~\ref{self_connected_prong} whose COT symbol is $\sigma$ into the open disk corresponding to $\B^2$. 
Then the resulting COT representation is $\beta_{\emptyset}(p_i(\sigma))$. 
Notice that the line field on the sphere $S_\beta^2$ as shown on the right in Figure~\ref{fig:ex_cot_line_01} also can be considered a line field on a closed disk $S_\beta^2 - \mathop{\mathrm{int}}D$. 

Moreover, we construct a line field whose COT representation is $\sigma_{\emptyset}(b_i(\sigma,\sigma))$ as shown on the right in Figure~\ref{fig:ex_cot_line_01} as follows: 
\begin{figure}[t]
\begin{center}
\includegraphics[scale=0.95]{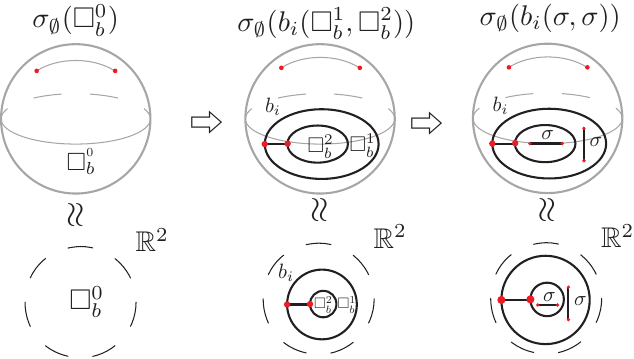}
\end{center} 
\caption{A construction of a line field whose COT representation is $\sigma_{\emptyset}(b_i(\sigma,\sigma))$.}
\label{fig:ex_cot_line_02}
\end{figure}
First, we consider a closed interval $I$ at the point at infinity $\infty \in S_\sigma^2$ shown in Figure~\ref{fig:boudnary} whose COT symbol is $\sigma_{\emptyset}(\B^0_{b})$, where $S_\sigma^2 := \R^2 \sqcup I$ is the sphere which is a compactification of the plain $\R^2$. 
Second, we insert a prong connection shown in Figure~\ref{self_connected_prong} whose COT symbol is $b_i(\B^1_{b},\B^2_{b})$ into the open disk corresponding to $\B^1$. 
Then the resulting COT representation is $\sigma_{\emptyset}(b_i(\B^1_{b},\B^2_{b}))$. 
Finally, we insert two prong connections whose COT symbols are $\sigma$ into the open disks corresponding to $\B^1$ and $\B^2$ respectively. 
Then the resulting COT representation is $\sigma_{\emptyset}(b_i(\sigma,\sigma))$. 
Notice that the line field on the sphere $S_\sigma^2$ as shown on the right in Figure~\ref{fig:ex_cot_line_01} also can be considered a line field on the plane $\R^2$ with the cyclic boundary condition. 

\section{Intermediate foliations between COT representations for generic levelable foliations on a sphere and a disk}\label{sec:cot_one}

We introduce coheight-one local topological structures and their COT symbols as follows.

\begin{table}
\begin{center}
\begingroup
\begin{tabular}{|l|c|c|}\hline
Class & Boxes  & Orbit structures  \\ \hline
class-$\sigma_{\emptyset}$ & $\B_{\sigma_{\emptyset}}$ & $\{ \beta_{\emptyset}, \sigma_{\emptyset} \}$   \\ \hline
class-$b$ & $\B_{b}$ & $\{ \sigma, \beta, p_o, p_i, b_\theta, b_o, b_i,\sigma_\pm,p_{\overline{o}}p_{\overline{o}},b_{\overline{\theta}},b_{\overline{o}},b_{\overline{i}},\sigma_\pm,b_{\pm\pm},b_{\pm\mp},\sigma_{3},$ \\
& & $p_{o3},p_{i3},p_{oo},p_{io},p_{ii}, b_{\theta m}, b_{\theta i}, b_{\theta o},b_{om}, b_{im+}, b_{im-},b_{oi},$ 
\\
& &
$b_{oo}, b_{io}, b_{ii},b_{i(o)}, b_{i(i)},b_{oj}, b_{oq},b_{ij}, b_{i(q)}, b_{i(j)},b_{qm}, b_{jm},$ 
\\
& &
$b_{\theta m}, b_{\theta j}, b_{\theta q}, b_{i\theta},b_{\theta^3}, b_{\theta \theta},b_{3}\}$ 
\\\hline
\end{tabular}
\endgroup
\label{tbl:COT_Structures_codim_1}
\caption{Classes of leaf structures used in COT symbols of at most coheight-one elements.}
\end{center}
\end{table}

\begin{table}
\begin{center}
\begingroup
\begin{tabular}{|l|c|l|l|}\hline
\multicolumn{2}{|c|}{Prong connection} & COT symbol & Figure \\ \hline
$b$-unstable& $p_{\overline{o}}$ & $p_{\overline{o}}$ &\ref{b_prong} \\ \cline{2-3}
& $b_{\overline{\theta}}$ & $b_{\overline{\theta}}(\B_b)$ &  \\ \cline{2-3}
 & $b_{\overline{o}}$ & $b_{\overline{o}}(\B_b)$ &  \\ \cline{2-3}
& $b_{\overline{i}}$ & $b_{\overline{i}}(\B_b)$ & \\ \hline
$c$-unstable & $\sigma_\pm$ & $\sigma_\pm$ &\ref{c_prong} \\\hline
$p$-unstable& $b_{\pm\pm}$ & $b_{\pm\pm}\{\B_b, \B_b \}$ &\ref{p_unstable_prong_connection} \\ \cline{2-3}
& $b_{\pm\pm}$ & $b_{\pm\pm}(\B^1_b, \B^2_b)$ & \\ \hline
$h$-unstable & $\sigma_3$ & $\sigma_3$ &\ref{h_unstable_prong_connection} \\ \cline{2-3}
& $p_{o3}$ & $p_{o3}(\B_b)$ & \\ \cline{2-3}
& $p_{i3}$ & $p_{i3}(\B_b)$ & \\ \cline{2-3}
& $p_{oo}$ & $p_{oo}(\B_b)$ & \\ \cline{2-3}
& $p_{\theta i}$ & $p_{\theta i}(\B_b)$ & \\ \cline{2-3}
& $p_{\theta m}$ & $p_{\theta m}(\B_b)$ & \\ \cline{2-3}
& $b_{\theta o}$ & $b_{\theta o}(\B^1_b,\B^2_b)$ & \\ \cline{2-3}
& $b_{\theta i}$ & $b_{\theta i}(\B^1_b,\B^2_b)$ & \\ \cline{2-3}
& $b_{\theta m}$ & $b_{\theta m}(\B^1_b,\B^2_b)$ & \\ \cline{2-3}
& $b_{om}$ & $b_{om}\{\B_b,\B_b\}$ & \\ \cline{2-3}
& $b_{im+}$ & $b_{im+}(\B^1_b,\B^2_b)$ & \\ \cline{2-3}
& $b_{im-}$ & $b_{im-}(\B^1_b,\B^2_b)$ & \\ \cline{2-3}
& $b_{oo}$ & $b_{oo}(\B^1_b,\B^2_b)$ & \\ \cline{2-3}
& $b_{oi}$ & $b_{oi}(\B^1_b,\B^2_b)$ & \\ \cline{2-3}
& $b_{io}$ & $b_{io}(\B^1_b,\B^2_b)$ & \\ \cline{2-3}
& $b_{ii}$ & $b_{ii}(\B^1_b,\B^2_b)$ & \\ \cline{2-3}
& $b_{i(o)}$ & $b_{i(o)}(\B^1_b,\B^2_b)$ & \\ \cline{2-3}
& $b_{i(i)}$ & $b_{i(i)}(\B^1_b,\B^2_b)$ & \\ \cline{2-3}
& $b_{oq}$ & $b_{oq}(\B^0_b, \{ \B_b, \B_b \})$ & \\ \cline{2-3}
& $b_{oj}$ & $b_{oj}(\B^1_b, \B^2_b, \B^3_b)$ & \\ \cline{2-3}
& $b_{ij}$ & $b_{ij}(\B_b, \{ \B^1_b, \B^2_b \})$ & \\ \cline{2-3}
& $b_{i(q)}$ & $b_{i(q)}(\B_b, \B^1_b, \B^2_b)$ & \\ \cline{2-3}
& $b_{i(j)}$ & $b_{i(j)}(\B^1_b, \{ \B^2_b, \B^3_b \})$ & \\ \cline{2-3}
& $b_{qm}$ & $b_{qm} \{ \B_b, \B_b, \B_b\}$ & \\ \cline{2-3}
& $b_{jm}$ & $b_{jm}(\B_b, \B^1_b, \B^2_b)$ & \\ \cline{2-3}
& $b_{\theta m}$ & $b_{\theta m}(\B^1_b, \B^2_b, \B^3_b)$ & \\ \cline{2-3}
& $b_{\theta j}$ & $b_{\theta j}(\B^1_b, \B^2_b, \B^3_b)$ & \\ \cline{2-3}
& $b_{\theta q}$ & $b_{\theta q}(\B^1_b, \B^2_b, \B^3_b)$ & \\ \cline{2-3}
& $b_{i\theta}$ & $b_{i\theta}(\B^1_b, \B^2_b, \B^3_b)$ & \\ & $b_{\theta^3}$ & $b_{\theta^3}(\B^1_b, \B^2_b, \B^3_b)$ & \\ \cline{2-3}
& $b_{\theta \theta}$ & $b_{\theta \theta}(\B^1_b, \B^2_b, \B^3_b)$ & \\ \cline{2-3}
& $b_{3}$ & $b_{3}\{ \B_b, \B_b, \B_b\}(\B^1_b, \B^2_b, \B^3_b)$ & \\  \hline
\end{tabular}
\endgroup
\label{tbl:All_structures_codim_1}
\caption{All leaf structures generated by the levelable line field and their COT symbols of coheight-one elements.}
\end{center}
\end{table}

\subsection{The $c$-unstable prong conneciton}

We introduce a coheight-one local topological structure corresponding to the $c$-unstable prong connection. 

\bigskip
\noindent
\textbf{($\sigma_\pm$ structure)} The simplest one-dimensional local orbit structure is a prong connection which consists of two $1$-prongs and a prong separatrix between them shown in Figure~\ref{self_connected_prong}.
We assign the COT symbol $\sigma_\pm$.

\subsection{The $b$-unstable prong connecitons}

We introduce coheight-one local topological structures corresponding to $b$-unstable prong connections. 

\bigskip
\noindent
\textbf{($p_{\overline{o}}$ structure)} 
It is a prong connection which consists of one $1$-prong, one $3$-prong, one prong separatrix between them, and one homoclinic prong separatrix on the boundary of the $3$-prong shown in Figure~\ref{b_prong}. 
We assign the COT symbol $p_{\overline{o}}$.

\bigskip
\noindent
\textbf{($b_{\overline{\theta}}$ structure)} The one-dimensional local orbit structure is a prong connection which consists of two $3$-prongs and three heteroclinic prong separatrices between them exactly two of which are contained in the boundary shown in Figure~\ref{b_prong}.
We assign the COT symbol $b_{\overline{\theta}}(\B_b)$.

\bigskip
\noindent
\textbf{($b_{\overline{o}}$, $b_{\overline{i}}$ structures)} They are prong connections which consists of two $3$-prongs, one heteroclinic prong separatrix between them, and two homoclinic separatrices of them exactly one of which is contained in the boundary shown in Figure~\ref{b_prong}.
When there is (resp. is not) a disk whose boundary consists of the $3$-prong and its homoclinic prong separatrix, the COT symbol is given by $b_{\overline{o}}(\B_b)$ and  $b_{\overline{i}}(\B_b)$. 

\subsection{The $c$-unstable prong connecitons}

We introduce coheight-one local topological structures corresponding to $c$-unstable prong connections. 

\bigskip
\noindent
\textbf{($\sigma_\pm$ structure)} The simplest zero-dimensional local orbit structure is a center shown in Figure~\ref{c_prong}.
We assign the COT symbol $\sigma_\pm$.

\subsection{The $p$-unstable prong connecitons}

We introduce coheight-one local topological structures corresponding to $p$-unstable prong connections. 

\bigskip
\noindent
\textbf{($b_{\pm\pm}$, $b_{\pm\mp}$ structures)} 
They are one-dimensional orbit structures consisting of a saddle and two homoclinic prong separatrices as shown in Figure~\ref{p_unstable_prong_connection}.
When the boundary component of the outer domain is (resp. is not) the prong connection, the COT symbol is given by $b_{\pm\pm}\{\B_b, \B_b \}$ (resp. $b_{\pm\pm}(\B^1_b, \B^2_b)$), where the boundary of the disk containing $\B^1_{b}$ (resp. $\B^2_{b}$) is (resp. is not) the prong connection. 

\subsubsection{The $h$-unstable prong connecitons}

We introduce coheight-one local topological structures corresponding to $h$-unstable prong connections. 

\bigskip
\noindent
\textbf{($\sigma_{3}$ structure)} The one-dimensional local orbit structure is a prong connection which consists of one $3$-prong, three $1$-prongs, and three heteroclinic prong separatrices between the $3$-prong and three $1$-prongs respectively shown in Figure~\ref{h_unstable_prong_connection}.
We assign the COT symbol $\sigma_3$.

\bigskip
\noindent
\textbf{($p_{o3}$, $p_{i3}$ structures)} The one-dimensional local orbit structure is a prong connection which consists of two $3$-prongs, two $1$-prongs, three heteroclinic prong separatrices between them, and one homoclinic prong separatrix shown in Figure~\ref{h_unstable_prong_connection}.
When there is a disk whose boundary is contained in the prong connection and which contains (does not contain) three prongs in the prong connection, the COT symbol is given by $p_{o3}(\B_b)$ (resp. $p_{i3}(\B_b)$).

\bigskip
\noindent
\textbf{($p_{oo}$, $p_{io}$, $p_{ii}$ structures)} The one-dimensional local orbit structure is a prong connection which consists of two $3$-prongs, two $1$-prongs, and four heteroclinic prong separatrices between them shown in Figure~\ref{h_unstable_prong_connection}.
When the boundary component of the outer domain contains exactly $four$ (resp. $three$, $two$) prongs, he COT symbol is given by $p_{oo}(\B_b)$ (resp. $p_{\theta i}(\B_b)$, $p_{\theta m}(\B_b)$), where $\B^2_b$ is contained in a disk whose boundary consists of two prongs and two prong separatrices. 

\bigskip
\noindent
\textbf{($b_{\theta m}$, $b_{\theta i}$, $b_{\theta o}$ structures)} The one-dimensional local orbit structure is a prong connection which consists of three $3$-prongs, one $1$-prong, and five heteroclinic prong separatrices between them shown in Figure~\ref{h_unstable_prong_connection}.
When the boundary component of the outer domain contains exactly four (resp. three, two) prongs, the COT symbol is given by $b_{\theta o}(\B^1_b,\B^2_b)$ (resp. $b_{\theta i}(\B^1_b,\B^2_b)$, $b_{\theta m}(\B^1_b,\B^2_b)$), where $\B^2_b$ is contained in a disk whose boundary consists of one prong and its homoclinic prong separatrix. 

\bigskip
\noindent
\textbf{($b_{om}$, $b_{im+}$, $b_{im-}$ structures)} The one-dimensional local orbit structure is a prong connection which consists of three $3$-prongs, one $1$-prong, three heteroclinic prong separatrices between them, and two homoclinic separatrices shown in Figure~\ref{h_unstable_prong_connection}.
We assign the COT symbols $b_{om}\{\B_b,\B_b\}$, $b_{im+}(\B^1_b,\B^2_b)$, and $b_{im-}(\B^1_b,\B^2_b)$ as shown in Figure~\ref{fig:list_self_connected_prong}.

\bigskip
\noindent
\textbf{($b_{oi}$, $b_{oo}$, $b_{io}$, $b_{ii}$, $b_{i(o)}$, $b_{i(i)}$ structures)} The one-dimensional local orbit structures are prong connections which consist of three $3$-prongs, one $1$-prong, four heteroclinic prong separatrices between them shown in Figure~\ref{h_unstable_prong_connection}.
When the boundary component of the outer domain contains exactly four prongs (resp. three $3$-prongs), the COT symbol is given by $b_{oo}(\B^1_b,\B^2_b)$ (resp. $b_{io}(\B^1_b,\B^2_b)$), where $\B^2_b$ is contained in a disk whose boundary consists of one prong and its homoclinic prong separatrix. 
When the boundary component of the outer domain contains exactly three prongs one of which is a $1$-prong (resp. two prongs), the COT symbol is given by $b_{io}(\B^1_b,\B^2_b)$ (resp. $b_{ii}(\B^1_b,\B^2_b)$), where $\B^2_b$ is contained in the innermost disk. 
When the boundary component of the outer domain contains exactly one prong and the innermost disk contains (resp. does not contain) a $1$-prong of the prong connection, the COT symbol is given by $b_{i(o)}(\B^1_b,\B^2_b)$ (resp. $b_{i(i)}(\B^1_b,\B^2_b)$), where $\B^2_b$ is contained in the innermost disk.

\bigskip
\noindent
\textbf{($b_{oj}$, $b_{oq}$, $b_{ij}$, $b_{i(q)}$, $b_{i(j)}$ structures)} The one-dimensional local orbit structures are prong connections which consist of four $3$-prongs, four heteroclinic prong separatrices between them, and two homoclinic separatrices shown in Figure~\ref{h_unstable_prong_connection}.
We assign the COT symbols $b_{oq}(\B^0_b, \{ \B_b, \B_b \})$, $b_{oj}(\B^1_b, \B^2_b, \B^3_b)$, where $\B^1_b$ (resp. $\B^3_b$) is contained in an innermost disk whose boundary intersects (resp. does not intersect) the boundary of the outer domain, and $\B^2$ is between the outer domain and the disk containing $\B^3$ as shown in Figure~\ref{fig:list_self_connected_prong}.
Moreover, we assign the COT symbol $b_{ij}(\B_b, \{ \B^1_b, \B^2_b \})$, where $\B^i_b$ for any $i \in \{ 1, 2 \}$ is contained in an innermost disk. 
In addition, we assign the COT symbol $b_{i(q)}(\B_b, \B^1_b, \B^2_b)$, where $\B^1_b$ (resp. $\B^2_b$) is contained in an innermost disk whose boundary contains two prongs (resp. one prong). 
Furthermore, we assign the COT symbol $b_{i(j)}(\B^1_b, \{ \B^2_b, \B^3_b \})$, where $\B^1_b$ is contained in an innermost disk, and where $\B^3_b$ (resp. $\B^2_b$) is (resp. is not) contained in a disk whose boundary intersects the boundary of the outer domain. 

\bigskip
\noindent
\textbf{($b_{qm}$, $b_{jm}$ structures)} The one-dimensional local orbit structures are prong connections which consist of four $3$-prongs, three heteroclinic prong separatrices between them, and three homoclinic separatrices shown in Figure~\ref{h_unstable_prong_connection}.
When the prong connection has (resp. does not have) a rotational symmetry, we assign the COT symbol $b_{qm} \{ \B_b, \B_b, \B_b\}$ (resp. $b_{jm}(\B_b, \B^1_b, \B^2_b)$), where $\B^1_b$ (resp. $\B^2_b$) is contained in an innermost disk such that the outer domain, $\B^1_b$, and $\B^2_b$ are arranged cyclically on the sphere in a counter-clockwise manner. 

\bigskip
\noindent
\textbf{($b_{\theta m}$, $b_{\theta j}$, $b_{\theta q}$, $b_{i\theta}$ structures)} 
The one-dimensional local orbit structures are prong connections which consist of four $3$-prongs, five heteroclinic prong separatrices between them, and one homoclinic separatrix shown in Figure~\ref{h_unstable_prong_connection}.
When the boundary component of the outer domain contains exactly one prong (resp. two, three, four), the COT symbols are $b_{i\theta}$, $b_{\theta m}$, $b_{\theta j}$, and $b_{\theta q}$. 
Moreover, we assign the inner structures as follows: $b_{i\theta}(\B^1_b,\B^2_b,\B^3_b)$ and $b_{\theta q}(\B^1_b,\B^2_b,\B^3_b)$, where $\B^2_b$ (resp. $\B^3_b$) is contained in an innermost disk whose boundary contains exactly three (resp. two) prongs. 
Furthermore, we assign the inner structures as follows: $b_{\theta m}(\B^1_b,\B^2_b,\B^3_b)$ and $b_{\theta j}(\B^1_b,\B^2_b,\B^3_b)$, where $\B^1_b$ (resp. $\B^2_b$) is contained in an innermost disk whose boundary contains exactly one (resp. three) prong. 

\bigskip
\noindent
\textbf{($b_{\theta^3}$, $b_{\theta \theta}$ structure)} 
The one-dimensional local orbit structure is a prong connection that consists of four $3$-prongs, and six heteroclinic prong separatrices which is a complete graph on the four vertices shown in Figure~\ref{h_unstable_prong_connection}.
We assign the COT symbol $b_{\theta^3}(\B^1_b, \{ \B^2_b, \B^3_b \})$, where $\B^1_b$ is contained in an innermost disk whose boundary contains exactly four prongs. 
Moreover, we assign the COT symbol $b_{\theta \theta}( \B^1_b, \B^2_b, \B^3_b)$ as shown in Figure~\ref{fig:list_self_connected_prong}.

\bigskip
\noindent
\textbf{($b_{3}$ structures)} 
The one-dimensional local orbit structures are prong connections which consist of four $3$-prongs, and six heteroclinic prong separatrices between them which is not a complete graph shown in Figure~\ref{h_unstable_prong_connection}.
We assign the COT symbol $b_{3}\{ \B_b, \B_b, \B_b\}$.


\subsection{COT representations of at most coheight-one line fields on the sphere}

We have the nest structures for at most coheight-one structures as shown in Table~\ref{tbl:All_structures} and Table~\ref{tbl:All_structures_codim_1}. 

As above, for any line field $\F$ of $\mathcal{L}^1(\mathbb{S}^2)$ with a root component, we  can assign the COT representation, which is a sequence of symbols, of $\F$ as follows: 
\\
(1) If the root component is of type $\sigma$ (resp. a periodic leaf on the boundary), then the intial structure of the sequence of symbols is $\sigma_{\emptyset}(\B_{b})$ (resp. $\beta_{\emptyset}(\B_{b})$). 
\\
(2) Since the Reeb graph has a root given by the vertex of degree one which corresponds to the root component, tracing adjacent vertices from the root uniquely determines the structure of the associated prong connections or boundary components. 
Therefore, for every vertex, we recursively substitute either $\sigma$, $\beta$, $p_o(\B_{b})$, $p_i(\B_{b})$, $b_\theta(\B_{b})$, $b_o(\B_{b})$, $b_i(\B_{b})$, or one of prong connections in Figure~\ref{fig:list_self_connected_prong}, into $\B_b$ according to the corresponding prong connection. 
The substitution proceeds in order from the vertices closest to the root, where possible.

The resulting sequence of symbols is called the partially cyclically ordered rooted tree ({\bf COT}) {\bf representation} of the line field $\F$ of $\mathcal{L}^1(\mathbb{S}^2)$ with the root component. 

By construction of the COT representation, Theorem~\ref{th:tree_rep_codim_one} implies Theorem~\ref{th:cot_02}. 

%



\section{Examples of singular foliations on surfaces}

We have the following examples of singular foliations on surfaces. 

\begin{example}\label{ex:liftable}
{\rm
Let $\mathbb{K} := ([0,1]\times \mathbb{S}^1)/\mathop{\sim}$ be a Klein bottle where $\mathbb{S}^1 = \R/\Z$ and $(x,y) \sim (x',y')$ if $\{x, x' \} = \{ 0,1\}$ and $1-y \equiv y' \pmod 1$. 
The smooth foliation $\F_{\mathbb{K}} := \{ \{ [t]\} \times \mathbb{S}^1 \mid t \in [0,1) \}$ is a non-orientable foliation. 
Then we have the following observation. 

\begin{lemma}
The foliation $\F_{\mathbb{K}}$  is not liftable to a flow. 
\end{lemma}
\begin{proof}
Assume that $\F_{\mathbb{K}}$ is liftable to a flow. 
Then there is a foliation $(M, f, \theta)$ in the sense of Bronstein and Nikolaev. 
Since $\F_{\mathbb{K}}$ has no singularities, by definition of foliation in the sense of Bronstein and Nikolaev, the involution $\theta$ is a homeomorphism between Klein bottles. 
This implies that the set of orbits of the restriction $f$ to a Klein bottle is $\F_{\mathbb{K}}$. 
Therefore, $\F_{\mathbb{K}}$ is orientable, which contradicts that it is non-orientable. 
\end{proof}
}
\end{example}

\begin{example}\label{ex:liftable02}
{\rm
Let $\mathcal{G}$ be a smooth foliation on a torus $\mathbb{T} := (\R/\Z)^2$ defined by the set of orbits of a Hamiltonian vector field as shown on the left in Figure~\ref{disk_Hamiltonian}. 
\begin{figure}[t]
\begin{center}
\includegraphics[scale=0.375]{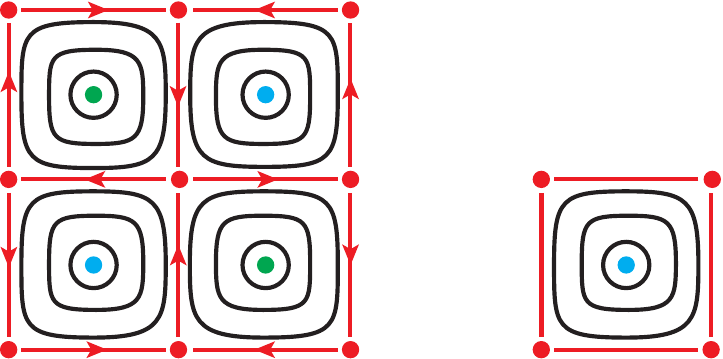}
\end{center} 
\caption{Left, a smooth foliation $\mathcal{G}$ on a torus defined by the set of orbits of a Hamiltonian vector field; right, a non-orientable foliation $\F_{\T^2}$ whose lift along a covering of degree $4$ is the foliation $\mathcal{G}$. }
\label{disk_Hamiltonian}
\end{figure}
Taking a $(\Z_2)^2$-action on $\mathbb{T}^2$ generated by $1/2$-rotations with respect to the $x$-axis and the $y$-axis, the quotient space $\left( \R/\left(\frac{1}{2}\Z \right) \right)^2$, which is also a torus, has the induced non-orientable foliation $\F_{\T^2}$ as shown on the right in Figure~\ref{disk_Hamiltonian}.  
Then we have the following observation. 
\begin{lemma}
The foliation $\F_{\T^2}$ is not liftable to a flow. 
\end{lemma}

\begin{proof}
The non-orientable foliation $\F_{\T^2}$ has exactly two singular points, which are a center and a saddle (i.e. $4$-prong). 
Therefore, $\F_{\T^2}$ has no non-orientable singularities. 
Assume that $\F_{\T^2}$ is liftable to a flow. 
Then, there is a foliation $(M, f, \theta)$ in the sense of Bronstein and Nikolaev. 
Since $\F_{\T^2}$ has no singularities, by definition of foliation in the sense of Bronstein and Nikolaev, the involution $\theta$ is a homeomorphism between tori. 
This implies that the set of orbits of the restriction $f$ to a torus is $\F_{\T^2}$. 
Therefore, $\F_{\T^2}$ is orientable, which contradicts that it is non-orientable. 
\end{proof}
}
\end{example}

\begin{example}
{\rm
Let $\F_{\theta}$ be a regular foliation on a torus $\T^2 := \R^2/\Z^2$ which is the set of orbits of a vector field $X = (1,\theta)$ on $\T^2$. 
Equip $\T^2$ with the induced Riemannian metric from the standard Riemannian metric on $\R^2$. 
For any convergence sequence $(\theta_n)$ of real numbers to $\theta_\infty \in \R$, the sequence $(\F_{\theta_n})$ of regular foliations converges to the foliation $\F_{\theta_\infty}$. 
This implies that the regular foliation $\F_{\theta}$ on a torus $\T^2$ for any $\theta \in \R$ is not structurally stable. 
}
\end{example}

\begin{example}\label{ex:nonHam}
{\rm
Let $\mathbb{K} := ([0,1]\times \mathbb{S}^1)/\mathop{\sim}$ be a Klein bottle where $\mathbb{S}^1 = \R/\Z$ and $(x,y) \sim (x',y')$ if $\{x, x' \} = \{ 0,1\}$ and $1-y \equiv y' \pmod 1$. 
The vector field $X = (1,0)$ on $\mathbb{K}$ is not a Hamiltonian vector field, but the set of orbits is a levelable foliation that is not transversely orientable. 
}
\end{example}

\begin{example}\label{ex:rot_disk}
{\rm
Let $p \colon \R/\Z \times [-1,1] \to D := \{(x,y) \in \R^2 \mid x^2 + y^2 \leq 1\}$ be a branched covering with two ramification points of degree $2$ such that $p(t, \pm 1) = (\cos 2\pi t, \sin 2\pi t)$ as in Figure~\ref{disk_2prongs02}, and $\F_D$ a foliation on $D$ as in Figure~\ref{disk_2prongs02}. 
\begin{figure}[t]
\begin{center}
\includegraphics[scale=0.375]{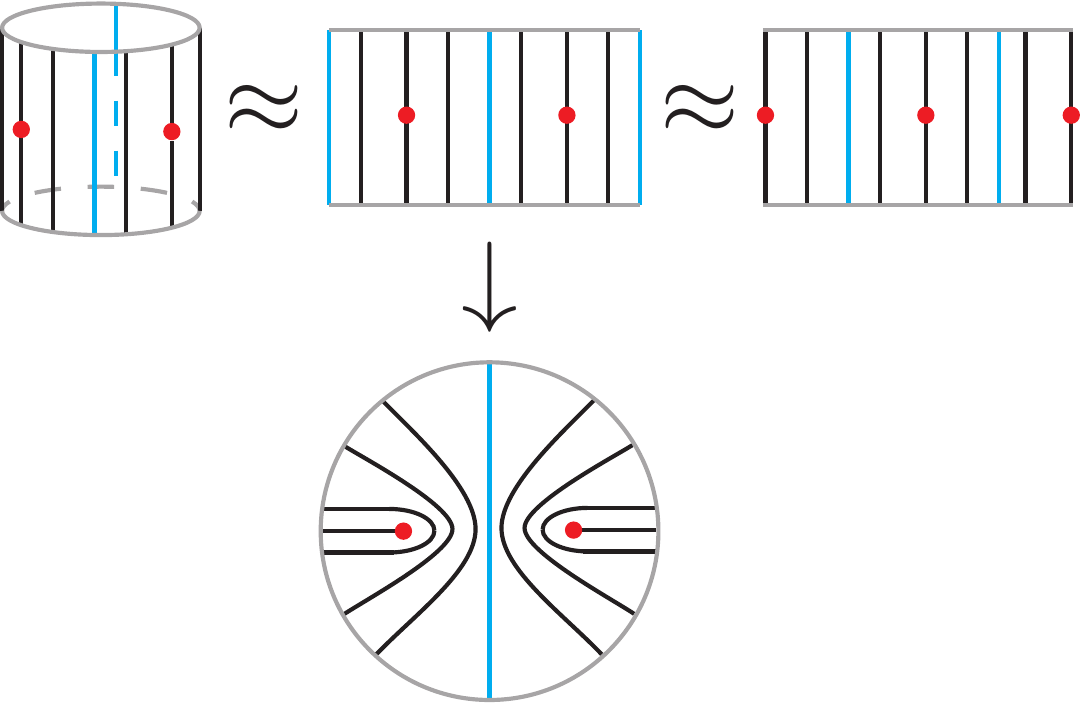}
\end{center} 
\caption{A closed disk with a foliation $\F_D$ whose branched covering $p \colon \R/\Z \times [-1,1] \to D$ with two ramification points of degree 2 lifts a closed annulus $\R/\Z \times [-1,1]$ with a foliation with two singular points, which are fake multi-saddles. }
\label{disk_2prongs02}
\end{figure}
Identifying $D$ with the hemisphere of the unit sphere $\mathbb{S}^2$, equip $D$ with the induced Riemannian metric from the standard Riemannian metric on $\mathbb{S}^2$. 
Pasting two copies of the unit disk $D$ equipped with the foliation $\F_D$ with $\theta$ rotation, the resulting surface is a sphere, and the resulting singular foliation is denoted by $\F(\theta)$. 
Then the singular foliations $\F(\theta)$ can be realized as divergence-free singular foliations because the singular foliations can be obtained by the resulting singular foliations by the quotient of singular foliations with four singular points $(0,0), (0,1/2), (1/2,0)$, and $(1/2,1/2)$, which are fake saddles, on a torus by the rotation $f \colon (\R/\Z)^2 \to (\R/\Z)^2$ by $(x,y) \mapsto (-x,-y)$ as in Figure~\ref{disk_3prongs}. 
\begin{figure}[t]
\begin{center}
\includegraphics[scale=0.375]{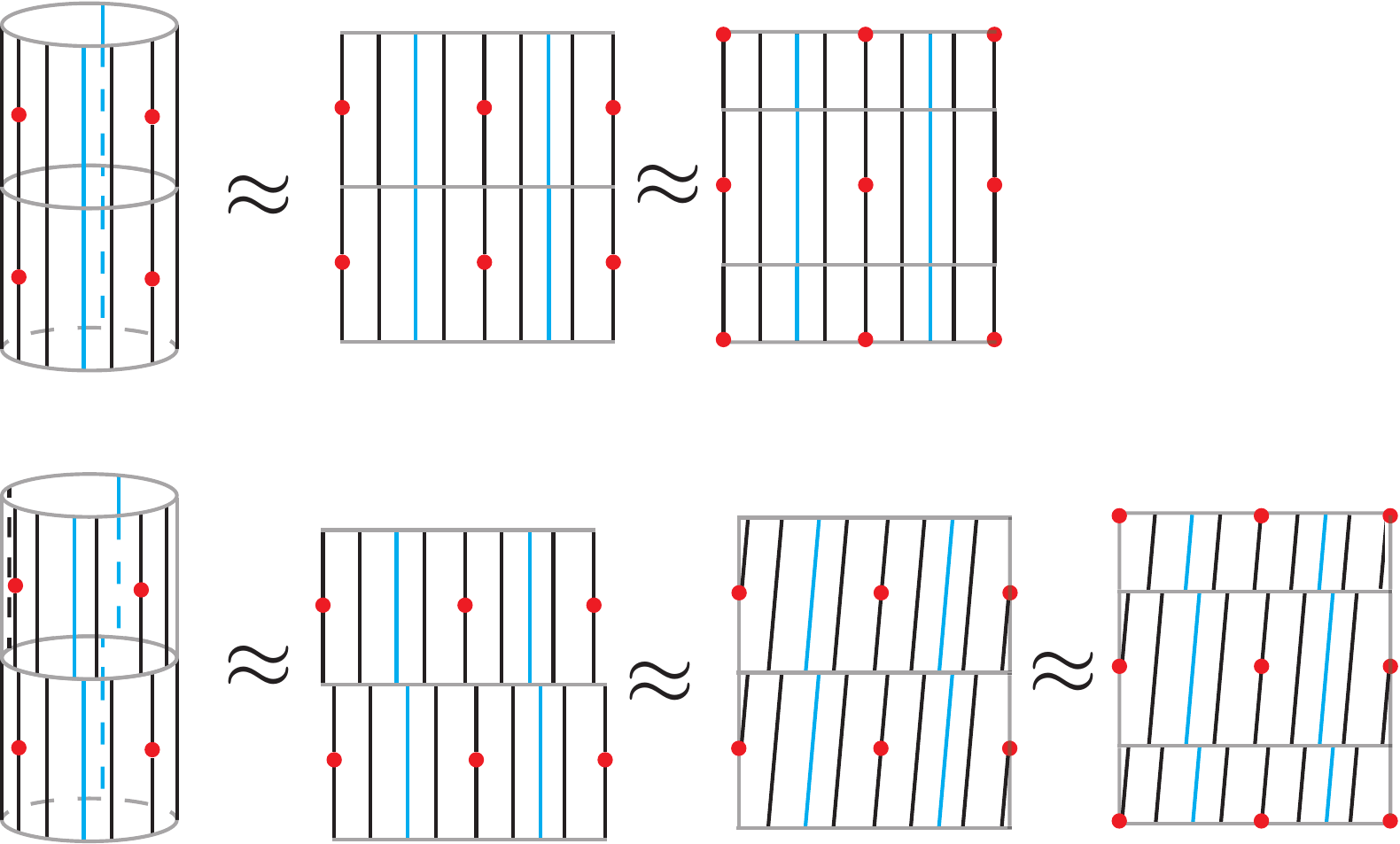}
\end{center} 
\caption{Spheres with foliations whose branched coverings $p \colon (\R/\Z)^2 \to {\mathbb{S}^2}$ with four ramification points of degree 2 lift tori with divergence-free foliations with four singular points  $(0,0), (0,1/2), (1/2,0)$, and $(1/2,1/2)$, which are fake multi-saddles.}
\label{disk_3prongs}
\end{figure}
For any convergence sequence $(\theta_n)$ of real numbers to $\theta_\infty \in \R$, the sequence $(\F(\theta_n))$ of regular foliations converges to the foliation $\F(\theta_\infty)$. 
This implies that the divergence-free regular foliation $\F(\theta)$ on the sphere $\mathbb{S}^2$ for any $\theta \in \R$ is not structurally stable.

}
\end{example}

\section{Final remarks}

Finally, we remark on Voronoi-like behaviors. 
We define Voronoi-like foliations to describe perturbations of Voronoi partitions. 

\begin{figure}[t]
\begin{center}
\includegraphics[scale=0.375]{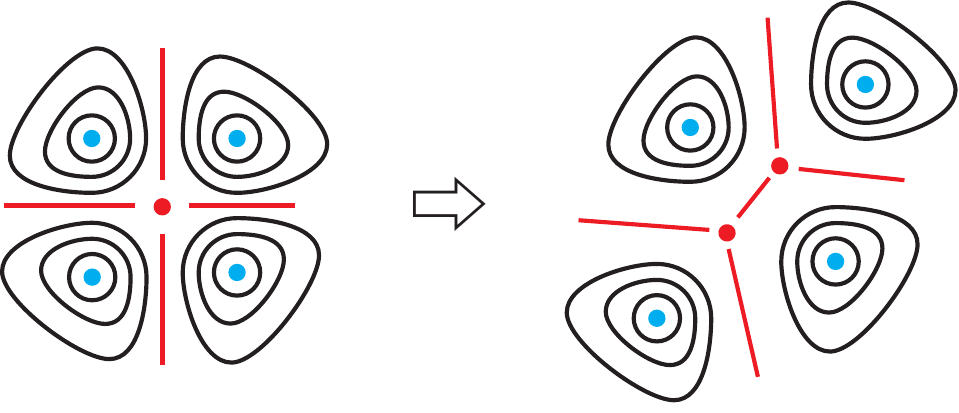}
\end{center} 
\caption{A perturbation by a symmetric breaking which implies a Whitehead move among Voronoi-like foliations.}
\label{Voronoi_prong}
\end{figure}

\begin{definition}
A levelable foliation $\F$ on a compact surface $M$ is {\bf Voronoi-like} if there is a finite set $V$ such that any connected components of the complement $M - (V \cup G_V)$ are periodic annuli, where $G_V := \bigcup_{p \in V} \partial \{ x \in M \mid d(x, p) \leq d(x, V - \{p \})\}$ is the Voronoi diagram of $V$, where $d$ is a Riemannian distance. 
\end{definition}
Notice that Whitehead moves appear in perturbations among Voronoi-like foliations as in Figure~\ref{Voronoi_prong}.

Moreover, to analyze various singular foliations on non-compact surfaces in nature, we will also discuss foliations on non-compact surfaces. 
Therefore, we will report various properties (e.g. openness, density, structural stability) of Voronoi-like foliations on surfaces and 
line fields on non-compact surfaces in future works.



\vspace{10pt}

{\bf Acknowledgement}: 
The author would like to thank Masashi Asaoka for his helpful comments. 

\appendix

\section{Flows and vector fields}\label{sec:flow}

\subsection{Flows}

Recall that a {\bf flow} on a manifold is a continuous $\R$-action on it.
Let $f \colon \R \times M \to M$ be a flow on a manifold $M$.
For $t \in \R$, define $f^t \colon X \to X$ by $f^t := f(t, \cdot )$.
For a point $x$ of $X$, we denote by $O(x)$ the orbit of $x$, $O^+(x)$ the positive orbit (i.e. $O^+(x) := \{ f^t(x) \mid t > 0 \}$), and $O^-(x)$ the negative orbit (i.e. $O^-(x) := \{ f^t(x) \mid t < 0 \}$).
A point $x$ of $X$ is {\bf singular} if $x = f^t(x)$ for any $t \in \R$, and is {\bf periodic} if there is a positive number $T > 0$ such that $x = f^T(x)$ and $x \neq f^t(x)$ for any $t \in (0, T)$.

\subsubsection{Non-wandering flows}

A point of a flow $f$ is {\bf wandering} if there are its neighborhood $U$ and a positive number $N$ such that $f^t(U) \cap U = \emptyset$ for any $t > N$. 
Such $U$ is called a {\bf wandering domain}. 

\begin{definition}\label{def:nw}
A point is {\bf non-wandering} if it is not wandering (i.e. for any its neighborhood $U$ and for any positive number $N$, there is a number $t \in \mathbb{R}$ with $|t| > N$ such that $f^t(U) \cap U \neq \emptyset$).
\end{definition}

Denote by $\Omega (f)$ the set of non-wandering points, called the non-wandering set. 
A flow is {\bf non-wandering} if the non-wandering set is the whole space. 

\subsection{Vector fields on surfaces}

\subsubsection{Divergence-free vector fields and Hamiltonian vector fields}
Recall the following concepts. 

\begin{definition}A $C^r$ vector field $X$ for any $r \in \Z_{>0}$ on a surface with a Riemannian metric $g$ is {\bf divergence-free} if $\mathrm{div} X = 0$, where $\mathrm{div}X := \mathop{*} d \mathop{*} g(X, \cdot )$. 
\end{definition}

\begin{definition}\label{def:df}
A flow {\bf divergence-free} if it is topologically equivalent to a flow generated by a smooth divergence-free vector field. 
\end{definition}

\begin{definition}\label{def:Ham_vf}
A $C^r$ vector field $X$ for any $r \in \Z_{\geq0}$ on an orientable surface $S$ is {\bf Hamiltonian} if there is a $C^{r+1}$ function $H \colon S \to \mathbb{R}$, called a {\bf Hamiltonian},  such that $dH= \omega(X, \cdot )$ as a one-form, where $\omega$ is a volume form of $S$.
\end{definition}

In other words, locally, the Hamiltonian vector field $X$ is defined by $X = (\partial H/ \partial x_2, - \partial H/ \partial x_1)$ for any local coordinate system $(x_1,x_2)$ of a point $p \in S$.

\begin{definition}
A flow is {\bf Hamiltonian} if it is topologically equivalent to a flow generated by a smooth 
\end{definition}

Hamiltonian vector field.
Note that a volume form on an orientable surface is a symplectic form.

\section{Transversality for singular piecewise $C^1$-foliations}\label{def:trans}

Let $\F$ be a singular piecewise $C^1$-foliation on a surface $S$. 
We define transversality as follows.  

\begin{definition}
A curve $C$ is {\bf transverse} to $\F$ at a point $p \in C \setminus \partial S$ if there are a small neighborhood $U$ of $p$ and a homeomorphism $h:U \to [-1,1]^2$ with $h(p) = 0$ such that $h^{-1}([-1,1] \times \{t \})$ for any $t \in [-1, 1]$ is contained in a leaf of $\F$ and $h^{-1}(\{0\} \times [-1,1]) = C \cap U$. 
\end{definition}

When $S$ has the boundary, denote by $\hat{S}$ the double of $S$ and by $\hat{\mathcal F}$ the induced foliation on $\hat{S}$ from $\F$.
Recall that $\hat{S} = S \times \{0,1\}/\mathop{\sim}$ where $(x,0)\sim (x,1)$ if $x\in \partial M$. 

\begin{definition}
A curve $C$ is {\bf transverse} to $\F$ at a point $p \in C \cap \partial S$ if there is a curve $\widetilde{C} \subset \hat{S}$ transverse to $\hat{\F}$ such that the restriction $\pi \vert_{\widetilde{C}} \colon \widetilde{C} \to C$ is homeomorphic, where $\pi \colon \hat{S}= S \times \{0,1\}/\mathop{\sim} \to S$ is the canonical quotient map define by $p(x,0) = p(x,1) = x$.  
\end{definition}

A simple curve $C$ is transverse to $\F$ if it is transverse at any point in $C$.  
A simple curve $C$ which is transverse to $\F$ is called a {\bf transverse arc}.
A simple closed curve is a {\bf closed transversal} if it transverses $\F$.

\bibliographystyle{abbrv}
\bibliography{yt20211011}

\end{document}